\documentclass[10pt, oneside]{article}   	
\usepackage{geometry}                		
\usepackage{graphicx}				
				
\usepackage{amssymb}

\title{Pavel Florensky and his world}
\author{Athanase Papadopoulos\footnote{\emph{Institut de Recherche Mathématique Avancée} and
\emph{Centre de Recherche et d'Expérimentation sur l'Acte Artistique  (ITI CREAA)}. Address:
Université de Strasbourg and CNRS,
7 rue René Descartes,
67084 Strasbourg Cedex France;
email: athanase.papadopoulos@math.unistra.fr}}
 \date{\today}							

\begin{document}
\maketitle

%

%

\hfill \emph{They call it ``Silver" Age, but it is comparable}

\hfill \emph{to the Renaissance in Europe}

\hfill  (Valentin Poénaru told me).

 \bigskip
  \bigskip
 
\begin{abstract}

This is an overview of the life and works of Pavel Florensky, an important and singular figure of the period rightly described as the
 \emph{Silver Age of Russian mathematics}, with a substantial overlap with the \emph{Silver Age of Russian literature, poetry and philosophy}. Florensky is certainly among the great scientists, philosophers, theologians and historians of art of the twentieth century, with a very atypical trajectory of life. He had a mathematical background and his work in philosophy, theology and history of art is imbued with mathematical ideas. Talking about his life  and works is also an opportunity to reflect upon the Russian mathematical school of the first third of the twentieth century, its  philosophical foundations and the conflicts it underwent. It is also an occasion for discussing poetry, literature and art during the Russian Silver Age.  
 This article will appear as a chapter in the Handbook of the History and Philosophy of Mathematical Practice, edited by Bharath Sriraman, published by Springer,  to appear in 2023.
 
\end{abstract}

\bigskip


\bigskip


\bigskip

 \section{Introduction}
 
Pavel Aleksandrovich Florensky (1882-1937)\index{Florensky, Pavel Aleksandrovich} is a complex and singular figure in the landscape of twentieth century culture. His complete works contain more than 1000 items (see \cite[Note p. 7]{Stupeur}).  In the last three decades, his books and manuscripts are continuously edited and published in several languages and the literature on him has been growing at an exponential rate.  Symposia on his case are regularly organized in Europe and elsewhere in the world, and a number of doctoral dissertations on his work have been written.  Biographers and people who knew him dubbed him the Russian Leonardo,\footnote{In \cite[p. 178]{Lossky-History}, the philosopher Nikolai Onufriyevich Lossky\index{Lossky, Nikolai Onufriyevich} (1870-1965) quotes an article by B. Filistinsky  titled \emph{A Russian Leonardo da Vinci in a concentration camp}, in which the latter says about Florensky:  ``A new Leonardo da Vinci was standing before us and we all were conscious of it".} others called him the Russian Pascal,\footnote{See Vassili Rozanov  \cite[p. 413]{Rozanov}.} sometimes the Gulag Pascal. Cédric Villani,\index{Villani, Cédric} in the preface to his book \cite{Imaginaires}, describes him as ``one of the most original thinkers of his times".   Florensky\index{Florensky, Pavel Aleksandrovich} is invariably presented as a mathematician, physicist, chemist, botanist, theologian, engineer, philosopher, art historian, specialist in electrodynamics, marine botany, linguistics,  languages, symbol, esoterism, gnoseology, logic, aesthetics, semiotics and  other fields of knowledge. Strictly speaking, I am not sure he may be called a mathematician, since I don't know of any mathematical discovery he made,  and I cannot really judge about most of the other fields I mentioned. But without any doubt his writings and thought encompass all these fields, and mathematics is always present there as a background. His highly unusual versatility led Nikolai Lossky to describe him as possessing ``inhuman erudition" \cite{Lossky1926}.  Perhaps this is where his strength and originality lie.

 When I started writing this essay, my intention was to concentrate on Florensky's\index{Florensky, Pavel Aleksandrovich} mathematical inclination, since the essay is addressed primarily  to mathematicians. But I soon realized that 
 it is impossible to isolate the mathematical dimension of his thought 
from the rest of his works and vision of life.  Florensky certainly had a mathematician's liking in his view of the world; it is this view, and some reflexions that go with it, that I would like to outline.

 Speaking of Florensky, I cannot avoid to begin by saying a few words about his life, and that is what I shall do in the next section.

\section{Vita} \label{s:Vita}

Pavel Aleksandrovich Florensky\index{Florensky, Pavel Aleksandrovich} was born on January 22 (Old Style: January 9), 1882, deep in the steppe, near the town of Yevlakh, today in Azerbaijan, to an Armenian mother,  Salomé Saparian,\index{Saparian, Salomé}\index{Saparian, Salomé} from an aristocratic family and a Russian father, Aleksander Ivanovich Florensky,\index{Florensky, Aleksander Ivanovich}  who was working as an engineer in the construction of the Caucasus railroads. Azerbaijan was part of the Russian Empire, and Aleksander Florensky was in charge of the railway section that passed through Yevlakh. Pavel was born before the city's railroad station was built, and the family lived in freight cars ``dressed in oriental carpets", he recalls.\footnote{In this section, unless otherwise specified, the sentences in quotation marks are taken from Florensky's autobiographical notes addressed to his children, which he wrote between 1916 and 1926, \cite{Florensky-bio}. Here and elsewhere, I have translated from the French (I cannot do it from the Russian).} Thus, it is possible that it is in a train carriage that Pavel was born.
 His childhood was spent in Tbilissi, then in Batumi,\footnote{Batumi, situated on the Black Sea, experienced its expansion with the construction of the Batumi-Tiflis-Baku Transcaucasus Railway, completed soon after Pavel Florensky's birth. Incidentally, Batumi is a town where Joseph Stalin lived, organizing strikes since 1901.} another Georgian town, under the guidance of a father
who had a ``distaste for convention". Family life, Pavel says, ``was like being on an island cut off from the world, that is, a desert island, because we didn't particularly like people and we always tried to keep ourselves out of the way. People would have taken hold of the purity, the unchanging calm and the rigorism of this insular paradise, which is why they were barely tolerated, and only occasionally." During these childhood years, Pavel was fascinated by the abundant nature of the Caucasus, of which he would later speak on several occasions. Especially, he recalled this landscape with particular emotion many years later, while he was a deportee in the Solovki Archipelago\index{Solovki Islands} of the White Sea, the Gulag where he spent the last four years of his earthly existence. It is conceivable that it was because of the powerful influence of nature, especially that of the Caucasus, which fascinated him as a child, that Florensky\index{Florensky, Pavel Aleksandrovich} remained deeply attracted, throughout his life, by the study of natural sciences and all the phenomena of the living.

At school, Florensky was a gifted pupil, interested in all subjects, even if attending school every day bored  him. Those who knew him at that age described him later on as a solitary child, 
extremely shy and very much an insomniac.  In his recollections, we read: 
``I would fight long and hard not to go to bed, and once in bed, I would go on for hours without falling asleep, turning from side to side and studying the pattern of the wallpaper and blanket for a million times. These were hours of almost torture, when I was in bed unable to sleep." 
Pavel reads a lot, sets up a chemistry laboratory in his room and works alone at night. He is interested in ancient languages. Later, he was able to speak most European languages, besides being proficient in Latin and ancient Greek, of course, but also in Hebrew and Syriac, which he learned at the Theological School.  Several years later, from his exile in the Gulag, he wrote to his 15-year-old son Mikhaïl \cite[Letter no. 77, p. 497, dated 23.X.1936]{Florensky-Letter}:   
\begin{quote}\small
I wish I could have been with you to teach you how to work properly and accumulate knowledge. [\ldots] When I was your age, every minute lost seemed to me like a misfortune or a crime, and I tried to fill my time as much as possible. I had diaries where I tried to record the essentials of my reading and reactions to books, notebooks of interesting quotations, albums of sketches from nature, copybooks where I methodically recorded the results of experiments, and others for on-site observations. Every evening, I would write myself a (motivated) note for the day's work. In this way, I acquired a wealth of knowledge, practical skills and, above all, the ability of making my own judgements based on what I had observed. I tried to compare and synthesize the data I collected in the form of tables, diagrams and curves, as this condensed aspect made them easier to understand, they become alive and acquire meaning. Immediately and automatically, you get an `empirical generalization'.
\end{quote}

In 1899, Florensky,\index{Florensky, Pavel Aleksandrovich} who was seventeen, arrived to Moscow, where he was to enter university. At his arrival, he  underwent a profound existential crisis. Having been brought up away from religion, in a positivist climate and an atmosphere of secular humanism, 
he realized the limitation of the natural sciences and became interested  in social,  religious and moral issues, in the life of the ordinary and poor people, in the human soul, striving for a meaning to life.  The crisis is recounted in the last part, titled \emph{Avalanche}, of the auto-biographical text he wrote for his children \cite{Florensky-bio}.
We read there:  ``On top of my own magical idealistic understanding of the world, a film had formed that prevented me from breathing freely; the idea of continuity,\index{continuity (philosophical notion)} fundamentally excluding the miraculous, was trying to take possession of my intelligence."  He was already seeing the world with the eye of a mature mathematician. The ideas of continuity and discontinuity,\index{discontinuity (philosophical notion)} which for him had become philosophical concepts, were to follow him for the rest of his life.
Florensky's\index{Florensky, Pavel Aleksandrovich} existential and religious crisis was often compared to that of Pascal,\index{Pascal, Blaise} who, at age twenty-nine, experienced  a similar experience which  led him to a life centered on philosophical and theological reflection, without leaving science. Pascal had, like Florensky, a strange personality. T. S. Eliot described him as  ``a man of the world among ascetics, and an ascetic among men of the world" \cite[p. 411]{Eliot}.

Under pressure from his parents, Florensky enrolled as a mathematics student at the Faculty of Physics and Mathematics of the Imperial University of Moscow. His mentor there was Nikolai Vasilievich Bugaev (1837-1903),\index{Nikolai Vasilievich Bugaev} a well-known mathematician who was the president of the Moscow Mathematical Society, of which he had been a co-founder. Bugaev  was engaged in the philosophy of mathematics and its repercussions in the social, historical and religious spheres. For him, questions like continuity,\index{continuity (philosophical notion)} discontinuity,\index{discontinuity (philosophical notion)} the finite and the infinite,\index{infinity (philosophical notion)} were intimately linked to philosophical issues such as free will, scientific and social revolutions, evolutions in history, etc.
I shall talk later in this chapter about Bugaev,\index{Nikolai Vasilievich Bugaev} whose ideas were to have a definite influence on Florensky.\index{Florensky, Pavel Aleksandrovich}  The other teachers of Florensky included Boleslav Kornelievich Mlodzievsky\index{Mlodzievsky, Boleslav Kornelievich} (1858-1923), a representative of the Moscow school of metric geometry  who, like Bugaev and several other Russian mathematicians  of that period,  spent some time in Western Europe on several occasions, in cities like Göttingen, Paris and Zürich, where he came into direct contact with Western mathematicians and philosophers. According to I. H. Anellis in \cite{Anellis}, Mlodzievsky gave the first course on set theory\index{set theory (historical development)}  and the theory of functions at Moscow University, two topics which became important parts of Florensky's background. Anellis adds that it was Mlodzievsky who established the Russian terminology in these subjects. According to Demidov \cite{Demidov2014}, Mlodzievsky's synopsis of his lectures was preserved thanks to notes compiled by Florensky \cite{Medvedev-Notes}. 

Florensky's\index{Florensky, Pavel Aleksandrovich} teachers also included Leonid Kuz'mich Lakhtin\index{Lakhtin, Leonid Kuz'mich} (1863-1927), a student of Bugaev, who studied  non-orientable surfaces,\index{non-orientable surface} in particular the question of representing them by explicit equations, see \cite{Lakhtin}. We  shall see below how non-orientable surfaces became part of Florensky's later interests, mathematically and symbolically. In a note in his booklet titled  \emph{The imaginaries in Geometry}, which we mentioned in the introduction, which we shall discuss at length below and in which the notion of non-orientability plays a key role, Florensky says that he wrote that book during the month of August 1902 while he was a student, and that he showed it to Lakhtin, see \cite[p. 79]{Imaginaires}. 
 Among the Russian mathematicians of the same epoch who studied non-orientable surfaces,\index{non-orientable surface} one should mention Nikolai Borisovich Delaunay\index{Delaunay, Nikolai Borisovich} (1856–1931) who addressed the problem of finding a locally Euclidean parametrization of the M\"obius band, in his memoir \cite{Delaunay} published in 1898 in the Bulletin of the French Mathematical Society, and which Florensky also mentions in a note to his work on the \emph{Imaginaries in geometry},  see \cite[p. 82]{Imaginaires}. In the same note, Florensky declares that in the \emph{Encyclopedia of mathematics} of which he was in charge, he wrote an article on the division of surfaces in orientable and non-orientable.

Below, I shall discuss the booklet titled \emph{Imaginaries in Geometry}, which was essentially written in 1902, at the time Florensky\index{Florensky, Pavel Aleksandrovich} was a mathematics student.

The young Florensky became quickly interested in formal logic. Questions of contradiction,\index{contradiction (philosophical notion)} antinomy\index{antinomy (philosophical notion)} and the nature of the relationship between mathematics and the physical world found an echo in him, both philosophically and theologically. 
He writes\index{Florensky, Pavel Aleksandrovich} \cite[p. 42]{Florensky-bio}: 
\begin{quote}\small
In mathematics, Fourier series and other topics which represent any complex rhythm as an infinite set of simple objects resonate with me, speak to me, correspond to my innermost feelings. Continuous\index{continuity (philosophical notion)} functions without derivatives and discontinuous\index{discontinuity (philosophical notion)} functions where everything unravels, where all elements are \emph{standing}, hold a dear language for me. If I listen to the depths of my being, I discover in the rhythm of my inner life, in the sounds that fill my consciousness, the rhythm of the waves that have been retained forever, inscribed in my memory for all time, and I know that they seek their conscious expression in me through the schema of these mathematical concepts.
\end{quote}

Leonid Sabaneeff,\footnote{Leonid Sabaneeff\index{Sabaneeff, Leonid Leonidovich} (1881--1968) graduated in mathematics from the Imperial University of Moscow and studied music under Rimsky-Korsakov and Taneyev. He became a composer and musicologist. He wrote an authoritative and influential book on Skriabin. He taught at the State University of Moscow and he founded the Moscow Institute of Musicology. He remained close to Florensky and to Nikolai Nikolaevich Luzin\index{Luzin, Nikolai Nikolaevich}  (1883-1950), whom I discuss below,\index{Florensky, Pavel Aleksandrovich} until the year 1926, when he left the Soviet Union.  He lived in Paris, London, the United States and Nice, where he is buried.} who knew Florensky since they were both mathematics undergraduates at the same university, writes in a set of recollections \cite{Sabaneeff} that Florensky rarely attended lectures,  preferring to
read books at home. He also recalls that Florensky seemed oddly old for his age, to the point that one of their fellow students said about him: ``he looks as
if he had already lived a thousand lives." Sabaneeff continues: ``One
 had an impression of genius, of unusual
depth and power of thought, and at the same time there was
something of black magic, dark, devoid of divine grace."  He had the feeling that Florensky possessed an immense spiritual experience
and was endowed with hypnotic power. He recalls Florensky\index{Florensky, Pavel Aleksandrovich} speaking of ``many facets of any true thought and of the compatibility of contradictions at the deepest levels", asserting  that ``every perceived law generates its own negation" and that  ``logic is valid for the earthly life, for the lower levels; but the true world is one where contradictions are compatible---a realm of antinomies."  Sabaneeff adds that Florensky ``obviously regarded antinomy\index{antinomy (philosophical notion)}  as the basic law of the universe, encompassing all others." In Florensky's mind, the notions of contradiction and antinomy have no negative connotation. Contradiction,\index{contradiction (philosophical notion)} non-provability and paradoxes in set theory fascinated him, and he often referred to them in his mathematical and non-mathematical writings.

At university, Florensky, besides studying mathematics, attended lectures at the philosophy department, by Lev Mikhailovich Lopatin\index{Lopatin, Lev Mikhailovich} (1855-1920), a psychologist  and philosopher  inspired by Leibniz' monadology, and by Sergei Nikolaevich Trubetskoy\index{Trubetskoy, Sergei Nikolaevich} (1862-1905),  who belonged to the Russian tradition of universal theism, trying to unite christianity with Platonism.\index{Platonism} Both philosophers were imbued with spiritual ideas, in search of universal\index{Florensky, Pavel Aleksandrovich!\emph{The pillar and ground of the truth}} wisdom,\footnote{Florensky, in the notes of his major theological work, \emph{The  pillar and ground of the truth} \cite{Flo1}, quotes several works by Lopatin, among which \emph{The positive tasks of philosophy}, \emph{The philosophical worldview of N. V. Bugaev} and \emph{The present and future of philosophy}. He also quotes books by Trubetskoy, among which \emph{The metaphysics of ancient Greece}, \emph{Foundations of idealism} and \emph{The doctrine of Logos}.} and their ideas had a definite impact on Florensky.\index{Florensky, Pavel Aleksandrovich}

  Bugaev\index{Nikolai Vasilievich Bugaev} died in 1903. Dmitri Fedorovich Egorov\index{Egorov, Dmitri Fedorovich} (1869-1931), who had been a student of Bugaev,  took over the direction of the mathematical school founded by his mentor.
At that time,
Florensky had already become friends with Bugaev's son,\index{Nikolai Vasilievich Bugaev} Boris Nikolaevich Bugaev (1880-1934), two years his senior, poet and writer, known by the pseudonym Andrei Belyi,\index{Belyi, Andrei} who was a leader and theorist of the Russian Symbolist movement.  Belyi's novel \emph{Petersburg } (first edition completed in 1913) made him later one of the greatest writers of the twentieth century.  Like Florensky,\index{Florensky, Pavel Aleksandrovich}  Belyi was interested in the impact of mathematics in life sciences and philosophy. During his study years at the Mathematics and Physics Department of the University of Moscow, Florensky initiated a ``Circle of mathematics students" within the Russian Mathematical Society, with the collaboration of his friend and fellow student Nikolai Luzin\index{Luzin, Nikolai Nikolaevich}, who was to become later one of the best-known Russian mathematicians of his generation, respected in Russia and France. On this occasion, he wrote a  \emph{Sketch for a discourse for the inauguration of the Circle of mathematics students of the Moscow Mathematical Society} which was published in 1990 \cite{Florensky-sketch}, discussing there, among other topics, Cantor's set theory\index{set theory (historical development)}  and Bugeav's arithmology (the theory of discontinuous\index{arithmology} functions). 
At about the same time, Florensky founded a ``Religious-Philosophical Society of Writers and Symbolists" and a ``Union of Christian Struggle", a movement dedicated to social Christianity, opposing  the then current subjugation of the Church to the State, which had started under Peter the Great. 

 I will say more about Luzin in \S \ref{s:Luzin}.

In 1904, Florensky\index{Florensky, Pavel Aleksandrovich} completed his mathematics studies with a dissertation titled \emph{On singularities of plane curves as loci violating their continuity}. Part of this dissertation is published in \cite{DP1986}.  According to Oppo \cite[p. 72 ]{Oppo},  this dissertation was part of a
wider doctoral project titled \emph{The Idea of discontinuity as part of a postulated worldview}. Florensky wrote an introduction to this work and published it under the title \emph{On a certain premise of a worldview} \cite{Florensky-certain}   in the symbolic journal \emph{Vesy}. Oppo gives further details on this publication in a note \cite[p. 87]{Oppo}.

Florensky's dissertation was well received by his teachers. Egorov \index{Egorov, Dmitri Fedorovich} considered him as one of the most brilliant and original mathematics students. One of Florensky's teachers,\index{Zhukovsky, Nikolai Yegorovich} Nikolai Yegorovich Zhukovsky,\footnote{Nikolai Yegorovich Zhukovsky (1847-1921)\index{Zhukovsky, Nikolai Yegorovich} was a mathematician, physicist and engineer. He is considered to be the founder of the Russian schools of hydromechanics and aeromechanics.  In 1904, he established, near Moscow, the first aerodynamical research institute in Europe. He is often called the ``father of Russian aviation", see \cite{Grigorian}. In the theory of conformal mappings, the transformation 
$z\mapsto \frac{1}{z}$ is called the \emph{Zhukovsky transformation.}
Zhukovsky's  work had a definite impact on Mikhail Alekseyevich Lavrentieff\index{Lavrentieff, Mikhail Alekseyevich}  (1900-980),  the Russian founder of the theory of quasiconformal mappings whom we shall mention later in this chapter.
Lavrentieff was a specialist of  partial differential equations and hydrodynamics and the founder of the
 Institute of Hydrodynamics of the Siberian Division of the Russian Academy of Sciences.}
offered him the possibility of pursuing a specialization curriculum which would lead him to a doctorate. But to the surprise of his teachers, 
Florensky turned down this opportunity and enrolled as a student at the Moscow Theological Academy. He delegated the function of secretary of the Students' Mathematical Circle to Luzin \cite{Demidov2014}.
The Theological Academy is situated in Sergiyev Posad, a city about 70 km northeast of Moscow which also hosts the Trinity-Saint Sergius Lavra.

Even though Florensky  distanced himself from the profession of mathematician, mathematics remained for him, for the rest of his life, a way of thinking and a basis for his philosophical and religious work.

Florensky's\index{Florensky, Pavel Aleksandrovich} strained relations with the state police started in 1906. That year, he was sentenced to three months' jail for having called for the abolition of the death penalty\footnote{Florensky made his statements against the death penalty in a homily, first delivered at the Lavra and then published clandestinely. Demidov  and Ford note in \cite{De-F} that ``the Moscow Theological Academy had the practice of allowing seminarians to preach in the Cathedral."} and for having protested, along with some fellow students, against the execution of Lieutenant Piotr Petrovich Schmidt,\index{Schmidt, Piotr Petrovich (Lieutenant)} a naval officer who had been one of the leaders of the Black Sea Fleet revolt in Sebastopol during the 1905 revolution.\footnote{Lieutenant Piotr Petrovich Schmidt, an idealist sympathizer of the 1905 revolution,  was  a typical representative of Russian intellectuals during this episode. He founded a \emph{Union of People's Friends Officers} and addressed several revolutionary meetings that culminated in a general insurrection of the Black Sea fleet crews on November 11, 1905.} 
Florensky spent a week in the Tagansky prison in Moscow (March 23 to 30, 1906), before being released after the intercession of several people, including the rector of the Theological Academy. During the week he spent in prison, he wrote a paper titled \emph{On the elements of the $\alpha$ number system}, see \cite[p. 601]{De-F}.

A few words on the 1905 social movements are in order, because these movements had a definite impact on several of the figures we will be talking about.

The origin of the 1905 social movements  goes back several years. In 1881, the year before Pavel Florensky\index{Florensky, Pavel Aleksandrovich} was born, Tsar Aleksander II,\index{ Aleksander II (Tsar)} who was ready to give his country a constitution, was assassinated by members of the \emph{Narodnaya Volya} (``People's Will"), a late 19th-century revolutionary socialist political organization. With this assassination, Russia started being in the grip of a political effervescence and a state of siege that culminated in the bloody events of 1905 that left their mark on Florensky, as they did on a large number of Russians. The 1905 revolt started with a peaceful demonstration by workers who had a request for the Tsar, and it ended in a terrible shoot-out. What happened during these events was later called the First Russian Revolution, and was regarded as the beginning of Imperial Russia's agony. 

In a sense, the 1905 revolution was a culmination of the idealistic traditions of the Russian intelligentsia.
   Chapter 2 of Pasternak's\index{Pasternak, Boris! \emph{Doctor Zhivago}} \emph{Doctor Zhivago},  titled \emph{A Girl from a Different World}, starts with the words ``The war with Japan was not yet over when it was unexpectedly 
overshadowed by other events. Waves of revolution swept 
across Russia, each greater and more extraordinary than the 
last."\footnote{English translation from \cite{Pasternak}.}
 In relation with the events we are recalling, it is also interesting to know that in 1925-1927, Pasternak\index{Pasternak, Boris! \emph{Lieutenant Schmidt}} wrote two long poems dedicated to the 1905 revolution, one of them titled \emph{Lieutenant Schmidt}.
 
Florensky obtained his theology diploma in the spring of 1908, with a dissertation titled \emph{On spiritual truth}, a work that was expanded later to a treatise titled\index{Florensky, Pavel Aleksandrovich!\emph{The pillar and ground of the truth}}  \emph{The pillar and ground of the truth: An essay in orthodox theodicy in twelve letters} \cite{Flo1},  which became his first doctoral dissertation, defended in 1914. This work had a definite impact on the Russian theological circles.

\begin{figure}[htbp]
\centering
\includegraphics[width=8cm]{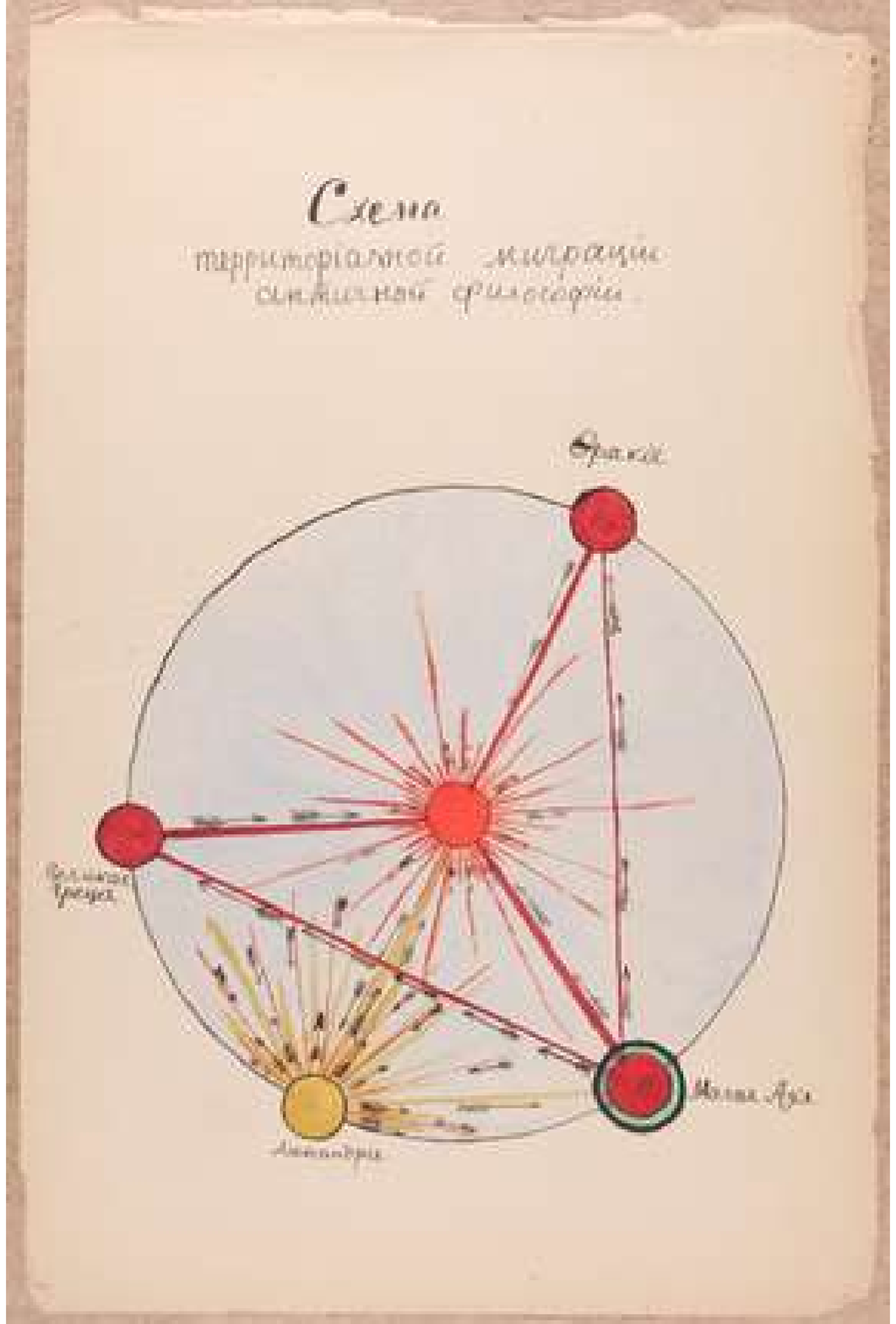}
\caption{From the album of illustrations of P. A. Florensky's course on history of classical Greek philosophy 1908-1909. The translation of the page title   is ``The scheme of territorial migration of ancient philosophy". Around the circle, starting from  the top, clockwise, one reads: Thrace, Asia Minor, Alexandria and Great Greece.
Paper, watercolor, pencil and ink. 
Multimedia Art Museum, Moscow.}
\label{fig:scheme}
\end{figure}

 Upon his graduation in 1908, Florensky \index{Florensky, Pavel Aleksandrovich}was entrusted with the teaching of history of philosophy at the Moscow 
Theological Academy, a position he held until the Academy's closure in September 1918. 
Figure \ref{fig:scheme} is extracted from the album of illustrations 
of his course on the history of classical Greek philosophy.

Florensky was ordained priest in 1911, after his marriage. 
Just as philosophy accompanied him during his mathematical studies, mathematics constituted an important element of his background during his theological studies. The discussions around the mathematical notion of discontinuity that were taking place at that time, the recent developments of set theory,\index{set theory (historical development)}   the\index{set theory (historical development)} language of mathematical logic and other mathematical subjects were at the basis of his treatise \emph{The pillar and ground of the truth}.\index{Florensky, Pavel Aleksandrovich!\emph{The pillar and ground of the truth}} Since the book was addressed to theologians rather than mathematicians, Florensky included there several sections containing technical mathematical explanations.
For instance, Section XV is titled \emph{Certain concepts from the theory of infinite}, 
Section XVII is titled \emph{Irrationalities in mathematics and in dogma}, 
Section XIX is titled \emph{The concept of identity in mathematical logic}, and 
Section XXX is titled \emph{The basic symbols and elementary formulas of symbolic logic},
etc. I shall say more about this work below.

In 1917, the teaching of  theology in Russia entered a period of upheaval as a result of the October revolution, and the Academy closed for good in September 1918.\footnote{The Moscow Theological Academy reopened in 1944. In 1930, the name of Sergiyev Posad, in reference to Saint Sergius of Radonezh, considered as the heart of Russia---It was from there that Russian culture expanded---was changed to Zagorsk, named after Vladimir Mikhailovich Zagorsky, a revolutionary who had been assassinated and whose  real name was Wolf Mikhelevich Lubotsky. Zagorsky was a party activist who had took part in the 1905 revolution. He was Secretary of the Moscow Committee of the Russian Communist Party (Bolsheviks) when he was killed, on September 25, 1919, by a bomb thrown by an anarchist. It was only in 1991 that the town of Zagorsk recovered its original name of Sergiyev Posad.}   Florensky\index{Florensky, Pavel Aleksandrovich} could no longer teach philosophy there. He became a teacher of history of Byzantine art at the MIKhIM.\footnote{This   is the acronym of Moskovskii Institut Istoriko-Khudozhestvennykh Izyskanii i Muzeevedeniia, that is, the Moscow Institute of
Historical and Artistic Research and Museum Studies.}
 In October 1919, he was appointed scientific secretary of the committee for the Preservation of Art Monuments and Antiquities of the Sergiyev Posad Lavra, which belonged then to the Commissariat of People Instruction. He prepared a treatise on religious art,\index{Florensky, Pavel Aleksandrovich!\emph{The reverse perspective}}  titled the \emph{Reverse perspective} \cite{Flo2}, which was the occasion for him to expand his ideas on space in drawing and painting.  I will say more about this work in \S \ref{s:Imaginaries} below.  
 This was also the epoch where Florensky became close to the writer Vassili Rozanov.\index{Rozanov, Vassili} The latter wrote about Florensky in a letter to the art critic Mikhail Spasovsky,\index{Spasovsky, Mikhail} quoted in part in \cite[p. 413]{Rozanov} and in the introduction of the French version of the \emph{Imaginaries} \cite[p. 15]{Flo2}: 
 \begin{quote}\small
 He is the Pascal\index{Pascal, Blaise} of our time, the Pascal of our Russia. You know, sometimes I even get the impression that he is a \emph{saint}, so exceptional is his mind, so out of the ordinary. \ldots I think and I am firmly convinced that he is even infinitely greater than Pascal, both in his essence---on the level of a Plato\index{Plato} with all that is so \emph{exceptional} about him---and in his intellectual discoveries, intellectual visions.
 \end{quote}

  By that time, Florensky had already produced a large amount of writings. He signed a contract with the Pomor'ye publishing house to have his works appear in 19 volumes. Eventually, only his book \emph{Imaginaries in geometry} \cite{Imaginaires}, which I shall review in Sections \ref{s:Imaginaries} and \ref{s:Dante} below, was published by Pomor'ye before its closing.

Pavel Florensky, with  his family, lived  in a small wooden house near the Sergiyev Posad monastery\index{Sergiyev Posad monastery}. For the rest of his life, including the years he spent in prison or in a Gulag, he used his scientific skills to work as a specialized technician in a plastics factory, as a teacher and encyclopedia editor, and as an engineer in several research institutes. He also taught mathematics and physics at the Zagorsk Pedagogical Institute. With all the troubles he had to deal with, work was what allowed him to survive psychologically, besides supporting his family.  At some point, he was called as an engineer on an important project to electrify Russia, which was known under the name GOELRO.\footnote{GOELRO is the transliteration of the Russian abbreviation for ``State Commission for the Electrification of Russia."} There, he had the occasion to meet Stalin,\index{Stalin, Joseph} who considered him,  like many intellectuals and clerics, an enemy of the people, but who tolerated him  because the state needed his skills. Trotsky,\index{Trotsky, Leon} at some point, took an interest in Florensky,\index{Florensky, Pavel Aleksandrovich} and to some degree protected him, before he was driven out himself.  The following episode is reported in the book by V.  Chentalinsky who collected KGB documents concerning literary figures when the KGB archives were made available, in the 1990s \cite[p. 174]{Chentalinsky}:

\begin{quote}\small
One day, long before his disgrace, Leon Trotsky\footnote{Trotsky was banished in 1929.} was astonished to see Florensky's white cassock. 

``Who's that?" he asked.

The leader of the world revolution, who thought he was as wise as King Solomon, held a multitude of positions. He was also head of the \emph{Glavelectro},\footnote{This is the Central Directorate of Electricity.} and as such visited the institute where the priest worked.

``This is Professor Florensky",  came the reply. 

-- ``Ah yes, I know..." 
Trotsky\index{Trotsky, Leon}  approached him and invited him to an engineering conference.

``Naturally,  he said, you wouldn't come dressed like that...

-- I haven't renounced my vows and can't wear any other clothes," Florensky\index{Florensky, Pavel Aleksandrovich} replied.

``-- Oh, you can't, \ldots  so come in this outfit \ldots." 

When Florensky stepped up to the podium to address the conference, a murmur of astonishment ran through the hall: a priest at the pulpit! And, more than by his brilliant communication, the audience was struck by the enigma he represented: a cleric, and therefore an obscurantist, who possessed such knowledge of the exact sciences!
\end{quote}

Trotsky's\index{Trotsky, Leon}  indulgence for Florensky would later prove fatal, as we shall recall below.

Florensky,\index{Florensky, Pavel Aleksandrovich} besides carrying out several assignments for the state, was teaching drawing and history of art, and most of all, he was writing and writing. At some point, he was offered to leave the country on one of the ``philosophers' boats"\footnote{The ``philosophers' boats"  is the name that was given to the steamships that transported intellectuals expelled by the Soviet government in 1922, soon before the creation of the Union of Soviet Socialist Republics. Leon Trotsky had decided that ``there was no pretext for shooting these people, but it was no longer possible to tolerate them", cf. \cite{Nikolski}.} that had already taken a large part of the Russian Intelligentsia.\footnote{Among those who were expelled and with whom Florensky was closely  associated, we mention the symbolist poet Vyacheslav Ivanov\index{Ivanov, Vyacheslav Ivanovich} (1866-1949),  the religious philosopher Nikolai Berdyaev\index{Berdyaev, Nikolai Aleksandrovich} (1874-1948) and Sergei Bulgakov\index{Bulgakov, Sergei Nikolaevich} (1871-1944), who was a theologian, philosopher and economist. We also mention the philosopher Siméon Frank\index{Frank, Siméon Lyudvigovich} (1877-1950), the theologian
Vladimir Nikolaievich Lossky\index{Lossky, Vladimir Nikolaievich} (1903-1958), the writer and  journalist
Mikhail Andreyevich Osorgin (1878-1942), the sociologist  Pitirim Alexandrovich Sorokin\index{Sorokin, Pitirim Alexandrovich } (1889-1968) and Valentin Fyodorovich Bulgakov\index{Bulgakovn Valentin Fyodorovich} (1886-1966), a Christian anarchist and antimilitarist who had been the last secretary and biographer of Leo Tolstoy. A large number of the persons expelled participated to the 1905 revolution and belonged to the Socialist Revolutionary Party.} He refused, declaring that he wanted to stay in his country. Of this refusal, Sergei Bulgakov, who was his disciple and friend and who had been expelled, in 1922, said in an address he gave in Paris in 1943, when the news of Florensky's death reached him \cite{Bulgakov}: ``He had not the shadow of naivety and idyllic primitiveness, [\ldots] but he loved this native land, like the Mother of all men, the Demeter of the Ancients." 
 In staying in his country, Florensky knew what he was expecting. He was not put straight into prison, and he was always given new tasks, even when he was in the Gulag.

At dawn on May 21, 1928, Florensky was arrested and taken to the Lubyanka.\footnote{This is the building in Moscow that hosts the headquarters of the Political Police, called State Political Directorate (known by the acronym GPU), with its affiliated prison.} 
 The episode is recounted in the book \cite{Chentalinsky}. He was required there to fill in a questionnaire. To the question of whether he had ever been prosecuted, he replied: ``Yes, in 1906, for a sermon against Lieutenant Schmidt's\index{Schmidt, Piotr Petrovich (Lieutenant)} death sentence". This could have played in his favor, but at the same time he added: ``I have practically nothing to say on political questions. Because of my character, my occupations and my conviction, based on the lessons of history, that events still develop according to as yet unclarified laws of social dynamics and not according to the will of the participants, I have always avoided the political struggle." In those times, it was dangerous to declare oneself politically neutral. Thus, Florensky was arrested along with 79 other people, accused of a plot in connection with the shooting of a revolutionary. The minutes of his arrest at that date and of later arrests are published in the book \cite{Chentalinsky}.

\begin{figure}
\begin{center} 
\includegraphics[width=10 cm]{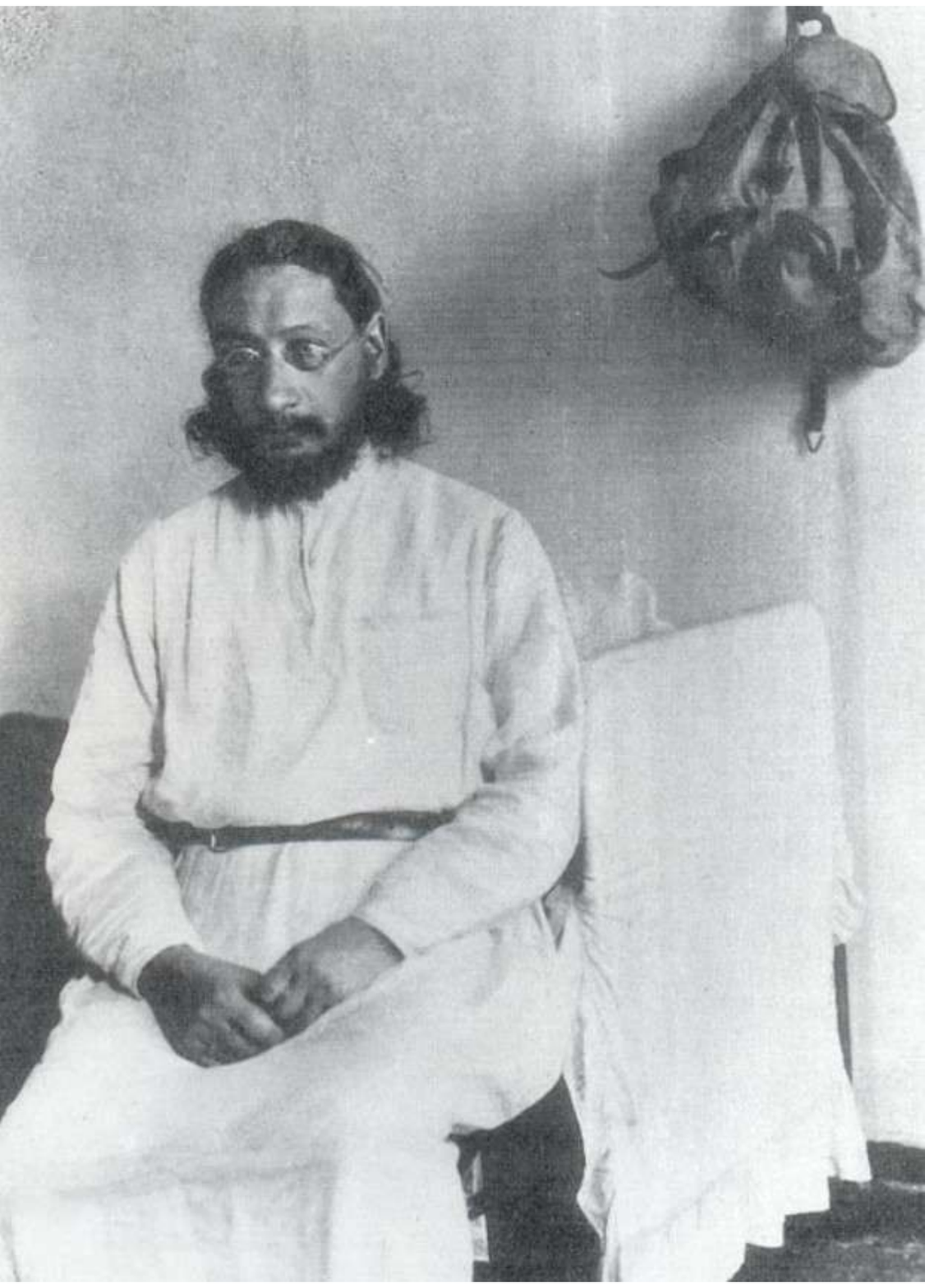} 
\caption{Pavel Florensky in exile in Nizhny Novgorod, 1928}
\label{fig:Nizhny}
\end{center}
\end{figure}

Florensky was placed under house arrest in Nizhny Novgorod. The photo is reproduced in Figure \ref{fig:Nizhny} was taken there. A few months later, he was granted early release, after to the intervention of influential people, among them Ekaterina Peshkova, Maxim Gorky's ex-wife. Peshkova was 
at that time director of the political Red Cross and  was on good terms with Stalin. Florensky\index{Florensky, Pavel Aleksandrovich} returned to Moscow and resumed his work for the state.
 In 1930 he was made deputy director of the Electrical Engineering Institute K. A. Krug in Moscow and in 1931  member of the central directorate for the study of electro-insulating material. 
 Meanwhile, he continued to teach philosophy clandestinely at the Danilovski and Petrovski monasteries, where the Theological Academy continued to exist informally, after it was officially closed  in 1918. Florensky continued to publish in scientific journals until his deportation. In  1932, he published several articles in the   journal \emph{Socialist Reconstruction and Science}, among them an article in which he describes an analogical calculator for finding approximate solutions of 
algebraic equations of high degree, see \cite[p. 608]{De-F} and the reference there.
 
 Between 1927 and 1933 (the year of his deportation), Florensky filed almost 50 patent applications, but most of them were rejected because he was a priest, automatically accused of obscurantism.

The existence of a cleric, always dressed in cassock, who was a scientist and useful to the nation, became more and more unbearable for the Soviet authorities, and at some point they ran out of patience. At dawn on February 26, 1933, Florensky was arrested at his home and beaten up at the Lubyanka, again accused of plotting and anti-revolutionary propaganda.   V.  Chentalinsky, based on the KGB documents, writes about his arrest \cite[p.164]{Chentalinsky}: 
``Today we can see Florensky's true face, but his case file presents him as a dangerous criminal for society, as an obscurantist placed under permanent and strict surveillance until his death."

This time, Florensk\index{Florensky, Pavel Aleksandrovich}y spent six months at the Lubyanka. During his captivity, he worked around the clock. On July 26 of the same year, after a sham trial based on false accusations, he was sentenced to ten years in the Gulag, with the obligation to continue his scientific work.\footnote{It is probably time to recall that the word Gulag is an acronym  for ``Glavnoye Upravleniye Lagerey", which means Main Administration of Camps. In theory, the Gulags were ``camps for re-education through work".} He was again offered the possibility to leave Russia. In fact, Tomáš Masaryk, the president of the Czechoslovak Republic, offered to take him in, so that he could continue his scientific work there. But once again, Florensky refused to leave his country. 
In August of the same year, he was transported to Eastern Siberia on a train full of common criminals. He arrived at his destination, a camp whose name translates as ``free," exhausted, starving and penniless: everything he possessed was stolen during the journey. 
In February 1934, that is, a few months after his arrival to the camp, he learned from a courier that his library had been confiscated during a search at his family home in Zagorsk. Shocked, he wrote to the OGPU:\footnote{The OGPU is an acronym for ``Joint State Political Directorate", which was 
 the Soviet Union's state police.} 
 \begin{quote}\small
 My whole life has been devoted to science and philosophy, and I have never known rest, leisure or pleasure. Not only did I devote all my time and energy to the service of humanity, but I also spent a large part of my modest salary on buying books, etc.  My library was not simply a collection of books, but selections for specific, established subjects. One might say that many of the works I intended to write were already half-finished in the form of notes in books whose arrangement on the shelves, the intellectual key, was known to me alone [\ldots] My whole life's work is now destroyed [\ldots] This is for me far worse than physical death. \cite{Misler} 
 \end{quote}

 In November 1934, after seeing for the last time his family who was visiting him, Florensky was forbidden visitation rights and sent to the Solovki Islands,\index{Solovki Islands} an austere place in the White Sea, 160 km from the Arctic Circle. There lies an ancient monastery, known under the name \emph{Solovetsky}, surrounded by an enclosure of walls and towers and encircled by immense frozen expanses. In the days of the Tsars, it was the place where recalcitrant monks were sent.
  In 1920, the monastery was transformed into a Gulag, and it remained so until 1937, the year Florensky\index{Florensky, Pavel Aleksandrovich} was shot. Prisoners there were crammed into the former monks' cells, several to a cell. Later, the monastery came to be regarded as the most terrible of the Soviet camps, where half a million people perished.  Aleksandr Solzhenitsyn\index{Solzhenitsyn, Aleksandr} called this place ``the mother of Gulags". 
  
  About Florensky, Solzhenitsyn wrote that he was ``perhaps the most remarkable person devoured by the Gulag". Chapter 2 of the second volume of his \emph{Gulag Archipelago}\index{Solzhenitsyn, Aleksandr!\emph{The Gulag Archipelago}} opens with a historical account of the Solovki Islands, where the story of the Gulag starts \cite{Solzh}: 
  \begin{quote}\small
On the White Sea, where the nights are white for half a year
at a time, Bolshoi Solovetsky Island lifts its white churches from
the water within the ring of its bouldered kremlin walls,\footnote{The Russian word ``kreml'' (kremlin) refers to a fortified complex or citadel, and not just the one in Moscow.} rusty-red from the lichens which have struck root there---and the grayish-white Solovetsky seagulls hover continually over the kremlin and screech.
\end{quote}

On 10.X.1934, Florensky sent his first letter to his family from his new exile. He was still in Kem, 200 km as the crow flies from the Solovki Islands.\index{Florensky, Pavel Aleksandrovich} He wrote \cite[Letter No. 1, p. 13]{Florensky-Letter}: 
\begin{quote}\small 
On arrival at the camp, I was robbed by an armed attack and threatened with three axes, but as you can see, I am safe and sound, albeit without luggage or money; some of my belongings have been found. All this time I have been cold and hungry. It was much harder than I had imagined when I left the Skovorodino train station.\footnote{Skovorodino is a city situated on the  Trans-Siberian Railway, 7306 kilometers from Moscow.} I was supposed to go to Solovki, which in itself would not be so bad, but I am held up in Kem where I am busy making and filling out control forms.
\end{quote}

\begin{figure}
\begin{center} 
\includegraphics[width=14 cm]{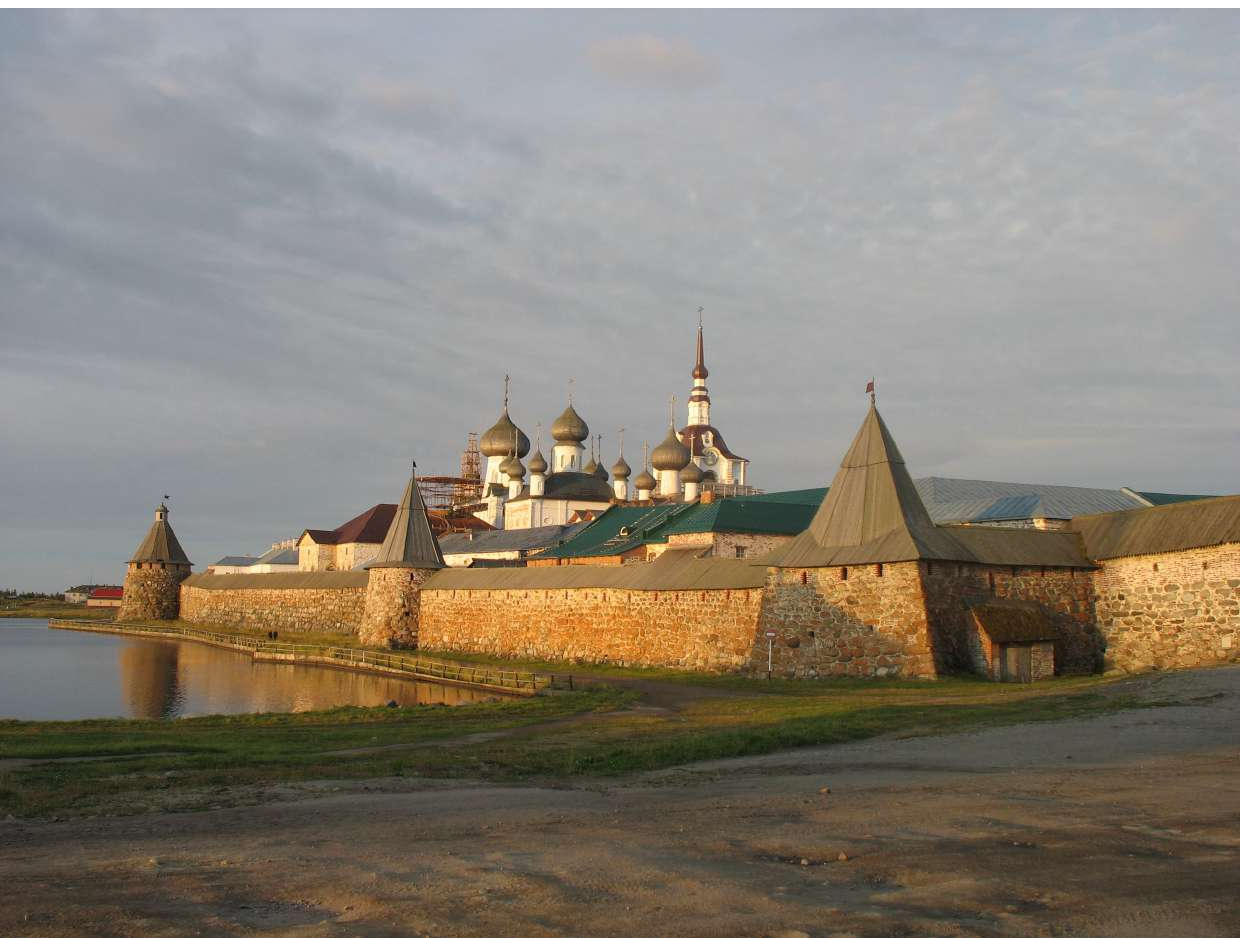} 
\caption{The Solovetsky Monastery, where the prisoners were crammed together.  This recent and romantic photograph does not give any idea of the atmosphere in the monastery during the time it was used as a Gulag,  between 1926 and 1939, but the sky above it and the sea surrounding it have not changed. Photo: Cultural heritage object in Russia, No. 2910179018, GitaOggg.}
\label{f:Solovki}
\end{center}
\end{figure}

In the Gulag of the Solovki Islands,\index{Solovki Islands} a new life for Florensky started, although he spent his time on exactly the same activity as before: working.  
 I have reproduced here a recent photograph of the Solovetsky monastery, see Figure \ref{f:Solovki}.

  In a letter dated 8-9.VII.1936, addressed to his wife \cite[Letter No. 67, p. 435]{Florensky-Letter}, he\index{Florensky, Pavel Aleksandrovich} writes: ``I never see the sun or nature, even through the window: I am constantly busy. All my time is spent developing and deepening technological processes, equipping the algae workshop as much as possible, given the increasing demands for quality and quantity." Besides,
Florensky also set up a theater for the prisoners and taught classes there, keen to passing on his knowledge to the other prisoners.   His work force sustained and guided him through all his hardship, as the enormous correspondence with his family shows. The prison authorities, aware of Florensky's talents and keen to exploit them to the full, arranged workshops for him and let him organize his courses and lectures. 
In a letter addressed to his son Vassili, dated 27-28.IV.1935, he writes \cite[Letter no. 17, p. 111]{Florensky-Letter}: ``I just gave a talk on Maxwell's theory and related questions. I got home at 11:30 p.m. and it was perfectly clear outside. [\ldots] The light here is transparent: everything seems to be seen in a dream, and in fact non-existent. This is particularly true at night, but also during the day, when the sun is shining." He would have liked to do mathematics, but this was not easy. In a letter to his mother, dated  24.VI.1936, Florensky writes \cite[Letter no. 65, p. 422]{Florensky-Letter}: 
\begin{quote}\small
I deal a bit with mathematical questions; for more serious work in this field, I have neither time, nor the necessary books, and finding everything in one's head is too unproductive a task. With the most advanced workers at Idoprom, where I work, i.e., with those who have completed eight years of school or technical institute, or a few years of higher education, I have organized a physics circle with the aim of learning how to reactivate one's knowledge of physics, mathematics and a little chemistry, really assimilate them, deepen them and complete them. So, the sum of the various tasks and activities occupies all the time, not only during the day, but also most of the night. I try, after dinner, from 7-7.30, to get about two hours of sleep, so that I can work with a fresh mind at night.
\end{quote} 

He accepted his new situation philosophically. Perhaps with a touch of humor, but also with great sadness,\index{Florensky, Pavel Aleksandrovich} he wrote to his wife on 7.XI.1934 \cite[p. 17]{Florensky-Letter}: 
\begin{quote}\small
I have long since come to the conclusion that desires in our lives are realized, but with great delay, in a caricatured way that makes them unrecognizable. In recent years, I wanted to live in the immediate vicinity of a laboratory, and that came true; in Skovorodino. I wanted to take care of the soil, and that came true here too. Before that, I wanted to live in a monastery, and that came true, but in Solovki. As a child, I used to fantasize about living on an island, seeing the tides, handling seaweed, and now here I am on an island, there are tides and I may soon be handling seaweed. But the realization of desires is such that one no longer recognizes them, and it happens when the desire has already passed. I send you a big kiss, my dear Annoulia, don't lose courage.
\end{quote}

In fact, early in 1920, he had written in his autobiography about this call of the sea \cite[p. 41]{Florensky-bio}:  
 \begin{quote}\small
 The call of the sea, the fractional sound of the surf, the infinite, glittering surface in which I could distinguish ever-smaller flakes, down to the tiniest particles, but which never merged, all this resonated endlessly in my soul. My body was craving the saltiness of the sea, the salty, iodine-soaked air, but also the fractionated air carrying the tiniest crystals: sometimes it feels good to hold at least one vial of tincture of iodine close to you. I crave the taste of the sea\ldots  It seems to me that if I could find a pile of seaweed, I would devour them all.
 \end{quote}

Florensky\index{Florensky, Pavel Aleksandrovich} set up a laboratory to process seaweed and extract iodine and agar-agar, a substance used as a glue, and also in cooking. He conducted research in 
 biology and organic chemistry. During his stay in Solovetsky, he registered seven patents.  Among the prisoners were engineers,
 scientists and academics, educated people as well as workers; he gave courses to all of them, working at an inhuman pace. Since he worked more than double the required time, he was allowed to send three letters a month to his family. He sent long letters, sometimes a dozen of pages long, to his wife, mother and children, and also, from time to time,  to friends such as the pianist Maria Yudina,\index{Yudina, Maria Veniaminovna} who was Stalin's favorite pianist and who had become Florensky's spiritual daughter. He also corresponded with Vladimir Ivanovich Vernadsky,\index{Vernadsky, Vladimir Ivanovich} a member of the  Russian and Soviet Academies of Sciences.\footnote{Vladimir Ivanovich Vernadsky (1863-1945) was a famous mineralogist and  geochemist. He is the author of a book titled  \emph{The Biosphere}, published in 1936,  in which he defends the ideas that life and human cognition are the forces that are at the origin of the Earth's development and the transformation of the biosphere. Vernadsky was one of the fathers of the Soviet atomic bomb. He never sought to hide his friendship and respect for Pavel Florensky.}
 In his letters to his family, he used to describe his everyday life in the Gulag, giving advice to each family member.

In some of his letters, Florensky\index{Florensky, Pavel Aleksandrovich} gave mathematics and physics lessons to his children. For instance, in the letter dated 5.XII.1935, to his son Vassili \cite[Letter no. 40, p. 245]{Florensky-Letter}, he formulated a problem in mathematical physics and he developed its solution. In the letter dated 24.III.1936 \cite[Letter no. 54, p. 341]{Florensky-Letter}, he explained a mathematical problem to his daughter Olga.
Several times, in his letters, he enquired about his confiscated books. On September 16, 1935, he wrote to his wife: ``On occasion, enquire about the fate of my work on permagel and my mathematical articles, and in particular the articles on four-color maps." The works were never found. On February 25, 1936, he wrote, again to his wife \cite[Letter No. 50, p. 311]{Florensky-Letter}: 
\begin{quote}\small
Now that I live in the laboratory, I can also work at night, especially as the work is done 24 hours a day and needs to be supervised all the time. At the moment, it is already nearly 3 o'clock. Next to me, a worker is cutting the agar into pieces for freezing, and on the other side, pumped water is gurgling in a large tub. [\ldots] Try not to get bored and to find interest in work in addition to that which is immediately essential. Life must be embellished.
\end{quote}

He tried to keep abreast of developments in all the sciences, and every small information he could get was important for him.
In the same letter, he wrote: 
\begin{quote}\small
In the Pravda of February 9, 1936, No. 39 (6645) p. 5, I came across the composition of a liquid for toothache, that of Professor Hartman, an American from Columbia University. As this liquid is very easy to obtain and, according to the newspapers, gives 80 \% positive results, you need to have some on hand.  Here is the formula: one part etyl alcohol, 2 parts ether, 1.25 parts timol.
\end{quote}

 He is worried again about his books: ``Have you received my books and materials from P. N.? Ask him to send them to you, because I have spent a whole year putting them together, and I want the boys to use them: I am already giving them nothing, so do not neglect the little they can get." He is thinking about his past, but also about the past of the universe. In one of his last letters (letter dated  20.II.1937, \cite[letter No. 71, p. 455, ]{Florensky-Letter}), Florensky writes,\index{Florensky, Pavel Aleksandrovich} to his son Kyrill: 
 ``A distant past absorbs my thoughts: the formation of the island's rocks some two billion years ago, the history of the formation of the Earth's crust in these very special places, the appearance of new rocks based on electronic waste, the movements of the glacier, its advance and retreat." With his meager salary, he is also sustaining his family. In the same letter, he writes to his son: ``  Did you receive the small sum I sent you?" 

  V. Chentalinsky describes  the last days of Father Florensky. We are in June or July 1937. The camp had been recently reorganized under a new management which had recently arrived.  The iodine production plant, where Florensky was the leading scientist, was closed. The detention conditions were tightened.
 Prisoners were herded off to an unknown destination.   Florensky\index{Florensky, Pavel Aleksandrovich} was arrested after a report, made by an informer. Chentalinsky  recalls this denunciation in the following terms \cite[p. 191]{Chentalinsky}:

\begin{quote}\small
    Florensky and another inmate walked to the library, talking loudly and gesturing forcefully. The informer ``Comrade" follows them closely and listens attentively. They are talking about the upcoming war:

``The predictions of Trotsky,\index{Trotsky, Leon}  the eminent Party strategist and ideologist, according to which war will soon break out seem justified to me, says Florensky. That is the rule: war breaks out periodically every fifteen or twenty years \ldots"

 This short, seemingly innocent report had unfortunate consequences. It was used to draw up a special attestation on Florensky's conduct: 
 
 ``Engaged in counter-revolutionary propaganda in the camp, praising the enemy of the people Trotsky." [\ldots]

\end{quote}
 
Very soon after this, Florensky\index{Florensky, Pavel Aleksandrovich} witnessed the destruction of his research facilities and learned that his discoveries on algae (in Solovki) and permagel (in Siberia) have been attributed to others. References to him in the \emph{Technical Encyclopedia} of which he was an editor were cancelled in the second edition, see \cite[p. 608]{De-F}. Florensky  was again sentenced, this time to death, by a three-judge emergency court. He was immediately transferred to the mainland, in the Leningrad region, with a large group of prisoners from Solovki.
 The execution act was signed, and on December 8 of that year,  he died riddled with bullets, along with 508 other people, on the outskirts of Leningrad. The execution act is reproduced in Chentalinsky's book \cite[p. 194]{Chentalinsky}: ``The ruling of the troika of the Leningrad region NKVD leadership concerning the death row inmate, Pavel Aleksandrovich Florensky was executed on December 8, 1937, as attested by the present deed."
 
  Florensky's body, along with those of the other shot, was dumped in a mass grave, which has since become a place of memorial. We were in the middle of the period known as the ``Great Stalin Terror'' (1937-1938), during which 750,000 people were executed in the USSR. 
  
Meanwhile, life went on at the mathematics departments of Leningrad university, a few kilometers from that place. At the beginning of his historical review \cite{Zalgaller}, V. A. Zalgaller writes:
  \begin{quote}\small
  
  In 1937, enrolment in the Mathematics and Mechanics Faculty of Leningrad
University (Mat-Mekh) was unusually high: 210 students. The ``big Chemistry
lecture hall" in the main building on the Neva embankment hardly contained us. At
the beginning, analytic geometry was taught by Professor A. R. Kulisher (1878-1937). But in September, he was already arrested and replaced by a 25 year old
Professor Aleksandr\index{Aleksandrov, Aleksandr Danilovich} Danilovich Aleksandrov. [\ldots]
\end{quote}

On May 5, 1958, the Moscow City Court overturned Florensky's conviction of July 16, 1933. The rehabilitation court order reads: \cite[p. 184]{Chentalinsky}:
\begin{quote}\small
Florensky's arrest (like that of the others judged for the same case) is not based on any document in the case file. Witnesses were not questioned, and those who took part in the investigation of the case were subsequently convicted of forgery. Florensky (and the others) were unjustly convicted, in the absence of proof of their guilt.
\end{quote}

On March 5, 1959, the Arkhangelsk Oblast Court overturned his conviction of Nov. 15, 1937.

     It was only in 1989 that Florensky's\index{Florensky, Pavel Aleksandrovich} family and descendants learned, through a letter from the KGB, of the date and place of his execution.

\begin{figure}
\begin{center} 
\includegraphics[width=10 cm]{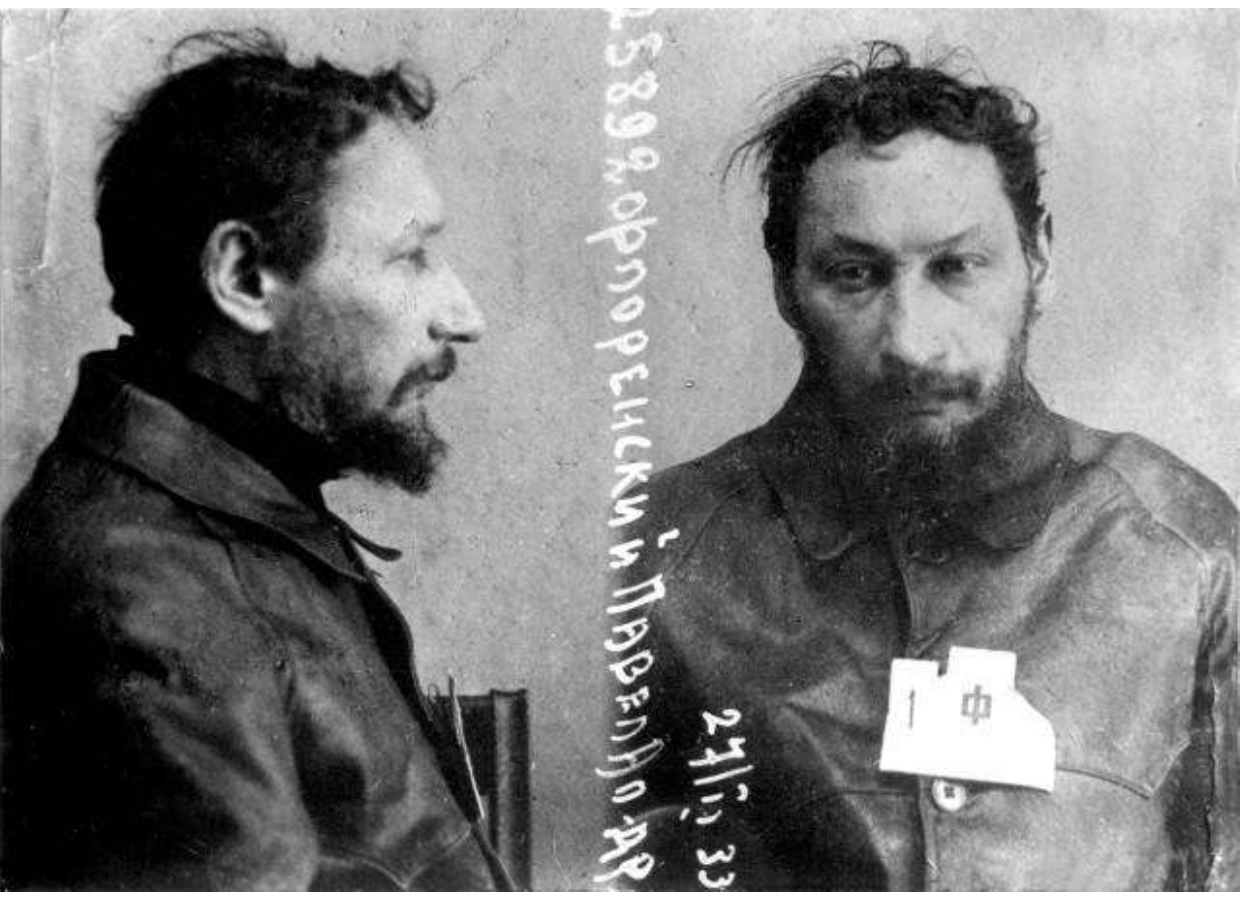} 
\caption{Pavel Florensky, shortly after his arrest in 1933}
\label{f:Florensky}
\end{center}
\end{figure}

 In this section, I tried to draw, with words, some features of Florensky's\index{Florensky, Pavel Aleksandrovich} face. For an extended biography, the reader can refer to \cite{Pyman}. 
 A complete biography should take several volumes, and is still to be written. There are a number of portraits of him. I have reproduced three of them here, in Figures \ref{fig:Nizhny} and \ref{f:Florensky}. 
 
 \section{Imaginaries} \label{s:Imaginaries}
 The book, for which I shall refer to the French version \cite{Imaginaires}, is titled 
 \emph{Imaginaries in geometry},
 with two subtitles: 
 \emph{Extension of the domain of geometric images to two dimensions} and \emph{An essay of a new materialization of the imaginaries}. The work 
  is complex, as are all  of Florensky's writings,\index{Florensky, Pavel Aleksandrovich} because of the continuous  intermingling of philosophy, science and symbolism.\index{symbolism} 
  In the following pages, I shall extract a few ideas from this book, which deserves a very thorough analysis that would go beyond the scope of the present essay.
  Villani, who wrote the preface to the French edition, declares there: 
  ``We find the exposition of  a brilliant dreamy dissertation on the notion of duality, illustrated with advanced mathematical concepts, all designed to lead the reader's mind to the point where it can absorb a daring metaphysical thesis evoking a duality between matter  and ideas".  
   The physicist Pierre Vanhove, in the introduction to the same edition, writes:
``This singular book should not be read as a work of mathematics or physics, but rather as an attempt to use the geometry of the imaginary to analyse the fundamental notion of duality between the real and spiritual worlds."
At the same time, the \emph{Imaginaries} is an illustration of the fact that Florensky's mathematical works do not consist of theorems and proofs, as one would usually expect from a mathematician; the mathematics that it is about is non-mainstream, close to art and to symbolism;\index{symbolism} it is the language of a new imaginary world. To mention the opinion of a contemporary to Pavel Florensky, let me quote a fragment from a letter that Luzin\index{Luzin, Nikolai Nikolaevich}  wrote to his wife about Florensky, after a visit he made to him in 1908; the fragment is reproduced in an article by C. E. Ford in \cite{Ford}: 
 ``[\ldots] But as soon as he showed me his works in mathematics, my old opinion came back to me: that all his works are of no value
in the area of mathematics. Suggestive hints, beautiful analogies---very attractive
and promising, provocative, beckoning, but without any results." On the connection between Luzin\index{Luzin, Nikolai Nikolaevich} and Florensky see the article \cite{Demidov1} by Demidov,  Parshin,  Polovinkin and Florensky\index{Florensky, Pavel Aleksandrovich} Jr., the discussion in \cite{Sinkevich} by Galina Sinkevich, and  \S \ref{s:Luzin} below.
Establishing a relation, which, in Florensky's \emph{Imaginaries}, takes the form of a mathematical duality,  between matter and ideas, between the rational and spiritual worlds, etc., are objectives that Florensky sought to achieve throughout his life.

  Let me start with a clarification concerning the word ``Imaginaries" in the title of the work.
 This word, in the first few pages of the book, has the meaning that mathematicians give it, when they talk about an ``imaginary number", that is, the ``imaginary" refers to one of the two coordinates of a complex number. Later in the book, this word designates the world of Platonic\index{Plato} ideas, which Florensky considers as living in a space which is ``inverted" with respect to our own. The fact that a complex number is the sum of a real and an imaginary number becomes a symbol of the fact that the world consists of a combination of concrete and abstract elements.\footnote{\label{f:Bombelli} Let me add here the following remark made to me by Stelios Negrepontis: ``The sensible entities according to Plato are the imperfect approximations of the true Being to which they participate. The eventual construction of an intelligible true being employs a quadratic conjugation between binomial and apotome lines, according to\index{Theaetetus}  Theaetetus' Book X of Euclid's \emph{Elements}, a conjugation that inspired Rafael Bombelli\index{Bombelli, Rafael} around 1558 in his discovery of  complex numbers, which Bombelli called `binomio', exactly because of this similarity."}

In the same way  as Florensky declares that the real and imaginary numbers constitute the two sides of the same surface, he\index{Florensky, Pavel Aleksandrovich} sees the two worlds, the physical and the spiritual, as the two sides of the same reality.   
  Gradually, in the book, the word ``imaginary" takes on a literal meaning, as something that does not exist in reality, but rather as a tale that belongs to the realm of our imagination.
As a result, the discussion quickly turns philosophical and symbolic, and the book becomes a poetic, historical and cosmic metaphor on space and time, and on man and mankind.
The book culminates in a fabulous reference to Dante's \emph{Divine Comedy},\index{Dante!\emph{Divine Comedy}} discussing the Florentine poet's geometrical vision. Talking about geometry in Dante's work has caused considerable ink flow in recent years, by authors who have taken up Dante's\index{Dante} allegory, without any reference to Florensky and missing an essential point of what he actually discovered; I will discuss this in \S \ref{s:Dante}.

 The  \emph{Imaginaries} has also a cosmological side, and it
 may come as a surprise to the reader that Florensky, in the last section of the work, invoking Einstein's theory of relativity, defends Ptolemy's geocentric system against that of Copernicus.\footnote{Florensky had to explain and defend himself before the state authorities on this issue, see Footnote \ref{n:Postface}.}  This section of the book, which is its culmination, should be taken (as several parts of the work) as a metaphor,\index{Florensky, Pavel Aleksandrovich} continuing there the exposition of Dante's poetic and spiritual worldview.\index{Dante}  In his commentary, Vanhove, restoring Florensky's considerations in context, writes that ``far from having misunderstood the relativity revolution, Florensky wants to make Einstein's formalism compatible with his spiritual vision of the world."

  Let me review the book in some more detail.
    
  I start with the part that concerns complex numbers.
  
  In this part, Florensky, showing his great erudition, refers to a host of mathematicians who worked on the foundations of the theory of complex numbers and their various uses. He surveys, in \S 1 and especially in the first and second notes, works of Gregory, Marie, Truel, Gauss, Argand, K\"uhn, Wessel,  
   l'Abbé Buée,  Français, Warren, Lambert, Mourey, Bellavitis, Cauchy,  Matzka, Scheffler, M\"obius, Riecke, Hoüel and Gomes Teixeira. I know of no other work on this topic that contains such a large number of references.\footnote{Rafael Bombelli,\index{Bombelli, Rafael} whom we mentioned in Footnote \ref{f:Bombelli} before, is however missing in this list.}  ``The various \emph{interpretations} of the complex variable are", he says, ``interpretations through which this variable manifests itself symbolically." This leads him to a literary and poetic digression \cite[p. 36]{Imaginaires}:
   \begin{quote}\small
We know that several translations of the same poetic work in one or more languages, far from hindering each other, actually complement each other, although none of them can totally replace the original. Scientific representations of the same reality can and should also be increased without detracting from the truth. Knowing this, we have learned not to blame this or that interpretation for what it does not give, but to be grateful to it whenever we have the opportunity to use it.
   
  Yet we are forced to point out the limitations of an interpretation as soon as we notice the hypertrophy of this or that translation which attempts to identify itself with the original and substitute itself for it, i.e., which monopolizes a certain entity jealously cutting off all other interpretations.

   \end{quote}

   I will address later in this paper the question of interpreting poetry. The problem of the ``limitations of an interpretation" will also appear in the discussion on theology. 
   
   Talking about \emph{functions} of a complex variable, Florensky\index{Florensky, Pavel Aleksandrovich} discusses the works of Lejeune-Dirichlet, Riemann, Weiertrass and others.   He notes in passing that the notion of function as a ``correspondence" is ``insufficient, for if we take into consideration only the {\bf content} (the `material cause') of the function, we miss the essential point, the function itself as a whole, as a {\bf form} linking this content into a whole (the `formal cause')".\footnote{The bold face and the quotation marks, here and in other quotes, are Florensky's.} \cite[p. 36-37]{Imaginaires}

Then comes a part on projective geometry, a geometry dear to the Renaissance mathematicians and artists,\footnote{Leonardo da Vinci and Albrecht Dürer come to mind first, but there are many others, and I would like to mention just a few of those discussed by Florensky in his important book on religious art,\index{Florensky, Pavel Aleksandrovich!\emph{The reverse perspective}}  \emph{The Reverse Perspective} \cite{Flo2}:  Leon Battista Alberti, Filippo Brunelleschi, Paolo Uccello, Piero della Francesca, Donatello, Masaccio, Fillippo Lippi, Raphael (Raffaello Sanzio), Michelangelo Buonarotti, and there are others. The history of perspective is developed by Florensky in this book and in other articles referred to there.}  and it is not surprising that it fascinated Florensky, especially later on, he taught perspective drawing, and he wrote articles and a book on this topic.  The importance of projective geometry for the rest of this work will also become clearer when Florensky will explain, in the last \S of the \emph{Imaginaries}, with Dante's\index{Dante}  poem in the background, that the projective plane is an important part of our physical universe. I will comment on this below. Florensky proposes to develop analysis (in the mathematical sense of the term) on two sides of the same surface, to which he refers as a real side and an imaginary side. He attempts to extend the duality of projective geometry
  to a philosophical context.  Villani, in his Preface to the book, insists on the fact that Florensky's interest lies not only in the mathematical side of things: ``We must seek the reason for this duality in something else, something deeper".  Florensky expounds then a metaphor concerning the real world and the imaginary one, and the visible world and the spiritual one,  which finds its apogee in the evocation of the \emph{Divine Comedy} that I shall review below.     
  
In a relatively long section of the book, titled \emph{Explanation of the cover}, Florensky comments on the book cover,\footnote{ This cover won a Silver Medal at the 1925 \emph{Exposition des Arts Décoratifs} held in Paris.} the work of the artist Vladimir Favorsky,\index{Favorsky, Vladimir Andreyevich} a well-known graphist, woodcut illustrator and painter who was his friend and colleague at the VKhUTEMAS.\footnote{VKhUTEMAS is an abbreviation for the Russian expression ``Higher Art and Technical Studios". This is the state art and technical school where\index{Florensky, Pavel Aleksandrovich} Florensky lectured  from 1921 to 1924.} In this section, Florensky, reviewing the phenomena of perception and the psychology of vision, talks about a duality in concrete visual perception which is an expression of the dualistic nature of the geometric plane, where the visual images correspond to the real side of this plane, and the abstractly visual images to its imaginary side. 

Since we are talking here of projective geometry,\index{projective geometry} let me also note that Florensky's\index{Florensky, Pavel Aleksandrovich!\emph{The reverse perspective}}  \emph{Reverse Perspective} \cite{Flo2} contains several mathematical and historical notes on this topic and its complex and long development, with mention of works of Ptolemy\index{Ptolemy} who, in his \emph{Geography}, described the stereographic projection from the sphere onto the plane, together with other projections that were used in the art of drawing geographical maps, as well as works of 
Lambert, Loria, Aschieri, Enriques, Chasles, Poncelet, von Staudt, Fiedler, Wiener, Kupfer, Burmeister, 
Wilson, and the writings of the historian of mathematics Moritz Cantor\index{Cantor, Moritz}  on this subject.

Florensky\index{Florensky, Pavel Aleksandrovich} completed his book in July 1922. At that time he was teaching at VKhUTEMAS. It was the occasion for him to explain his ideas on space and time in art. Vassily Kandinsky\index{Kandinsky, Vassily} was teaching at the same institute. The scientist and the artist, who were both art theorists as well, were friends. 
Kandinsky was influenced by Florensky; there are several articles on this topic, see \cite{Taroutina, Sokolova}. The book \cite{Misler} contains a collection of writings of Florensky on art. 
 
 \section{Dante}\label{s:Dante}

  Section 9 of the \emph{Imaginaries}, which is the last section of the booklet, is dedicated to the 600th anniversary of the death of the founder of Italian poetry. This event was celebrated in Russia on September 14, 1921.  I will comment now on this section of Florensky's work.
  
A certain number of modern mathematicians and other authors have quoted passages from the \emph{Divine Comedy}, in order to point out a mathematical idea expressed by Dante\index{Dante}  there, though in a hidden way, namely, that the universe is a 3-dimensional sphere. Let me mention some of these authors, approximately in a chronological order: 
 The physicists I. Ozsv\'ath, and E. Schücking,  in the article  \emph{An anti-Mach metric}, published in the book \emph{Recent Developments in General Relativity}  \cite{O-Sch} (1962), and the same authors 34 years later, in 
the article \emph{The world viewed from outside} \cite{O-Sch-Dante}  (1998); 
  M. Peterson, in his paper \emph{Dante and the 3-sphere} \cite{Peterson-Physics} (1979),  
  and the same author some  30 years later, in his paper  \emph{The geometry of paradise} \cite{Paterson-Intelligencer} published in the Mathematical Intelligencer (2008); 
   R. Osserman, in his book  \emph{Poetry of the Universe: A Mathematical Exploration of the Cosmos} \cite{Osserman} (1995);  W.  Egginton, in his paper  \emph{On Dante, Hyperspheres, and the Curvature of the Medieval Cosmos} \cite{Egginton}  (1999);
M. Wertheim, in her book \emph{The pearly gates of cyberspace: A history of space from Dante to the internet} \cite{Wertheim} (1999);  G. Mazzotta,  in his article  \emph{Cosmology and the Kiss of Creation (Paradiso 27-29)} \cite{Mazzotta}  (2005); 
  S. L. Lipscomb in his book  \emph{Art Meets Mathematics in the Fourth Dimension} \cite{Lipscomb}  (2011); 
Curt McMullen in his survey, \emph{The evolution of geometric structures on 3-manifolds} \cite{McMullen-Dante}, published in the Bulletin of the American Mathematical Society (2011);\footnote{McMullen's survey starts  with the sentence:
 ``In 1300, Dante\index{Dante}  described a universe in which the concentric terraces of hell---nesting down to the center of the Earth---are mirrored by concentric celestial spheres,
rising and converging to a single luminous point. Topologically, this
finite yet unbounded space would today be described as a three-dimensional sphere."}
 M. Bersanelli, in his article \emph{From Dante's\index{Dante!\emph{Divine Comedy}}  universe to contemporary cosmology} \cite{Bersanelli}  (2016);
C. Rovelli, in his article \emph{Michelangelo's stone: an argument against Platonism in mathematics}, published in a philosophical journal 
\cite{Rovelli} (2017);
E. Ghys, in  his book \emph{A singular mathematical promenade} \cite{Ghys} (2017);
K. Sunada, in the article \emph{From Euclid to Riemann and beyond: how to describe the shape of the universe}, published in a book I edited  \cite{Sunada}  (2019). 

Let me also mention  the  interesting and beautifully illustrated article \cite{Dreyer} by J. L. E. Dreyer, \emph{The cosmology of Dante}, published\index{Dante!cosmology}  in 1921 in \emph{Nature},  in which the author surveys Dante's cosmological ideas expressed  in his poem (with no mention of the topology of the space described there).
The article ends with the words: 
``To the student of the history of science it is a never-failing source of pleasure to find medieval cosmology so beautifully illuminated in the writings of the great Florentine poet" \cite{Dreyer}.

I think that all these authors failed in two ways:
First of all, none of them refers to the work of Florensky\index{Florensky, Pavel Aleksandrovich} which was written between 40 years and almost a century earlier. For several among them, this is understandable because Florensky's booklet, until its first translation in French in 2016, existed only in Russian.\footnote{First edition, Moscow, Pomor'ye 1922 (1000  printed copies), new edition, Moscow, Lazur, 1991.} Let me mention incidentally that one rarely finds authors on this list who quote others on the same list who made equivalent statements. Among
the few authors who cited others are Rovelli, Sunada and Ghys, who quote Peterson.

The second way I think these authors failed is that all of them declared, omitting an important part of Dante's\index{Dante!\emph{Divine Comedy}}  poem which was highlighted by Florensky, that the space which is described by the Florentine poet is the 3-dimensional sphere. In fact, they missed the part on the non-orientability of the surface\index{non-orientable surface!Dante} on which they traveled.
The only paper I am aware of where the authors mention non-orientability is a commentary on Florensky's \emph{imaginaries} by Bayuk and Ford \cite{BF}, in which the mathematical phrasing is very shaky but where the authors reach the conclusion that ``[\ldots] All these considerations brought Florensky to the
conclusion that the geometry of Dante's cosmos is very similar to that of the single-sided elliptical surface described by Riemann and Klein."
  
Before developing this, I would like to point out that
 Florensky,\index{Florensky, Pavel Aleksandrovich} with his acute intellectual honesty, noted that three authors before him had  observed that the geometry of the \emph{Divine Comedy} is non-Euclidean, namely, the American mathematician George Bruce
Halstead (1853-1922), the German mathematician Heinrich Martin Weber (1842-1913), and the German historian of mathematics Maximilian Simon
(1844-1918).\footnote{Personally, I looked for written texts where these authors made such a statement, but I could not find any.} For the fact that Dante's world is non-Euclidean (which presumably means, in this setting, that it carries a geometry of positive curvature), Florensky refers 
to a passage in the \emph{Divine Comedy}\index{Dante!\emph{Divine Comedy}}  where the author speaks of ``triangles whose  angle sum is not equal to two right angles". Indeed, at some point of the exposition, Dante\index{Dante!\emph{Divine Comedy}}  quotes King Solomon who asks the Lord whether there exists a triangle inscribed in a semi-circle, with one side being the diameter, and no right angle  (\emph{Paradiso},  Canto XIII, v. 101-102):\footnote{Transl. Henry Wadsworth Longfellow (1807–1882).}  
\begin{verse}
non \emph{si est dare primum motum esse}, \hfill
nor \emph{si est dare primum motum esse},\footnote{``nor if we have to admit a prime motion" (Dante wrote this line in Latin.)}\\

o se del mezzo cerchio far si puote \hfill  	nor if in a semicircle a triangle can be formed \\

triangol s\`\i \  ch'un retto non avesse.\hfill   without its having one right angle. \\
 
\end{verse}

The question of non-orientability needs a longer quote and is based on Florensky's analysis, which I will review now.

Florensky\index{Florensky, Pavel Aleksandrovich} starts by recalling Dante's description of the Ptolemaic model of the universe, composed of spheres: the Earth, the celestial bodies and the spheres on which they stay, together with the Empyrean stratum, surrounding the  whole universe and occupied by the element of fire. 
 Dante\index{Dante!\emph{Divine Comedy}}  recounts a long imaginary journey he makes on his way  to visit God, who dwells in the last circle after the Empyrean. The path he takes leads successively through Hell, Purgatory, and Heaven. In this journey Dante is 
accompanied by Virgil, the Roman poet from the I${}^{\mathrm{st}}$ century {\sc bc}, author of the \emph{Aeneid}, who guides him through Hell and most of the Purgatory.  On their way,  the two poets start by going downwards, traversing one after the other the spirals of Hell.  It is important to note that during this descent, they always maintain a standing position, their heads pointing towards the point from which they started, that is, the city of Florence, and their feet directed towards the center of the Earth.
When the two poets find themselves at the bottom of Hell, they see Lucifer with legs up and head down. It is conceivable (and mathematically equivalent) that at this point they are themselves turned upside down, with their feet directed towards the surface of the Earth and their heads on the opposite side. Here is the passage of the\index{Dante!\emph{Divine Comedy}}  \emph{Divine Comedy} quoted by Florensky\index{Florensky, Pavel Aleksandrovich}  (\emph{Inferno}, Canto XXXIV, v. 76-96):

\begin{verse} \small\small
Quando noi fummo là dove la coscia  \hfill When we were come to where the thigh  revolves \\

si volge, a punto in sul grosso de l'anche,  \hfill  Exactly on the thickness of the haunch,\\

lo duca, con fatica e con angoscia,    \hfill  The Guide, with labour and with   hard-drawn  breath, \\   

\medskip                        

volse la testa ov'elli avea le zanche,  \hfill Turned round his head where he had had his legs, \\

e aggrappossi al pel com'om che sale, \hfill    And grappled to the hair, as one who mounts, \\

s\`\i \  che 'n inferno i' credea tornar anche.   \hfill    So that to Hell I thought we were  returning. \\    

\medskip                  

``Attienti ben, ché per cotali scale",\hfill    ``Keep fast thy hold, for by such stairs as   these," \\

disse 'l maestro, ansando com'uom lasso, \hfill   The Master said, panting as one  fatigued, \\

``conviensi dipartir da tanto male".  \hfill  ``Must we perforce depart from so much  evil." \\   

\medskip                              

Poi usc\`\i \  fuor per lo f\'oro d'un sasso, \hfill  Then through the opening of a rock he  issued,  \\     

e puose me in su l'orlo a sedere;  \hfill  And down upon the margin seated me; \\     

appresso porse a me l'accorto passo.  \hfill     Then tow'rds me he outstretched his  wary  step.  \\      

\medskip

Io levai li occhi e credetti vedere \hfill   I lifted up mine eyes and thought to see  \\     

Lucifero com'io l'avea lasciato,\hfill     Lucifer in the same way I had left him;  \\     

e vidili le gambe in sù tenere;       \hfill       And I beheld him upward hold his  legs.     \\

\medskip

e s'io divenni allora travagliato,\hfill   And if I then became disquieted,  \\     

la gente grossa il pensi, che non vede \hfill   Let stolid people think who do not see  \\     

qual è quel punto ch'io avea passato.   \hfill      What the point is beyond which I had passed.      \\       

\medskip

``Lèvati sù", disse 'l maestro, ``in piede: \hfill  ``Rise up," the Master said, ``upon thy  feet;  \\     

la via è lunga e 'l cammino è malvagio,  \hfill The way is long, and difficult the road,  \\     

e già il sole a mezza terza riede".  \hfill   And now the sun to middle-tierce returns."

\end{verse}

\begin{figure}[htbp]
\centering
\includegraphics[width=8cm]{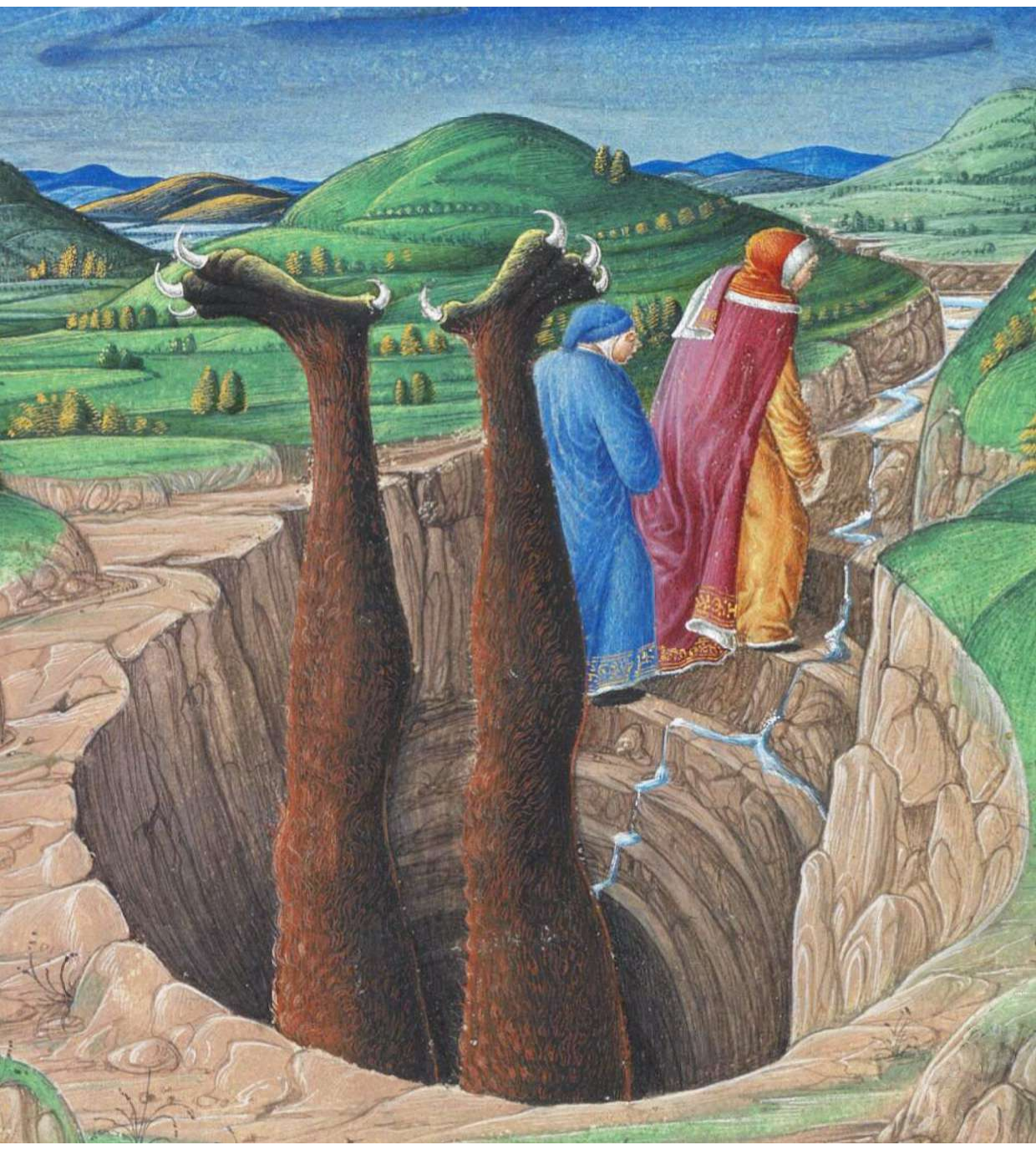}
\caption{Illustration of canto XXXIV of Dante's \emph{Inferno}, with Lucifer upside down, from a manuscript in Rome, Vatican Library, Ms Urb. Lat. 365, fol. 95v. Ferrara, 1474-1482.}
\label{fig:Ferrara}
\end{figure}

The two poets, after having avoided the center of the world, continue,\index{Florensky, Pavel Aleksandrovich} following the same way, but now they find themselves moving upwards, until they traverse a crater track, to find themselves in the clear world. This is illustrated here in Figure \ref{fig:Ferrara}, from a manuscript in the Vatican Library (Inferno, Canto XXXIV, v. 130-139): 

\begin{verse}\small

Lo duca e io per quel cammino ascoso \hfill The Guide and I into that hidden road \\

intrammo a ritornar nel chiaro mondo;    \hfill   Now entered, to return to the bright world; \\

e sanza cura aver d'alcun riposo,     \hfill     And without care of having any rest   \\    

\medskip                        

salimmo sù, el primo e io secondo,    \hfill  We mounted up, he first and I the second,  \\

tanto ch'i' vidi de le cose belle     \hfill  Till I beheld through a round aperture \\

che porta 'l ciel, per un pertugio tondo.     \hfill  Some of the beauteous things that Heaven 
 doth bear; \\

\medskip

E quindi uscimmo a riveder le stelle.   \hfill Thence we came forth to rebehold the  stars.  \hfill

\end{verse}

Beyond this limit, Dante,\index{Dante!\emph{Divine Comedy}} begins the ascent towards the celestial spheres, always following the same line, passing through the Empyrean where he is able to contemplate the Glory of God, and eventually ending up again in Florence, where he began his journey.
In other words, the poet, having traveled along a (non-Euclidean) straight line in the 3-dimensional world, finds himself at the same place and in the same position, but after having undergone a complete reversal.  The surface on which he has traveled is non-orientable.  Obviously, Dante did not have the technical means to express this fact in the way a mathematician would today (and, likewise, he did not have the technical means to say that the universe is a three-dimensional sphere).  Florensky commented on Dante in his own way, using the tools he had, and letting his imagination do the rest.

Florensky\index{Florensky, Pavel Aleksandrovich} calls Dante\index{Dante}  ``the greatest master of world design" \cite[p. 70]{Imaginaires}, and 
 the last sentence of the \emph{Imaginaries} reads:
 ``Thus, beyond time, the \emph{Divine Comedy} is not behind, but in front of contemporary science" \cite[p. 78]{Imaginaires}. I said earlier that all the modern authors I mentioned failed to understand Dante's cosmology,\index{Dante!cosmology}  and I explained why. In fact, Florensky himself considered that his predecessors failed to understand Dante's cosmology, for another reason. Let me quote him on this. Before he gives his own description to which I have just referred, he  writes, on p. 71 of \cite{Imaginaires}: 
 \begin{quote}\small
 The following image is usually given in commentaries on the \emph{Divine Comedy}. The terrestrial sphere is surrounded by the sphere of celestial bodies, by the sky of unchanging, immovable stars, by the crystal sky and finally by the Empyrean, and Dante's path from the bosom of the Earth follows a broken line spiraling through the concentric spheres and curving 180$^{\mathrm{o}}$ towards the zenith of Sion.  But this drawing corresponds neither to Dante's\index{Dante!\emph{Divine Comedy}}  narrative, nor to the foundations of cosmology. 
 \end{quote}

 I noticed long ago that Dante occupies a special place in Russian culture. Talking with some Russian colleagues and friends, I realized that some of them 
can recite passages of the \emph{Divine Comedy}; they learned it at school.
I cannot resist mentioning the first three lines of this poem, which my Russian colleague, Mitia Millionschikov, to whom I recently told about the article I was writing, spontaneously recited:

\begin{verse}
Nel mezzo del cammin di nostra vita, \hfill
Midway upon the journey of our life,
\\

mi ritrovai per una selva oscura \hfill 
I found myself within a forest dark,\\

ché la diritta via era smarrita.
\hfill   
For the straightforward pathway had been   lost.

\end{verse}

I will mention Andredi Tarkovsky later in this chapter.
The opening of the \emph{Divine comedy} is recited (twice) in his film \emph{The mirror}.

Let me end this section by mentioning that in his major theological work, \emph{The pillar and ground of the truth},\index{Florensky, Pavel Aleksandrovich!\emph{The pillar and ground of the truth}} which I will quote several times, Florensky\index{Florensky, Pavel Aleksandrovich} refers to another grand epic poem in which a long journey is recounted, the \emph{Odyssey}.

\section{Time, relativity and literature} \label{s:time}

In the \emph{Imaginaries} \cite{Imaginaires}, after the passage on Dante, comes a discussion of  Einstein's\index{relativity theory}  theories of special and general relativity, an important component of the book.  Florensky's goal is to show that, in the same way as ``space must be represented as an elliptic space and recognized as limited, time is limited and is closed in on itself." \cite[p. 73]{Imaginaires}

 Sabaneeff, who kept notes of his conversations with Florensky, recalls that the latter had the following thoughts: 
 \begin{quote}\small
 Time, like space, can
be not only straight but also curved (Einstein's time is curved).
It can be infinite, but it can also be closed (analogous to a
 closed line). The concept of a single future is nothing but a
theory of predestination. Furthermore, at every moment, the
future divides itself into a multitude (infinity) of actual new
futures. The ``closed" time, after having run its cycle, ``returns"
and is inevitably drawn into the ``sphere of eternal repetition,"
into a vicious circle, which he was inclined to identify with
the idea of ``hell" in religion.
 \end{quote}
 
I think that every sentence in this quotation could be analysed in the context of scientific and philosophical culture of those times, and of today. In fact, one could make relations with the idea of eternal return which is associated with  Pythagoreanism and Stoicism in Ancient Greece, then Proclus, Maximus the Confessor and many others. Closer to us, the idea of eternal return was revived by Nietzsche, late in his life, in his search for a way to combat nihilism.  

Florensky\index{Florensky, Pavel Aleksandrovich} understood early on that the theory of relativity has changed the worldview of scientists, and that in this context, Copernicus' heliocentric system was no more valid than the geocentric theory of the Ancients. More than that, he maintained the following (\cite[p. 75]{Imaginaires}): 
\begin{quote}\small
It would be a serious error to declare that the two systems, Copernicus' and Ptolemy's, are equally legitimate: they are so only in terms of abstract mechanics, but in terms of all the data, it is the latter that is true and the former that is false. This is a direct confirmation of the great poem, albeit more than six hundred years later. 
\end{quote}

Florensky is talking here about ``data". Ptolemy collected a huge amount of data for the drawing of  geographical maps, in his \emph{Almagest}, for celestial maps, and his \emph{Geography}, for terrestrial maps. It is true that Ptolemy's data, and the maps that were drawn using them,  are amazingly precise, despite the rudimentary means used for collecting the needed astronomical observations, and despite the fact that the theory is are based on a geocentric system.
Florensky's defence of Ptolemy's system, which was not been accepted by some of his contemporaries, was the cause of subsequent troubles with the authorities who accused him of obscurantism. These problems are reported and commented on in an interesting addendum to the French edition of  the \emph{Imaginaries} \cite[p. 95-118]{Imaginaires}.\footnote{\label{n:Postface}The Postface to the French edition of the book, written by Pavel Vassilievich Florensky,\index{Florensky, Pavel Vassilievich} the author's grandson, contains the translation of a letter sent by Pavel Aleksandrovich Florensky to the \emph{Political Department}, dated September 13, 1922.  The letter starts by the words: ``This is not to claim everyone's `right' to print whatever comes to his mind, and while recognizing that, in the name of goals concerning the whole state, the authorities cannot yield to the interests of literature, but I assume that in the case at hand, of my little brochure on  \emph{Imaginaries}, there is a \emph{quid pro quo}, which can be explained by the novelty of the object, which is why I will take the liberty of disturbing the censors in order to set out the following considerations. \ldots" The passage that was censored is the one about Dante,\index{Dante}  and Florensky recalls in his petition the context in which he wrote this passage (the 600th anniversary of the death of the Florentine poet) while defending the use in geometry of poetic images as the expression of a certain \emph{psychological factor}. He ends his letter with a reminder of the work he is conscientiously doing to serve the State. His request not to have the book censored was eventually accepted.}
Although at some point Florensky gives\index{Florensky, Pavel Aleksandrovich} a reasonable exposition of some elements of the theory of relativity, he introduces other thoughts which are part of the expression of his atypical view of science, which irritated the authorities, accusing him of anti-modernism, and which were termed irrational by readers of that work. 
At the end of his book, after his excursus on the \emph{Divine Comedy},  Florensky returning to the \emph{imaginaries}, writes:
\begin{quote}\small
The domain of the imaginary is real, conceivable, and in Dante's language it is called the {\bf Empyrean}.\footnote{The bold face letters are Florensky's.} We can represent all of space as {\bf double}, being formed of surfaces with Gaussian coordinates, real surfaces and imaginary ones corresponding to them, but the transition from the real to the imaginary surface is only possible through the {\bf fracture} of space and the {\bf turning} of the body inside itself. In the meantime, we only represent this process through the acceleration of motion, perhaps through the motion of certain small parts of the body, beyond the speed limit $c$, but we have no proof that other means are impossible.
\end{quote}

\begin{figure}
\begin{center} 
\includegraphics[width=10 cm]{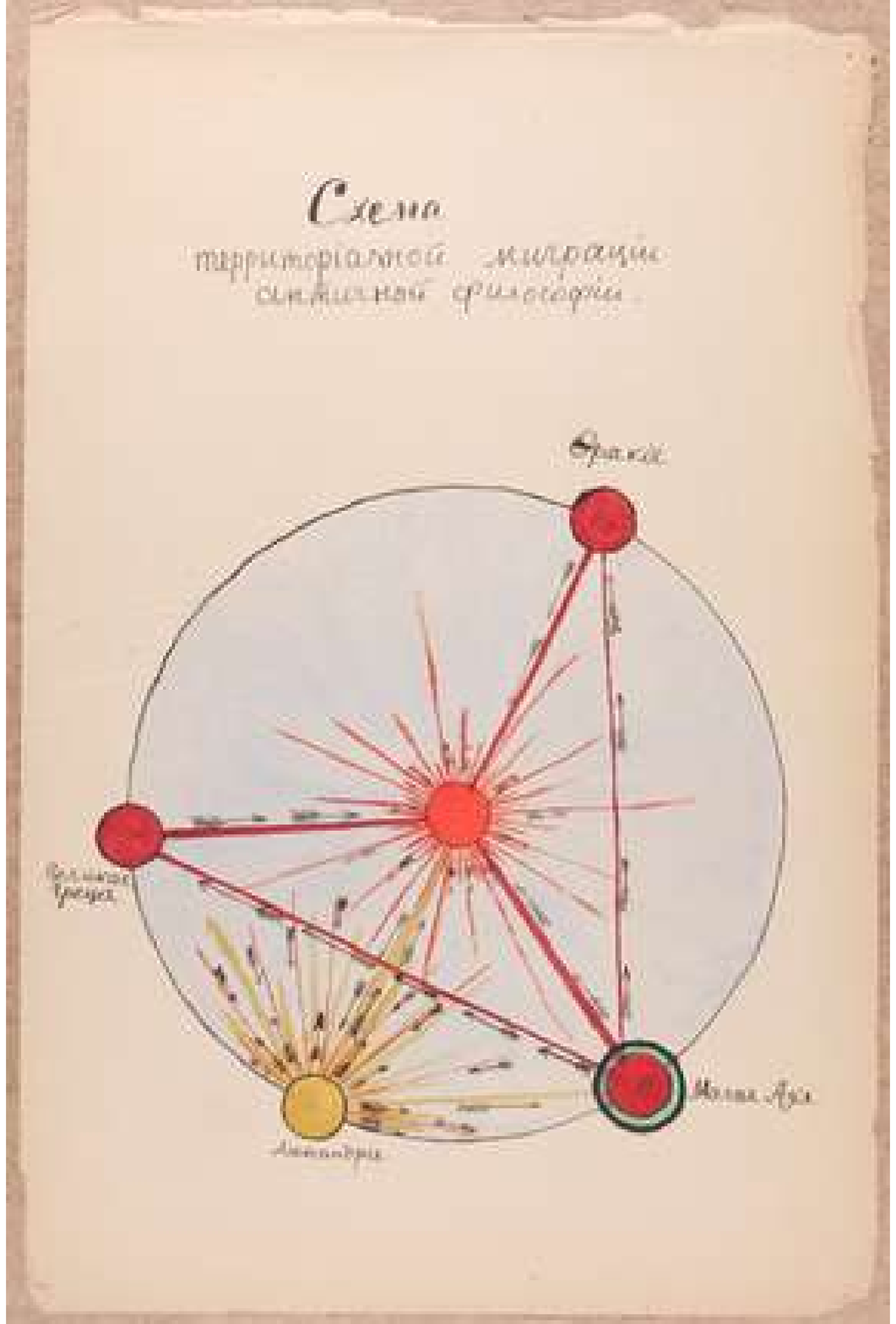} 
\caption{Illustration from Florensky's \emph{Imaginaries in geometry}.} 1922. Paper, ink
\label{fig:Imaginary}
\end{center}
\end{figure}

In this passage, the expressions \emph{fracture} and \emph{turning}, that Florensky emphasizes, refer to the discontinuous operation of orientation-change that occurs after one describes a complete geodesic in the projective plane. In Figure \ref{fig:Imaginary}, I have reproduced an illustration of Florensky's \emph{Imaginaries in geometry}, of a man walking in a non-orientable world. As to the mathematical notion of discontinuity, I already mentioned that it was already at the service of the philosophy of Bugaev,\index{Nikolai Vasilievich Bugaev} Florensky's first mentor who exerted a strong influence on him. I shall elaborate on this  in \S \ref{s:Bugaev}. Let me quote  here another passage of the \emph{Imaginaries} in which  Florensky\index{Florensky, Pavel Aleksandrovich} talks about discontinuity \cite[p. 76]{Imaginaires}:

\begin{quote}\small
 
What does the speed limit $c=3\times 10^{10}$ cm/second actually mean? It does not mean that velocities equal to or greater than $c$ are impossible, but it does mean that with them would appear entirely new conditions of life that we cannot yet visually represent, and, perhaps, forms transcendent to our Kantian terrestrial experience. But this does not at all mean that such conditions are impossible, and perhaps as the realm of experience expands, they will become representable. In other words, the life of a world with velocity equal to $c$ is qualitatively different from that observed with velocities below $c$, and \emph{the transition between the domains of this qualitative difference can only be thought of as discontinuous}.
\end{quote}

   In a manuscript on religious art, written in 1922 and published in  Russia in 1999 and in French translation a little bit later (see the historical note in \cite[p. 121]{Flo2}), Florensky expresses a similar idea: 
   \begin{quote}\small
 [\ldots] While everyone agrees, even without knowing the principle of relativity, that every system (at least in relation to the case under consideration) has its own time, its own rhythm and its own measure, very few people have thought about the possibility of time running at an {\bf infinite} speed, and even turning around completely, once this infinite speed has been exceeded, and {\bf reversing} its course. But time can really be instantaneous, flowing from the future to the past, from consequences to causes, teleologically speaking, and this is what happens when our lives move from the visible to the invisible, from the real to the imaginary.\footnote{Here, Florensky refers to the work of Carl Du Prel\index{Du Prel, Carl} (1839, 1899), who wrote about instantaneous time (without mentioning inverted time).}
   \end{quote}

   The last three quotes express similar ideas, in different contexts, poetry, science and the imaginary. The reversal of time\index{reversal of time}  is a counterpart to the idea of reversal of space,\index{reversal of space} expressed\index{reverse perspective} in the \emph{Reverse perspective}.\footnote{It may be useful to recall here that ``reverse perspective", or ``inverted perspective", is a technique of drawing, used in Byzantine and Russian icon painting, which is in opposition with the rules of perspective that are followed in Western classical art and theoreticized  by the Renaissance artists. In reverse perspective, objects depicted that are supposed to be further away from the viewer of the painting  are drawn as larger, whereas closer objects are drawn as smaller. 
In reverse perspective, parallel lines that move away of the viewer of the picture are drawn as diverging instead of converging lines, as they are drawn in the usual linear perspective drawing.}

       The \emph{Imaginaries} had an impact on the symbolist and Neo-Platonic trend in Russian philosophy. In the Postface, written by Florensky's\index{Florensky, Pavel Aleksandrovich} grandson (p. 104-105 of the French edition), there is a passage on the influence of this work on the philosopher Aleksei Fedorovich Losev\index{Losev, Aleksei Fedorovich} (1893-1988), who referred to the last section of the  \emph{Imaginaries} for the possibility of a motion at a velocity greater than that of light, where ``the body collapses into itself; its surface folds up within itself, turning inside out, through itself, and time begins to flow backwards, i.e., the consequence precedes the cause \ldots The cause is transformed into a goal, an ideal."

       I would like to end this section by a few words on the influence of the \emph{Imaginaries} on  literature, and for this I will quote another passage from the Postface of the French edition \cite[p. 105-106]{Imaginaires}. It concerns Mikhaïl Bulgavov, one of the best Russian writers of the twentieth century:\footnote{Mikhail Afanasyevich Bulgakov (1891-1940) is best known for his novel \emph{The master and margarita}, one of the important pieces of literature of the twentieth century. It was written under the Stalin regime. Bulgakov was a medical doctor. During his lifetime, the novel belonged to the underground literature in the Soviet Union. It was published posthumously, for the first time as a book, in Paris, in 1967.}
\begin{quote}\small
Mikhaïl Bulgakov's widow, Elena\index{Bulgakova, Elena Sergeyevna} Sergeyevna Bulgakova,\index{Bulgakov, Mikhail} 
testifies that the writer liked Florensky, and that  his library contained a copy of the \emph{Imaginaries}.\index{Bulgakov, Mikhail!\emph{The Master and Margarita}} This book was one of the writer's first purchases in Moscow, and he would hide it from onlookers. One of the most inspired and lyrical passages of the novel \emph{The White Guard}, written in 1923-1924, appears as a paraphrase of the concluding paragraph of \emph{The Imaginaries}. [\ldots] While working on his novel \emph{The Master and Margarita}, Bulgakov reread the \emph{Imaginaries} more than once.  His copy of the book contains underlined words and exclamation marks in the margins of the pages of paragraph 9. Elena Bulgakova recounts that when the first listeners of \emph{The Master and Margarita} heard, dumbfounded, the novel's finale, Bulgakov pointed out to them the pages of Florensky's book that scientifically confirmed his ideas. In Florensky's\index{Florensky, Pavel Aleksandrovich} mathematical and philosophical interpretation of Dante's\index{Dante}  journey with Virgil, Bulgakov saw the basis for the final chapters of his novel.
\end{quote}
 
We have been naturally led to literature, and I will say more about this in \S \ref{s:Silver-age}.  I would like to end this section by quoting again Sabaneeff, from his article \cite{Sabaneeff}, on the mathematical tendency of Florensky: 
 \begin{quote}\small
[Florensky was] an excellent mathematician who concentrated his attention
primarily on the latest trends in mathematics, he expressed
time and again the idea that mathematics was broader and
more comprehensive than the human mind. Mathematics, according
to him, ``is more intelligent than the human brain; it
leads man farther ahead into regions inaccessible to the brain
and not corroborated by sensual images." It was clear to me
that mathematics was his guide even in the area of mystic
speculations: it helped him not only through the elementary
language of numbers (as in the case of many earlier mystics)
but by means of the whole panoply of the modern mathematical
apparatus: analysis, the theory of sets, and all the latest theories
on the boundary between physics and mathematics.
\end{quote}
 
\section{Bugaev and Florensky: Mathematics and philosophy of discontinuity}\label{s:Bugaev}
 
 If one wishes to mention a single mathematician who exerted a definite influence on Pavel Florensky\index{Florensky, Pavel Aleksandrovich} during the years he was a student at the Imperial University of Moscow, Bugaev's\index{Nikolai Vasilievich Bugaev} name comes naturally. Nikolai Vasilievich Bugaev Bugaev taught mathematics at Moscow University from 1867 till his death. Besides being a founder and president of the Moscow Mathematical Society, he was also an active member of the Moscow Psychological Society. He is considered as the main founder of the Moscow school of function theory, about which I will say a few words in \S \ref{s:Luzin}. Some of his students became first-order mathematicians; these include Nikolay Yakovlevich Sonin\index{Sonin, Nikolay Yakovlevich} (1849-1915), Konstantin Alekseevich Andreev (1848-1921), and Dmitri Egorov,\index{Egorov, Dmitri Fedorovich} whom we already mentioned.  He  also exerted a strong influence on Boleslav Kornelievich Mlodzievsky, who was one of Florensky's teachers and whom we already mentioned as well. 
 
 Besides being a mathematician, Bugaev\index{Nikolai Vasilievich Bugaev} was a philosopher steeped in humanism and Christian philosophy. For him, mathematics was a means of expressing philosophical ideas.  Notions like continuity and discontinuity, which formed the basis of the theory of functions that was at the heart of his mathematical investigations,  were for him, at the same time, part of the language of  metaphysics,\index{continuity (philosophical notion)} used in theories concerning free will, continuity in history,\index{discontinuity (philosophical notion)} social revolutions, etc. It was through his influence that Florensky made some of these mathematical notions the basis of his philosophical thought. Sabaneef writes about Florensky, in the article we mentioned \cite[p. 321]{Sabaneeff}: ``It seems to me that the dominant motive of his life was the
idea of the alliance and fusion of science and revelation, the
termination of the antagonism between these two spheres that
had developed in the course of the historical process."

To give an idea of Bugaev's\index{Nikolai Vasilievich Bugaev} world and his influence on Florensky,\index{Florensky, Pavel Aleksandrovich} I will make a short overview of the talk he delivered at the 1897 International Congress of Mathematicians held at  Zurich.\footnote{Initially, the first International Congress of Mathematicians was scheduled for the year 1900, and the Zurich congress of 1897 was planned as a trial congress. But eventually  the latter was so successful that it was considered later as the first ICM.} The  talk was naturally in French, and is titled \emph{Les mathématiques et la conception du monde au point de vue philosophie scientifique} (Mathematics and the conception of the world from the point of view of scientific philosophy) \cite{Bugaev-ICM}. The aim of this talk is to develop the idea that some fundamental principles of mathematics are at the basis of important philosophical conceptions of the world. Thus, the subject of the talk involves science, philosophy and the description of the world.\footnote{\label{f:Poincare} This kind of philosophical topic, chosen as the subject of a lecture delivered at an ICM, was not exceptional. Poincaré, in his lecture to the same congress, titled \emph{On the relationship between pure analysis and mathematical physics}, makes a comparison between mathematics cultivated for its own sake and that which has applications and is used as a language for physics.  He discusses there analogies in mathematics, the reason why mathematics needs other sciences to develop, questions of rigor, and, like Bugaev,\index{Nikolai Vasilievich Bugaev} issues of continuity and discontinuity.} 
 The mathematical  concepts that are referred to are those that accompanied Florensky in his later philosophical writings:  discontinuity,\index{continuity (philosophical notion)}  contradiction,\index{contradiction (philosophical notion)} the infinite,\index{infinity (philosophical notion)} antinomy,\index{antinomy (philosophical notion)} etc.

Bugaev\index{Nikolai Vasilievich Bugaev} declares in his talk:\footnote{My translation from the French.} 
\begin{quote}\small
It is very important to examine presently the influence of these parts of mathematics on the scientific and philosophical conception of the world. [\ldots] At present, this need for number and measurement seems to be the vital issue of the day, not only in science, but also in art and human relations. Finding the measure in the realm of thought, will and feeling is the problem of the contemporary philosopher, politician and artist.

\end{quote}
 Bugaev,\index{Nikolai Vasilievich Bugaev} in his speech, wants to show how man, with the help of ``number and measure", can get out of the realm of instinct. He considers that what is needed for that purpose is the theory of functions.  He makes a distinction between \emph{analysis}, the theory of continuous functions, and \emph{arithmology},\index{arithmology} the theory of discontinuous ones.\footnote{In his article at the same congress, Poincaré also talks about the ideas of continuity and discontinuity. In particular, he says: ``The only natural object of mathematical thought is the [natural] number. It was the external world that imposed the continuous on us, which we undoubtedly invented, but which it forced us to invent. Without it, there would be no infinitesimal analysis; all of mathematical science would be reduced to arithmetic or the theory of substitutions."}  He states that mathematical analysis is distinguished by its generality and universality while arithmology bears the stamp of an original individuality, and is attractive because of its mysterious character and striking beauty. He believes that this explains the fact that some thinkers have linked various questions of mystical philosophy to questions related to the integers. He adds that the elegant truths of arithmology awaken in the scientist a sense of scientific beauty that fully satisfies him. He then talks about the applications of mathematical analysis in philosophy.  Continuous functions explain mathematically some natural phenomena. It is interesting that he quotes a poem by the poet Sichler (translated from the Russian), expressing the analytic conception of the world. Having a closer look at natural phenomena, he concludes that the point of view of discontinuity is necessary, in particular in atomic physics, chemistry, biology and sociology. ``Discontinuity always appears,\index{discontinuity (philosophical notion)} he says, where independent and autonomous individuality manifests itself," and further: ``The highest human interests are often linked to problems of arithmology. The arithmological worldview indicates that final causes also play a role in the world's manifestations."  Bugaev\index{Nikolai Vasilievich Bugaev} refers to Leibniz, as a source of inspiration for such questions.
 He explains why until present times only the analytical point of view was used in the scientific and philosophical considerations on the world and nature. The problem, he says, is that until now it was assumed that any scientific question assumes only one answer. But in arithmology, one finds functions that are multiple-valued, and sometimes these values are infinite. Such functions, he says, arise in nature, and he gives examples from sensation theory. He is led then to talk about antinomy.\index{antinomy (philosophical notion)} From the mathematical point of view, it is impossible to reconcile the concepts of continuity\index{continuity (philosophical notion)} and discontinuity.\index{discontinuity (philosophical notion)} But to have a complete understanding of facts of nature, we have to take equally into account the two notions. He says that in logic, biology, psychology, history of philosophy and sociology, the universal and the individual, the concrete and the abstract, and person and society mutually complete each other. These thoughts are  very much within the cosmic Russian stream of that epoch, according to which man and the cosmos are in harmony; this is why they called man a microcosmos, a terminology that was used in  Ancient Greece by Anaximander, Plato, the Stoics, and others.

  The use of the idea of discontinuity\index{discontinuity (philosophical notion)}  and its necessity in the description of nature, in Bugaev's\index{Nikolai Vasilievich Bugaev} teaching,  did not escape Florensky.\index{Florensky, Pavel Aleksandrovich}
The title of the dissertation he presented at the Faculty of Mathematics and Physics in 1904, first under the supervision of Bugaev and, after the latter's death, under that of Egorov,\index{Egorov, Dmitri Fedorovich}  is informative: \emph{The idea of discontinuity as an element of a worldview}. He quotes this work several  times in his later writings; see in particulier Note 183 of the \emph{Pillar} \cite[p. 485] {Flo1}, where he mentions the numerous existing works ``in the domain of the formal investigation of the idea of discontinuity, i.e., in mathematics and in logic",  referring to the bibliography
 given in the above manuscript. He then 
  refers to the ``works of the Moscow mathematical school", and in particular to Bugaev's \emph{Introduction to the theory of numbers}  and his \emph{Mathematics and the scientific-philosophical worldview} (1898) \cite{Bugaev-world}.
%
  In the same note, Florensky talks about the idea of 
 discontinuity in thermodynamics,
leading to the construction of the discontinuous thermodynamic surface, and in the works of Gibbs on phases, of Korzhinsky on heterogenesis and of neo-Lamarckian and neo-vitalist researchers in biology, ``opening broad horizons of the discontinuous\index{discontinuity (philosophical notion)} evolution of organisms and their discontinuous adaptation." He also mentions researches in psychology and on the subconscious
and superconscious psychic life, where one analyses discontinuous changes in consciousness, in
 creativity, in inspiration, etc.       
   In a note to Letter X of his \emph{The pillar and ground of the truth} titled \emph{Sophia}, Florensky refers again to two papers by Bugaev,\index{Nikolai Vasilievich Bugaev} \emph{Geometry of arbitrary quantities} and \emph{Discontinuous geometry}. 
   
   In \S \ref{s:Vita}, I mentioned a passage from Florensky's\index{Florensky, Pavel Aleksandrovich} autobiographical notes \cite{Florensky-bio} where, talking about his 1899 crisis, he declares that  the notion of continuity\index{continuity (philosophical notion)} ``suffocated him". This shows the importance of this topic for him, and at the same time,  it shows how he experienced mathematics physically,  in the same way he experienced sound---we shall talk about that at the \S \ref{s:Silver-age}.

With Florensky, the notion of discontinuity entered the realm of theology, which is the main subject of \emph{The pillar and ground of the truth}.\index{Florensky, Pavel Aleksandrovich!\emph{The pillar and ground of the truth}}   In Letter V, titled \emph{The Comforter}, Florensky  talks about ``the discontinuity of the metahistorical revelation of the Spirit" \cite[p. 101]{Flo1}. 
 In Letter VI, titled \emph{Contradiction},\index{contradiction (philosophical notion)}  talking about the Apostle Paul, Florensky\index{Florensky, Pavel Aleksandrovich} writes: ``His brilliant religious dialectic consists of a
series of discontinuities;\index{discontinuity (philosophical notion)} it jumps from one affirmation to another, where each successive affirmation is antinomic with respect to the preceding one. Sometimes an antinomy\index{antinomy (philosophical notion)}  is even embodied in a stylistic discontinuity of exposition, in an external asyndeton.\footnote{By the word \emph{asyndeton}, Florensky refers to a form of speech in which coupling conjunctions are dropped or replaced by commas, in order to speed up speech. An example, in the writings of Saint Paul to whom  Florensky refers, given by F. Godet in \cite{Godet}, is the following passage (\emph{2nd Epistle to the Romans}, v. 9 and 10):  ``Tribulation and anguish, upon every soul of man that doeth evil, of the Jew first, and also of the Gentile; but glory, honour, and peace, to every man that worketh good, to the Jew first, and also to the Gentile."  Godet comments: ``As always, the asyndeton indicates a more energetic reaffirmation of the preceding idea." The notion  of \emph{asyndeton} is present at several places in the explanations that Godet gives in his book  \cite{Godet}.} Rationally contradictory and mutually exclusive judgments have their sharp edges directed against each other." Antinomy and discontinuity, are present in all of Florensky's theological works; I will say more about this in \S \ref{s:Contradiction}.
Regarding the philosophical implications of the concept of discontinuity\index{discontinuity (philosophical notion)} in Florensky's works, and, more generally, in Russian philosophical thought, let me mention the article \cite{Oppo} by Andrea Oppo, titled \emph{Conceptualising discontinuity: Pavel Florenskii's \emph{Preryvnost'} as a universal paradigm of knowledge}, in which the author discusses, among other things, Florensky's interpretation of Shakespeare's \emph{Hamlet} in his essay \cite{Florensky-Hamlet}, and Florensky's passage on Dante in the \emph{Imaginaries} \cite{Imaginaires}.

\index{Manin, Yuri Ivanovich!\emph{Mathematics as a metaphor}}
To conclude this section on Bugaev,\index{Nikolai Vasilievich Bugaev} and talking again about the ideas of continuity\index{continuity (philosophical notion)} and discontinuity\index{discontinuity (philosophical notion)} beyond mathematics, let me mention a paper by Yuri Manin  titled\index{Manin, Yuri Ivanovich} \emph{Time and periodicity from Ptolemy to Schr\"odinger: paradigm shifts vs. continuity in history of mathematics}, in which the author discusses the idea of ``scientific revolution", central in Thomas Kuhn's philosophy, opting rather to a  theory of continuity for what regards the history of pure mathematics. In the same article, Manin mentions Ptolemy's system, also dear to\index{Florensky, Pavel Aleksandrovich} Florensky.\footnote{For a related discussion, the interested reader is also referred to the article on Manin in the present Handbook \cite{Manin}. For a recent point of view on (the philosophy of) continuity and discontinuity in mathematics and in the relation between mathematics and physics, let me mention Plotnitsky's article \emph{Continuity and discontinuity in, and between, mathematics and physics} \cite{Plot}. Let me also take this opportunity to mention another paper by Plotnitsky, titled \emph{The Image of time and the specter of the clock, from modernist literature to cinema} \cite{Plot-Cinema}, in which he discusses continuity and discontinuity in Belyi's \index{Belyi, Andrei} famous symbolist novel \emph{Petersburg} (1913), and also in
cinema, and in particular in the work of Andrei Tarkovsky,\index{Tarkovsky, Andrei} who was very close spiritually speaking and also in his films, to Florensky, as we shall recall in \S \ref{s:conclusion}. Plotnitsky tells me that Belyi's novel \emph{Petersburg} is not only a founding novel of symbolism and of modernism in general, but it is also a novel about anarchism. In fact, the action takes place during the 1905 revolution. Plotnitsky points out the following passage from this novel which makes a relation with cinema, which in those times was one of the rising arts, seen by some (Lenin was among them) as the most important communist art: 
``We [our brains] are the unrolling of a cinematographic film subject to the minute action of occult forces; should the film stop, we will be fixed forever in an artificial pose of terror."}

 We have  mentioned Bugaev's\index{Nikolai Vasilievich Bugaev} influence on Florensky. If, in another direction, we want to name a mathematician on whom Florensky had a definite influence, the first name that comes to mind is Luzin, of whom I will talk in the next section.\index{Luzin, Nikolai Nikolaevich}

\section{Luzin: Spiritual crisis, function theory, Luzitaniya and disgrace} \label{s:Luzin}

 Nikolai Nikolaevich Luzin\index{Luzin, Nikolai Nikolaevich} and his mentor, Dmitri Fedorovich Egorov,\index{Egorov, Dmitri Fedorovich} 
are usually considered as the founders of the Moscow school of functions of a real variable. Demidov, in his paper \emph{On an early history of the Moscow school of theory of functions} \cite{Demidov1999}, writes that this school was born with Egorov's \emph{Comptes Rendus} note \emph{Sur les suites des fonctions mesurables} (On the sequences of measurable functions) \cite{Eg1}.  In this paper, Egorov\index{Egorov, Dmitri Fedorovich} proves that any convergent sequence of measurable functions is uniformly convergent up to neglecting a set of measure $\epsilon$, where $\epsilon$ may be chosen arbitrarily. The theorem was stated without proof by Lebesgue\index{Lebesgue, Henri} in his paper \cite{Lebesgue1} (1903) and in his book \cite{Lebesgue2} (1904).\footnote{Egorov attended Lebesgue's course at the Collège de France in 1903, see \cite{Kuznetsov}.}

The development of the theory of functions  of a real variable relies in an essential way on set theory,\index{set theory (historical development)}  and more precisely, on the field known as ``descriptive set theory,"\index{descriptive set theory} whose goal was to set up firm foundations for the topological and set-theoretical methods used in the theory of functions.  Descriptive set theory may be considered as a sub-area of mathematical logic, to a large extent based on the works of Georg Cantor\index{Cantor, Georg} (1845-1918), whose name is closely associated with the use of infinity in mathematics and who, incidentally, was Russian-born.\footnote{See the article \cite{Sinkevich-ICCM} for a concise biography of Cantor, together with notes on Florensky,\index{Florensky, Pavel Aleksandrovich} Mlodzievsky and other actors of the Silver Age.} Luzin\index{Luzin, Nikolai Nikolaevich}  was one of the architects of this emerging area, and an extensive development of this theory is contained in his 1930 book \emph{Leçons sur les ensembles analytiques et leurs applications} (Lectures on analytic sets and their applications) \cite{Luzin-Lecons}. The French mathematicians that were thoroughly involved in this field include \'Emile Borel (1871-1956), René Baire\index{Baire, René} (1874-1932), Jacques Hadamard,\index{Hadamard, Jacques} (1865-1963) Henri Lebesgue\index{Lebesgue, Henri} (1875-1941), Arnaud Denjoy (1884-1974) and Paul Painlev\'e (1863-1933). Florensky was a mathematics student at the time of the development of the foundations of function theory, with set theory in the background, the  theory which remained for him a constant subject of reflection.

  Luzin\index{Luzin, Nikolai Nikolaevich}  became familiar with the main problems of set theory and measure theory during a stay he made in Paris after the first Russian revolution and the bloody events in its aftermath (1905-1906). In Paris,  he attended lectures on function theory by the leading mathematicians there. Functions were classified in terms of set-theoretic properties and limit operations. New sets emerged and were given names. Since open and closed sets are not stable under simple operations involving infinity, $G_{\delta}$- and $F_{\delta}$-\emph{sets} were introduced as countable unions (respectively intersections) of open (respectively closed) sets. \emph{Borel sets} are obtained from open or closed sets under the operations of countable union, countable intersection, and taking complement. Among the other classes of sets that appeared, we mention \emph{A-sets}\index{analytic set}\index{analytic set} (A stands for ``analytic"), later called \emph{Suslin sets},\index{Suslin set} whose existence arose with a remark that Suslin\footnote{Mikhail Yakovlevich Suslin\index{Suslin, Mikhail Yakovlevich} (1894-1919) was Luzin's student. He died at age 25 from typhus, an epidemic that followed the Russian Civil War.} made, saying that a continuous image of a Borel set is not necessarily Borel, correcting a  mistake made by Lebesgue\index{Lebesgue, Henri} on this metter.\footnote{The story is recounted by Lebesgue in his preface to Luzin's monograph \emph{Leçons sur les ensembles analytiques et leurs applications} \cite{Luzin-Lecons}.} Suslin\index{Florensky, Pavel Aleksandrovich}  introduced the class of \emph{analytic sets}\index{analytic set} in his note \cite{Suslin}. Luzin,\index{Luzin, Nikolai Nikolaevich} who gave Suslin the impetus for the introduction of these sets, continued the latter's work on this topic after his death.   \emph{Baire functions}\index{Baire function} arose, as functions obtained from continuous functions by  a transfinite repetition of simple limits. Baire\index{Baire, René} first and second category sets\footnote{A \emph{Baire first category set}\index{Baire first category set} is a subset of a Baire topological space contained in a countable union of closed subsets which all have empty interior. A \emph{Baire second category set}\index{Baire second category set} is a subset of a Baire topological space which is not a Baire first category set.}  and \emph{Luzin sets}\footnote{A \emph{Luzin set}\index{Luzin set} is a subset of the real numbers which is uncountable but whose intersection with every set of the first Baire category is countable.} were introduced. Luzin\index{Luzin, Nikolai Nikolaevich} introduced the notion of \emph{projective set}, that is, a set obtained from a Borel set in $n$-dimensional Euclidean space  by operations of projection and taking complement. The notion of \emph{effective sets}\index{effective set} (a set whose definition does not use the axiom of choice) was born. Lebesgue\index{Lebesgue, Henri} introduced the notion of \emph{nameable set},\index{nameable set} a set that can be ``named". Luzin\index{Luzin, Nikolai Nikolaevich}  presents this notion, in his \emph{Mémoire sur les ensembles analytiques et projectifs} (1926) \cite{Luzin-memoire}. He writes : \begin{quote}\small
  One of the fundamental ideas that we owe to Mr. Lebesgue is the precise distinction between the two notions of a nameable set and   unnameable set.\index{unnameable set}  We already know what, according to Mr. Lebesgue,\index{Lebesgue, Henri} is meant by a nameable set: it is the set of points that can be named, i.e., characterized without possible confusion, by means of a suitable definition. ``Nameable" sets\index{nameable set}  are precisely those which, according to Mr. Lebesgue, are the necessary object of Mathematical Analysis.  As we know, in the current state of Science, unnameable sets are those that can be obtained by means of Mr. Zermelo's Axiom or similar reasoning. It is these unnameable sets that Mr. Lebesgue\index{Lebesgue, Henri} considers doubtful: they can only be introduced into Science under the necessary condition that they show their interest, applicability, usefulness and efficiency in solving some classical problem. But, in general, it is the domain of nameable beings\index{nameable being} that is, par excellence, the content of Mathematical Analysis. This is the conception of Mathematical Analysis proposed by Mr. Henri Lebesgue.
  \end{quote}

Luzin\index{Luzin, Nikolai Nikolaevich}  was aware of the difficulties which affected the theory of functions that arose from  set theory, including the use of transfinite constructions and the continuum hypothesis, the discovery of hierarchies between infinities, etc. 
Lebesgue\index{Lebesgue, Henri} refused the use of the axiom of choice, considering it as counter-intuitive.  He considered that in order for a mathematical object to exist, one has to able to explicitly ``name" a property that defines it in a unique way. Bourbaki in \cite[p. 53ff]{Bourba} has an interesting passage on this matter:
\begin{quote}\small

For Lebesgue, who broadens the debate, it all boils down to knowing what we mean when we say that a mathematical being ``exists''; for him, we have to ``name'' explicitly a property that defines it in a unique way (a ``functional'' property, we would say), and already the question of choosing a single element in a set seems to him to raise difficulties: we need to be sure that such an element ``exists'', i.e., that we can ``name'' at least one of the elements in the set. Can we then speak of the ``existence'' of a set in which we do not know how to ``name'' each element? Already, Baire has no hesitation in denying the ``existence" of the set of parts of a given infinite set; in vain does Hadamard observe that these requirements result in abandoning even the idea of talking about the set of real numbers: it is indeed this conclusion that \'E. Borel eventually rallies to. With the exception of the fact that the countable seems to have acquired the right to exist, we have pretty much returned to the classical position of the opponents of ``actual infinity".\footnote{The clear distinction between actual and potential infinity was made by Aristotle. According to Aristotle, actual infinity cannot exist because it is  paradoxical. Only potential infinity exists.} 
\end{quote}

On the questions of existence and non-existence of certain mathematical objects, Lebesgue,\index{Lebesgue, Henri} in the preface he wrote to Luzin's book, \emph{Leçons sur les ensembles analytiques et leurs applications} \cite{Luzin-Lecons}, recounts the following story:  ``When I was a student, we used to have coffee and talk about general ideas; the discussion was getting hot and seemed bound to go on forever, when one of us shouted: `First of all, you, do you exist? I say \emph{you} for convenience, but even, do I exist?'. \ldots"

 In France, like in Russia, a group of religiously inclined  mathematicians\footnote{In Russia, these mathematicians include Bugaev,\index{Nikolai Vasilievich Bugaev}  Luzin and Florensky.\index{Florensky, Pavel Aleksandrovich} In France, the religious tendency among mathematicians is represented by Cauchy, Hermite and Appell.} raised questions such as how, and to what extent, infinity---which\index{infinity (philosophical notion)} is an attribute of God---can be used in human and mathematical contexts.  It was a period where mathematicians realized that  the theory of transfinite numbers and the use of the axiom of choice may lead to contradictions. Lebesgue did not exclude the existence of sets that are neither finite nor infinite.\footnote{\label{Leb} Lebesgue writes in a letter to  Borel:  ``Although I very much doubt that we will ever name a set that is neither finite nor infinite, the impossibility of such a set does not seem to me to have been proved," see \cite{Cinq}. The same sentence is quoted in Luzin's ICM talk \cite{Luzin-ICM}.} This indicates again the fact that the notion of infinity, in the mathematical context, must be used with care.  It was natural that Florensky, even after he quitted the Mathematics Department of Moscow University, given his interest in all areas of science, and especially mathematics,\index{set theory (historical development)} kept abreast of all these developments which involve at the same time science, philosophy and religion.   
 
 Luzin\index{Luzin, Nikolai Nikolaevich}  was also tormented by philosophical and religious thoughts.
At the time he was a mathematics student,\index{Luzin, Nikolai Nikolaevich},  he underwent, like Florensky a few years before him, a spiritual crisis that shook what he called his ``materialistic worldview,"  an expression which, at that time, meant atheism. Eventually, his relation with Florensky\index{Florensky, Pavel Aleksandrovich} became an important factor in the resolution of this crisis and his return to mathematics. C. E. Ford, in his article  \emph{The influence of P. A. Florensky on N. N. Luzin}, writes \cite[p. 332-333]{Ford}: ``The early history of the Moscow school has been greatly illuminated by the recent discovery that a close relationship
had developed between Luzin\index{Luzin, Nikolai Nikolaevich} and Florensky \ldots This enabled Luzin to launch his career as a mathematician". The article \cite{Ford} is based on the correspondence between Luzin and Florensky which has been preserved, and I would like to quote some passages from it. This will confirm the fact that Florensky had a definite influence on one of the main founders, if not \emph{the} main founder, of the 20th century Moscow mathematical school.

In a letter to Florensky sent from Paris, dated May 1st, 1906, quoted in \cite[p. 335]{Ford}, and translated from \cite{Demidov1}, Luzin writes: 
\begin{quote}\small
Life
is too depressing for me, sometimes agonizingly depressing. I am left with nothing, no solid
worldview; I am unable to find a solution to the ``problem of life." My self-image is so frequently
changing that life has become pure torment. \ldots You would not believe how precious your letter and your concern for my inner life were
to me. It is horrible, horrible, infinitely horrible to feel yourself surrounded by egoism, nothing
but unrelieved egoism. \ldots By now you will have understood how precious to me, against the background of the vast
and empty night that oppresses the soul, the heart, and the brain to madness, how precious
to me is your concern for my inner life. \ldots  I cannot be satisfied any more with analytic functions and Taylor series. \ldots  To see the misery of people, to see the torment of life, to
wend my way home from a mathematical meeting, to wend my way through the Aleksandrovsky
Garden, where, shivering in the cold, some
women stand waiting in vain for dinner purchased with horror—this is an unbearable sight.
It is unbearable, having seen this, to calmly study (in fact to enjoy) science. \ldots Then came the terrible days for us all
[the 1905 revolution], for me in particular. \ldots  I could
not work in science, and I seem to have begun to lose my mind as a result of the impossibility
of living quietly and understanding where, where the truth is.  \ldots Yes, I now understand that ``science," in essence, is metaphysical and based on nothing.
The scientific quest for ``differential resolvents for equations of degree six" is absolutely
incomprehensible. \ldots At the moment my scholarly interests are in principles, symbolic logic,
and set theory.\index{set theory (historical development)}  But I cannot live by science alone. \ldots I have nothing, no worldview, and no
education. I am absolutely ignorant of the philological sciences, history, philosophy. \ldots It is
painfully clear to me that it makes no sense to settle for a ``scientific" education. 
\end{quote}

  Luzin's\index{Luzin, Nikolai Nikolaevich} crisis lasted three years and was eventually resolved with the decisive help of\index{Florensky, Pavel Aleksandrovich}  Florensky,\index{Florensky, Pavel Aleksandrovich} who persuaded Luzin to return to mathematics \cite{Ford}. Upon his return to Russia, Luzin frequently visited the Sergiyev Posad monastery and valued his meetings with Florensky. In a letter to his wife, dated June 29, 1908, talking about Florensky's \emph{Pillar},\index{Florensky, Pavel Aleksandrovich!\emph{The pillar and ground of the truth}} of which he read a preliminary version, he writes \cite[p. 337]{Ford}: 
  \begin{quote}\small
  I read it all at once in a single day---skipping a lot, but the impression was
overwhelming. As I read it, I was STUNNED the entire time by blows from a battering ram
against a stronghold.\ldots Our life, the life of the XIX and XX centuries, is full of confusion, full of contradictions.\index{contradiction (philosophical notion)} The work of Florensky is one solution to this impossible situation.
\end{quote}

  Luzin\index{Luzin, Nikolai Nikolaevich} was interested in all of Florensky's\index{Florensky, Pavel Aleksandrovich} writings. In a letter dated June 19, 1919, he asks him: ``I
would like to have, if possible, all your published works. I do have many, but do not have
your big works, the latest lectures on philosophy (it is said that they are already published)
and above all, I do not have your thesis \emph{The Pillar and Foundation of Truth.}"

   On his return to Moscow,  Luzin\index{Luzin, Nikolai Nikolaevich} joined the Faculty of Mathematics of Moscow University. He was appointed professor there in 1917, just before the revolution.

 A few words on Luzin's\index{Luzin, Nikolai Nikolaevich} mathematical interests are in order. Luzin was strongly influenced by French mathematics. In 1911 and 1912, Egorov\index{Egorov, Dmitri Fedorovich} and Luzin\index{Luzin, Nikolai Nikolaevich} published a series of papers in the \emph{Comptes Rendus} which confirmed the French school's support for the Russian \cite{Eg1, Lu1, Lu2}. In particular, Egorov's result in \cite{Eg1}  on sequences of measurable functions was a substantial improvement  of previous works by Borel and Lebesgue.\index{Lebesgue, Henri} This result says that any sequence of measurable functions is uniformly convergent up to a set of arbitrarily small measure. 
  Luzin's\index{Luzin, Nikolai Nikolaevich} result in his paper \cite{Lu2}, titled \emph{Sur un problème de M. Baire} (On a problem of Mr. Baire), answers a question raised by Baire.

  Luzin was an atypical teacher.
 Mikhail Alekseyevich Lavrentieff\index{Lavrentieff, Mikhail Alekseyevich}, who was his student, in a memorial article on his teacher  \cite{L2}, writes that around the years 1922-1926, Luzin introduced to the mathematical division of Moscow university a new style of unprepared lectures, far from the standard ``good lecture," and that the result was a brilliant, deep and fascinating way of teaching, encouraging the development of originality and independent thinking. Lavrentieff also writes that  the essential characteristic of the school founded by Luzin\index{Luzin, Nikolai Nikolaevich}  was ``the fostering of independent thought, the capacity to crack problems, to find new methods and to pose new problems".  He reports that Luzin  prepared his lectures only in outline, that  he was often late to his class and that he believed that having a strict schedule was not compatible with doing mathematics.

Luzin's\index{Luzin, Nikolai Nikolaevich} 1926 \emph{Mémoire sur les ensembles analytiques et projectifs}  (Memoir on analytic and projective sets)  \cite{Luzin-memoire} starts with the words: 
\begin{quote}\small
The current state of mathematical analysis provides ample proof of the desirability of establishing a precise demarcation between mathematical beings regarded as ``existing" and others whose reality is only apparent. \ldots This distinction was the subject of the famous \emph{Five letters on set theory}\footnote{Luzin refers to the letters reproduced in \cite{Cinq}.}  by\index{Hadamard, Jacques} M. Hadamard, Baire,\index{Baire, René} Lebesgue\index{Lebesgue, Henri} and Borel, on whose necessity Mr. \'Emile Borel insisted with extreme exactitude in his later writings.
This distinction between mathematical beings that are truly real and others that appear to be real, but which have no substratum and to which no intuition corresponds, is currently very necessary, but its imminent realization seems unlikely in the present state of Science.
\end{quote}

In 1927, Luzin\index{Luzin, Nikolai Nikolaevich} was elected corresponding member of the USSR Academy of Sciences. The following year, he was elected full member, but, surprisingly, at the Department of Philosophy. In 1928, he was Vice-President of the ICM (Bologna). He presented there a\index{set theory (historical development)}  communication titled \emph{Sur les voies de la th\'eorie des ensembles} (On the paths of set theory) \cite{Luzin-ICM} 
 in which he discussed the problems that arose in the theory of functions as a consequence of the questions that were raised in set theory regarding infinity and the definition of real number.  He used in this communication the philosophical tone that Bugaev had, 31 years before him, in his ICM communication.\index{Nikolai Vasilievich Bugaev} The problems he addressed are the same as those that have preoccupied Florensky.
In the introduction, Luzin\index{Luzin, Nikolai Nikolaevich} writes: 
\begin{quote}\small
In order not to get bogged down in metaphysical discussions on the meaning of the word infinite,\index{infinity (philosophical notion)} I shall confine myself here to discussing a particular point that no doubt attracts the attention of all mathematicians: this is an attempt to satisfactorily resolve the so-called\index{Florensky, Pavel Aleksandrovich} continuum problem,\index{continuum problem} posed by Georg Cantor\index{Cantor, Moritz} some sixty years ago.
\end{quote} 

In the conclusion of the talk, after a tour of several problems concerning the foundations, Luzin\index{Luzin, Nikolai Nikolaevich}  ends with the sentence: ``Over the centuries, the great debate between d'Alembert \index{Alembert@d'Alembert, Jean le Rond} and Euler\index{Euler, Leonhard} on the
notion of arbitrary function now becomes the debate on arbitrary real numbers
and even arbitrary positive integers."\footnote{Luzin refers to a controversy which 
 involved some of the
brightest mathematicians of the eighteenth century, among whom are Euler,
d'Alembert, Lagrange and Daniel Bernoulli, on the question ``what is a function?". I have surveyed this question in \S 10 of the article \cite{Riemann}.}

  In 1930, Luzin\index{Luzin, Nikolai Nikolaevich} published his monograph \emph{Leçons sur les ensembles analytiques et leurs applications} \cite{Luzin-Lecons}, with an appendix by  Sierpi\'nski and a preface by Lebesgue.\index{Lebesgue, Henri} The book was published in Paris. 
 Lebesgue is full of praise for Luzin in the preface. At the end of this preface, he writes: 
  \begin{quote}\small
  [Luzin] has philosophical worries, and he admits this. For him, mathematical and philosophical imperatives are constantly associated, one might even say merged. Although his book is a mathematical work, written by a mathematician for mathematicians, this intimate association of philosophical and mathematical thought is clearly evident on almost every page, giving the volume exceptional importance and appeal.
  \end{quote}
  
At Moscow University, Luzin was a charismatic figure. A mathematical community was formed around him, affectionately called Luzitaniya,\index{Luzitaniya} whose members were linked not only by mathematics, but also by culture. S. S. Demidov and V. D. Esakov write in \cite[p. 3]{Demidov1991-1}  that Luzin  was ``the acknowledged leader of the Moscow School of function theory,
the heart and soul of the famous Luzitaniya, [\ldots] a
true brotherhood united around their adored teacher." They add \cite[p. 10]{Demidov1991-1} that ``[Luzin] had
a charismatic personality that was attractive to young people" and that ``[he] felt most in
his element when working in a tight circle organized like a family gathered around
an adored father-teacher." Alexey Sossinsky tells me that no such mathematical community has ever existed since, except perhaps for the early Bourbaki members.

 Luzin's\index{Luzin, Nikolai Nikolaevich} last years were saddened as he became  persecuted by the state and rejected by most of his former students. His fate was not comparable in severity to that of his old friend\index{Florensky, Pavel Aleksandrovich} Pavel Florensky, but still, he was betrayed, he fell into depression, and never recovered.  I will say a few words on this now.

The violent campaign against Luzin started in the summer of 1936. The purpose was to expel him from the Academy of Sciences. The period was that of a relative political quiteness, just preceding Stalin's Great Purges, which started a  year before Florensky was shot.  A document reproduced in the book   \cite{Demidov1991-1} by Demidov and L\"evshin \cite[p. 12-13]{Demidov1991-1} shows that as early as in 1931, a confidential letter signed by Ernst Yaromirovich Kol'man (1892-1979),\footnote{Kol'man was a virulent Marxist who, in his own way, brought order to Moscow's mathematical community in the 1930s. He is among those who attacked Florensky's \emph{Imaginaries}, ridiculing ``mathematics in the service of religion", ``mathematics in the service of priestcraft" and ``mathematics in the service of obscurantism" it contains \cite[p. 608]{De-F}. He published an article in the Moscow journal \emph{Bolchevik} titled \emph{Against the recent revelations of bourgeois obscurantism}. See pages 114 ff. of \cite{Imaginaires}.} who became president of the Moscow Mathematical Society  after Egorov was arrested in 1930, was sent to the divisions of the Central Committee, in which he  accuses Luzin\index{Luzin, Nikolai Nikolaevich} of several offences, among them, that of being closely connected with \'Emile  Borel, ``an active member of the military [French] establishment."\footnote{Borel had been the Minister of Naval affaires, between April 17 and November 28, 1925, in the government of Paul Painlevé.} In the same letter, Luzin is accused to be an idealist militant, and as a proof of that, Kol'man quotes a passage from a report in which Luzin\index{Luzin, Nikolai Nikolaevich} writes: ``The series of natural numbers does not seem to be an absolutely objective structure. It seems to be an artifact of the brain of the mathematician who happens to be speaking about the natural numbers. There seem to me, among the problems of number theory, some that are absolutely unsolvable." \cite[p. 12]{Demidov1991-1}

  The official public campaign against Luzin\index{Luzin, Nikolai Nikolaevich} started with two anonymous articles in the \emph{Pravda}, with the titles \emph{Answer to Academician N. Luzin} (July 2, 1936) and \emph{On enemies in Soviet masks} (July 3, 1936). The same year, a ``Committee concerning N. Luzin" was created at the Presidium of the Academy of Sciences of the USSR. The committee was headed by the Vice-president of the Academy, G. Krzhizhanovsky, and its members included several well-known mathematicians, among them I. Vinogradov,  L. Shnirelman, S. Sobolev,\index{Sobolev, Serge\"\i\   Lvovich}
P. Aleksandrov and  A. Khinchin. Kolmogorov\index{Kolmogorov, Andrei Nikolaevich} and Lyusternik,\index{Lyusternik, Lazar Aronovich} without being members of the Committee, were invited
to attend the meetings.   Most of the committee members were Luzin's\index{Luzin, Nikolai Nikolaevich} former students or students' students.        
     
    The academic accusation\index{Florensky, Pavel Aleksandrovich} claimed that Luzin,\index{Luzin case}\index{Luzin case} by the end of the 1920s, showed a lack of interest in the  mathematical schools that were emerging in the Soviet Union, and was only interested in his own field. 
   The  prosecution also claimed that he was putting pressure on his students to include his name in their publications, or to make them acknowledge his help in their work. We already mentioned that Luzin continued the work of his student Suslin,\index{Suslin, Mikhail Yakovlevich} after the latter's death. In fact, he was accused of stealing this work. No one will ever know if there is any truth in this, because no one will ever know to what extent Suslin shared unpublished ideas with Luzin. I learned from Sossinsky that some of Luzin's students tried to appropriate Suslin's results after the latter's death.
In any case, the claims that Luzin stole others' ideas contradict what Lebesgue\index{Lebesgue, Henri} wrote in the preface to Luzin's\index{Luzin, Nikolai Nikolaevich} book \cite{Luzin-Lecons}: ``Mr. Luzin is perfectly happy only when he has managed to attribute his own discoveries to someone else," adding that this is a ``strange mania, forgivable, it seems, since it is hardly to be feared that Mr. Luzin will be a model in this respect." Luzin was also blamed, by the Academic court, for the fact that he published his papers in French rather than Russian journals, which was claimed to be a disloyalty to his country.

In the end, and surprisingly, Luzin\index{Luzin, Nikolai Nikolaevich} was not condemned as an  ``enemy of the people", as was planned by his detractors. The common explanation  is that the physicist 
Piotr Leonidovich\index{Kapitsa, Piotr Leonidovich} Kapitsa\footnote{Kapitsa was awarded the Nobel prize in physics in 1978.
In 1946, he refused to work on the Soviet atomic weapons, and this led to his dismissal from the Academy, a position he took up again in 1955, two years after Stalin's death. At his death in 1984, Kapitsa was the only member of the Russian Academy of Sciences who was not a member of the Communist Party.} (1894-1984), who was the head of the Physics section of the  Russian Academy of Sciences, wrote a letter to Stalin, saying that Luzin might yet be
useful to the government, see  \cite{Hoffman} where the author comments: ``It is not clear why
Stalin listened, but his whim ensured the future
of a discipline." 
Still, after this affair, Luzin remained rejected by an important part of the mathematical establishment and the case left a deep wound in the Russian mathematical community.\footnote{In an email I received on April 29, 2019, Alexei Sossinsky writes:  ``When the Committee concluded its work, it was expected that it would declare that Luzin is ``an enemy of the people', which meant that he would be sent to the camps or even be given the death penalty. But it didn't! Until a few years ago, we did not know why this happened. The publication of certain documents now gives us the correct story, which involves Kapitsa's letter, the falling of Kol'man into disfavor, Krzhizhanovsky, Molotov, and Stalin, the final decision apparently being taken in a tête-a-tête between Krzhizhanovsky and Stalin."}  He continued working, but in isolation.\footnote{I learned from Yuri Neretin that in 1937--39, Luzin\index{Luzin, Nikolai Nikolaevich} published an ingenious work on some types of bendings of two-dimensional surfaces, solving a problem that had been discussed  during 50 years, see \cite{Luzin1939}. Unfortunately for him, the solution was negative and  his work closed the topic and had no direct continuation. Sabitov wrote a review of  this problem  in   \cite{Sabitov}, in which he says: ``a truly dramatic end
to this direction of geometry came with Luzin's paper \cite{Luzin1939}."}  In 1946, Luzin\index{Luzin, Nikolai Nikolaevich} published a 900 pages textbook on calculus for technological universities. In 1953 (the year Stalin died),  300,000 copies of the book were already sold.

  This so-called ``Luzin case,"\index{Luzin case}  a combination of a generation dispute and personal interests, is considered to be one of darkest episodes of the history of Soviet mathematics, if not the darkest one, see the details in \cite{Demidov1991-1}.
Discussion of these events remained a taboo in the Soviet Union, until their re-examination several decades later.   In
   2012, the Russian Academy of Sciences reversed the decision that was taken against Luzin\index{Luzin, Nikolai Nikolaevich} in 1936. A reconstruction of the minutes of the various trials that took place is carried out in  book 
 \cite{Demidov1991-1} by Demidov and Lëvshin. The reader may also consult Demidov's article \emph{The Moscow school of the theory of functions  in the 1930s} \cite{Demidov1993} which constitutes an excellent concise report on the school founded by Luzin.\index{Luzin, Nikolai Nikolaevich}   The Lavrentieff\index{Lavrentieff, Mikhail Alekseyevich}  memorial book \cite{Lavrentieff-time} contains a chapter dedicated to the Luzin case.
    
  Today, it is clear that the reason behind Luzin's persecution was his religious convictions. The same thing had happened to Egorov. The latter, who had not been politically active, fell into disfavor with the Soviet authorities six years before Luzin. Demidov and Esakov write in \cite[p. 10]{Demidov1991-1} that ``[Egorov] disliked the new regime
and made no attempt to hide his negative attitude."     The campaign against him started in 1929 when, Egorov, at the first All-Union Congress of Mathematicians in Kharkov, opposed a proposal to send greetings on behalf of the Mathematical Congress to the Sixteenth Party Congress, see \cite[p. 11]{Demidov1991-1}. The\index{Florensky, Pavel Aleksandrovich} same year, after having been President of the Moscow Mathematical Society for seven years, he was dismissed from his position at the University.
  The next year he was jailed for being a ``religious sectarian." About Egorov's religiousness, see \cite{Kolyagin}.
  In fact, Egorov was accused of being part the so-called Catacomb Church  \cite[p. 12]{Demidov1991-1}.\footnote{The ``Catacomb Church" is a generic name for underground Christian communities in the Soviet Union who no longer felt part of the official Orthodox Church after the state took over its activities. The name makes reference to the Christians worshipping in the catacombs of Rome during the period of the Roman persecutions against Christians.}  Little is known about this sad story and the information about Egorov's  the last years contains contradictions. Most probably, he died after a hunger strike. According to Demidov, he died in the hospital of the Institute for the Continuous Education of Doctors
(a branch of Kazan University) \cite{Demidov1999}. I learned from V. N.~Berestovski\u\i \ that Egorov was buried in the same cemetery as
N. I. Lobachevsky, in Kazan.

 On the Luzin case,\index{Luzin case} I have benefited a few years ago from correspondence with Yuri Neretin and Alexey Sossinsky. There are many comments on this in Neretin's article \cite{Neretin}.
  Let me end this section by a last information about Luzin, which I learned from Sossinsky.  The year before his death, Luzin thought he had a counterexample to the Goldbach conjecture. He sent his result to  the  \emph{Doklady Akademii Nauk}. The example was numerical, with huge numbers involved. It took a large amount of time for the persons in charge of checking the result to find an error in the calculations.

 The book \emph{Naming infinity: A true story of religious mysticism and mathematical creativity} by  Graham and Kantor \cite{Graham}, despite many historical, mathematical and theological inaccuracy,\footnote{Regarding the inaccuracy, see the corrective reviews by  Hoffman \cite{Hoffman}, Trimble \cite{Trimble},  Katasonov \cite{Katasonov}, to mention just a few. See also the interesting bachelor thesis by Eriksson  \cite{Eriksson}.     Katasonov \cite{Katasonov} argues against the claim which is repeated all over the book that Name Worshippers consider that naming things create them, whereas these worshippers never associated naming with creating.  Trimble \cite{Trimble} notes that  Graham and Kantor have no knowledge and a false reading of Orthodox theology. For instance, they conflate the hesychastic tradition as a whole with Name Worshippers. Eriksson writes that the book by Graham and Kantor with, ``lamentably, very narrow view of the
cultural context, attempts to analyse the connection between culture and mathematics." She adds that in their treatment, 
 ``Bugaev is left out, and in this way the mathematical roots of the tradition
are cut. Graham and Kantor even seem to be unaware of Bugaev's contributions, as Bugaev is not
mentioned at all in the discussion about continuity and discontinuity in Florensky's work."\index{Florensky, Pavel Aleksandrovich}
  Personally, I do not see any reason for Graham and Kantor's insistance on homosexuality in the Russian mathematical community they describe, and for me their suggestion that Florensky was homosexual is bewildering. I think  that it is based on a misinterpretation of his religious style of writing and a misunderstanding of his personality. In any case, homosexuality was severely condemned in Russia during the Soviet era, and Florensky was being very closely monitored. If there had been the slightest hint of his homosexuality, this would have been known and recorded by the KGB, but there seems to be nothing to suggest this in the documents found and published in \cite{Chentalinsky} and elsewhere.} contributed in the rising of interest among Western readers in the history of Soviet mathematics.

\section{Silver Age: Mathematics, symbolism, art and poetry} \label{s:Silver-age}

\bigskip 
 ``They call it `Silver' Age, but it is comparable to the Renaissance in Europe", Valentin Poénaru told me.

In the chapter \emph{The Practice of Mathematics in a Totalitarian Society}  published in the present Handbook, Alexei Sossinsky talks about his past as a student at Mekhmat, the Faculty of Mechanics and Mathematics  of Moscow State University, and then as a mathematician  at this faculty, during the period called the ``Golden Age of Russian mathematics".\index{Golden Age of Russian Mathematics} This is the epoch of an unprecedented development of mathematics in the Soviet Union, and it involves fields such as geometry, topology, number theory, function theory, probability theory, mathematical logic, and others.   With this period are associated the names  of Pavel Sergeyevich Aleksandrov\index{Aleksandrov, Pavel Sergeyevich} (1896-1982),  Mikhail Alekseyevich Lavrentieff\index{Lavrentieff, Mikhail Alekseyevich}  (1900-1980), Andreï Nikolaevich Kolmogorov\index{Kolmogorov, Andreï Nikolaevich} (1903-1987),  Israel Moiseyevich Gelfand\index{ Gelfand, Israel Moiseyevich} (1913-2009), Yuri Ivanovich Manin\index{Manin, Yuri Ivanovich} (1937-2023)  and several others, forming at least two generations of first-rank Russian mathematicians. The impact of their work on the development of mathematics worldwide was enormous.  The reader interested by the history of this period of Russian mathematics may refer to the book \cite{Golden} edited by  Zdravkovska  and Duren. 

The Golden Age\index{Russian Golden Age!mathematics}  of Russian mathematics was preceded by another, sometimes referred to as the ``Silver Age of Russian mathematics",\index{Russian Silver Age!mathematics} which has been the subject of much attention in recent years among historians of mathematics, see e.g. the articles \cite{Demidov2021, Demidov1993, Demidov1999}, by Demidov and \cite{Ford} by Ford. This is the period dominated by the figure of Luzin,\index{Luzin, Nikolai Nikolaevich} which we presented in  \S \ref{s:Luzin}. Luzin mentored and had a significant impact on several of the key players of the Golden Age.\index{Russian Golden Age!mathematics} This Silver Age\index{Russian Silver Age!mathematics}  of mathematics overlaps, by a few years, with a period also known as the ``Silver Age", for  poetry, literature and art. It is a period which saw the birth of the Russian Symbolist movement, one of whose main protagonists was Andrei\index{Belyi, Andrei} Belyi (Bugaev's son), whom we mentioned in \S \ref{s:Bugaev}. In 1910, Belyi published an important collection of articles bearing the characteristic name \emph{Symbolism},\index{symbolism} to which Florensky\index{Florensky, Pavel Aleksandrovich} refers in several writings and in particular in his major theological\index{Florensky, Pavel Aleksandrovich!\emph{The pillar and ground of the truth}} work, \emph{The Pillar and ground of the truth}  \cite{Flo1}. 

Let me add that Andrey Andreyevich Markov (1856-1922) the outstanding Russian mathematician who worked in Saint Petersburg, best known for his work on stochastic processes, and Andrei\index{Kolmogorov, Andreï Nikolaevich} Kolmogorov, the student of Luzin\index{Luzin, Nikolai Nikolaevich} whom we already mentioned and who became the world leader of probability theory, algorithmics, complexity and information theory, were both interested in the statistical approach to verse analysis in Russian poetry, a 
topic which was started by Andrei Belyi.\index{Belyi, Andrei} Both Markov and Kolmogorov applied this approach to Pushkin's poetry, see \cite{Markov}.
More recently, Sossinsky wrote a paper on this topic, titled \emph{Belyi--Kolmogorov theory} \cite{Sossinsky2}.

  Florensky, who belonged to the group of mathematicians of the Silver Age,\index{Russian Silver Age!symbolism} was connected at the same time with the Russian symbolist\index{symbolism} movement, with its renewed interest for classical literature and for the significance of language. Since the times he was a student at Moscow University, he had close relations with Moscow's Symbolist literary circle.\index{Moscow Symbolist literary circle}   He published several articles in symbolist periodicals such as \emph{Novyi Put'} (New Way) and \emph{Vesy} (Scales).       
  In 1903, he published in \emph{Novyi Put'} an article titled \emph{Superstition} \cite{Florensky-Superstition},  which is concerned with supernatural phenomena and the psychology of their perception. In 1904, he published in the same journal an article titled \emph{On symbols of infinity: Sketch of the ideas of G. Cantor} \cite{Florensky-Symbols},  to which he referred later in the \emph{Pillar}.\index{Florensky, Pavel Aleksandrovich!\emph{The pillar and ground of the truth}} This article is the first essay ever published in Russian on Cantor's ideas on infinity.\index{infinity (philosophical notion)} In the same year, Florensky published in \emph{Vesy} a series of poems, and an essay titled \emph{On a certain premise of a worldview} \cite{Florensky-certain}.  In 1905, he wrote the essay we mentioned, with symbolist inclinations, \emph{Hamlet}, on Shakespeare's play which many have seen as a work in which the author anticipates in several ways the modern world. This essay was meant to be published in the symbolist journal \emph{Vesy} but was never published before it was included in Volume I of Florensky's Collected Works published in Moscow in 1994. It was recently translated into French, in the form of a booklet, see \cite{Hamlet}. 
   
   Talking about symbolism\index{symbolism} in the Russian Silver Age,\index{Russian Silver Age!symbolism} one should mention the names of some leading symbolist poets of that period:
 Vyacheslav Ivanovich Ivanov\index{Ivanov, Vyacheslav Ivanovich} (1866-1949), a very complicated figure with a long and serpentine biography,
Valery Bryusov\index{Bryusov, Valery } (1873-1924) and
Aleksander Blok\index{Blok, Aleksander} (1880-1921). 
I like the following description of symbolism\index{symbolism} by Ivanov who,  according to Sabaneeff \cite[p. 314]{Sabaneeff}, knew Florensky\index{Florensky, Pavel Aleksandrovich} and had a high esteem for him: Symbolism\index{symbolism} is ``a cult of eternity and all-embracing unity of the multicolored reflections of the moments" \cite{Strakhovsky}.\footnote{Strakhovsky in \cite{Strakhovsky} is translating Ivanov, from his article \emph{Of the Merry Craftsmanship and Smart Revelry}, published in the ephemeral periodical \emph{Zolotoe Runo} (1906-1909) \cite[p. 54]{Ivanov}).} Finally, in referring to the Russian Silver Age\index{Russian Silver Age!art} in art, one may  mention the names 
of  Chagall, Kandinsky,\index{Kandinsky, Vassily} Stravinsky, Prokofiev and Skriabin, to name but a few of the best known figures.

Florensky\index{Florensky, Pavel Aleksandrovich} described symbolism\index{symbolism} as ``the necessity of seeing the concrete incarnation of the soul". In his autobiography dedicated to his children, he writes: ``All my life I have pondered a single problem, that of the symbol.  [...] Positivism repelled me, but abstract metaphysics repelled me just as much. I wanted to see the soul, but I wanted to see it incarnate. Someone will say that this is materialism, but it is not materialism that is at stake, but the necessity of the concrete, of symbolism."\index{symbolism}
In the same autobiography, we can read: ``I have always been a symbolist. Symbols are created by the mind according to precise laws, and with an inner necessity, and this happens whenever certain aspects of the mind begin to function in a particularly lively way. What symbolizes and what is symbolized are not fortuitously interrelated."  With all this, it is not surprising that Florensky's major theological and philosophical work, \emph{The Pillar and ground of the truth}\index{Florensky, Pavel Aleksandrovich!\emph{The pillar and ground of the truth}}, is written in the language of symbols. He writes there, in the X$^{\mathrm{th}}$ letter \cite[p. 236]{Flo1}: 
``I am
compelled to use a metaphysical terminology, but in my speech these
terms have not a strictly technical sense but a conventional or rather a
symbolic one. They have the significance of colors by means of which
inward experience is painted."
The language he uses is also the one of the \emph{Imaginaries}.

Speaking of art and poetry, I would like to mention Florensky's\index{Florensky, Pavel Aleksandrovich} inclination towards poetry and music.
Sabaneeff writes, in an article I have already mentioned several times, that Florensky ``was too far removed from music and from art in general." 
I think Sabaneeff was wrong.\footnote{Sabaneeff bases his opinion on the fact that when the poet Vyacheslav Ivanov introduced Florensky\index{Florensky, Pavel Aleksandrovich} to the world of Aleksander Skriabin, Florensky did not fit in there, see \cite[p. 315]{Sabaneeff}.} In his autobiography and correspondence, Florensky regularly talks about music and poetry. As a child, he received a musical education, and everybody in his family listened to and played music. Let me quote a small passage from his recollections \cite[p. 80]{Florensky-bio}: ``As a child, my musical inclinations also led me to poetry. I was less interested in the meaning of verses than in their sound and rhythm. Endowed with an almost absolute memory, I was able to memorize everything that interested me, especially poems, instantly and with total precision." Florensky's writings are so full of musical references, thoughts on music and poetry, that one could write a book on the subject: \emph{Music\index{Florensky, Pavel Aleksandrovich} and poetry in the works of Pavel Florensky}.

Florensky also wrote poetry. Many of his poems were confiscated and are still lost, probably forever, but the volume of correspondence\index{Florensky, Pavel Aleksandrovich} \cite{Florensky-Letter} contains an unfinished 60-page poem, ORO, which he wrote during his exile. It is translated into French by Fran\c coise Lhoest.

Regarding sound and music, let me quote a passage from Florensky's\index{Florensky, Pavel Aleksandrovich} correspondence (letter addressed to his mother dated 28-29.VII.1936) \cite[Letter no. 70, p. 447]{Florensky-Letter}):
  
\begin{quote}\small
An acquaintance of mine asks me why I never write about the sounds of Solovki, but only about colors and shapes. This is because everything here is soundless, like a dream. This is the realm of silence. Not in the literal sense, of course; there is more than enough unwanted noise, and one would want to get away from it. But you cannot hear the inner sound of nature, you cannot perceive the inner words of beings. Everything glides along, in a shadow theater, and sounds are added from outside, like an appendix or a disturbing noise. It is hard to explain why nothing produces sound, why there is no music in things and in life, and I do not really understand, but all the same, there is no music.

\end{quote}

Finally, talking about sound and music, I would like to quote another passage, this time from Florensky's autobiography, in which he describes how he experiences, in both a physiological and symbolic way, Chladni's\index{Chladni figures} figures, the  acoustic lines\index{acoustic lines} formed on a plate rubbed by a bow.\footnote{Ernst Florens Friedrich Chladni\index{Chladni, Ernst Florens Friedrich} (1756-1827) was a German physicist, sometimes considered as the founder of modern acoustics.  He conducted experiments consisting in making a plate vibrate by scraping it on the side with a bow, after covering it with a layer of fine sand, and observing the patterns the sand produces: geometric figures, always regular, but very different, depending on where and how the plate is excited. These are the \emph{Chladni figures} mentioned by Florensky  in the following passage.}   He\index{Florensky, Pavel Aleksandrovich} writes \cite[p. 66]{Florensky-bio}:

\begin{quote}\small
I was feeling almost physically like a string, or rather like a Chladni plate, over which nature passes like a bow: in my soul, or rather in my whole organism vibrated, almost audible to the ear, a pure, high, elastic sound, and in my thoughts were ordered schematic figures, like Chladni\index{Chladni figures} figures, like symbols of worldly phenomena. I write and I am almost certain to be misunderstood. People will want to hear in these words comparisons and poetic essays, whereas I want to express, to draw from myself, the soberest and most literal description possible, a kind of physiological picture. It consisted of this: everything in me, every vein, was full of an ecstatic sound that concentrated within itself \emph{my} knowledge of the world, generating patterns of a rather mathematical order: my categories of knowledge. 
\end{quote}

 \begin{figure}
\begin{center} 
\includegraphics[width=10 cm]{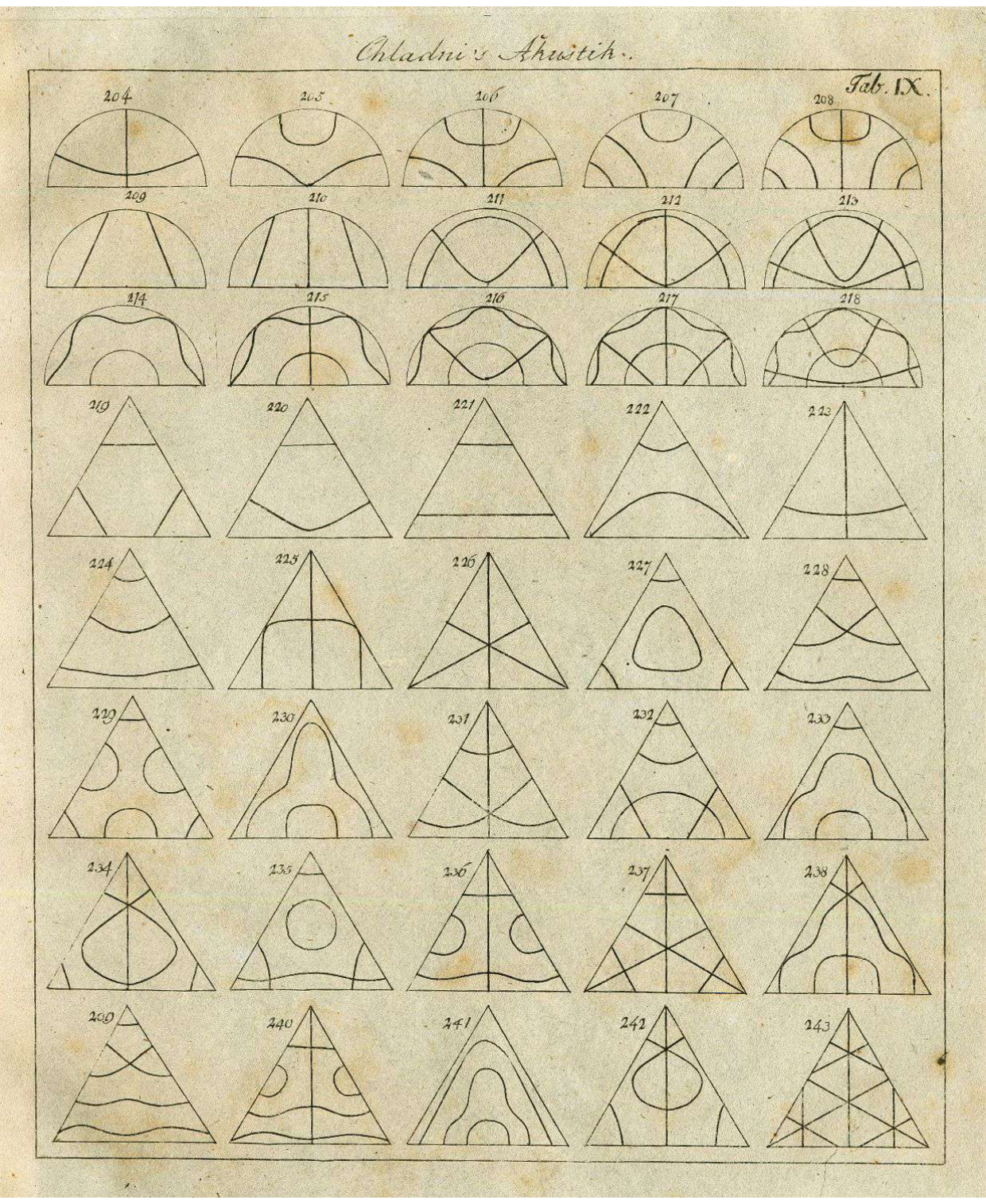} 
\caption{Chladni acoustic figures, from Chladni's \emph{Akustik} \cite{Chladni}}
\label{fig:Chladni}
\end{center}
\end{figure}

Talking about poetry, let me mention Pasternak\index{Pasternak, Boris} again.
Chentalinsky, in \cite[p. 164]{Chentalinsky}, quotes him describing Florensky as a ``hostage of eternity; captive of time." Arkady Plotnitsky tells me that the line comes from a well-known and much quoted poem from the 1950s by Pasternak, and that this poem is addressed to an artist in general, and not particularly to Florensky. But the poem fits  Florensky very well. Let me quote here the corresponding stanza and the one before:

\begin{verse}
  Don't sleep, don't sleep, work
  \\
Don't interrupt your labor,
  \\
Don't sleep, struggle with a slumber,
  \\
Like a plain pilot, like a star.

Don't sleep, don't sleep, an artist 
  \\
Do not submit to sleep,
  \\
You are a hostage of eternity,
  \\
In time's captivity.
 
\end{verse}

 In this section, I have talked about poetry in the Silver Age. In the next section, I will continue with, arguably, the greatest poet of the Silver Age.\index{Russian Silver Age!poetry}

\section{Florensky and Akhmatova}\label{s:Akhmatova}

 If the expression ``Silver Age\index{Russian Silver Age!mathematics} of Russian mathematics" is used in reference to the Golden Age\index{Russian Golden Age!mathematics} that followed it, the expression ``Silver Age of Russian poetry" is used in comparison with the Golden Age\index{Russian Golden Age!poetry} that \emph{preceded} it, that is, the age of Pushkin and Lermontov.\index{Golden Age of Russian poetry} Some consider that the expression ``Silver Age\index{Russian Silver Age} of poetry" was coined after certain verses in Akhmatova's\index{Akhmatova, Anna} poem \emph{Poem without a hero},\index{Akhmatova, Anna!\emph{Poem without a hero}} namely,
  \begin{verse}
  And the silver month is bright
\\
Frozen over the Silver Age.
 \end{verse}
 
 Akhmatova,\index{Akhmatova, Anna} who was dubbed ``the Soul of the Silver\index{Russian Silver Age!poetry} Age,"
  was close to Florensky\index{Florensky, Pavel Aleksandrovich} in many respects, and I will talk about this in the present section. 
  Let me start by quoting a poem she wrote on Dante\index{Dante}  \cite{Akhmatova}:\footnote{The following translation is by Alexei Sossinsky.}

\begin{verse}
From death there will be no return
\\
To the old Florence he adores.
\\
And leaving, back forever turned,
\\
For him I sing this song, this lore.

\bigskip

Dark night, fire, torches, last embrace,
\\
Beyond the Gate, the howl of Fate.
\\
From Hell he sent his final curse,
\\
In Heaven Florence he will not forget.

\bigskip

And barefoot, wearing shirt of shame,
\\
With flaming candles lightly carried,
\\
He walks in Florence all the same,
\\
So trech'rous, base, but so desired.

 \end{verse}

In the middle section of her \emph{Poem without a hero}, Akhmatova\index{Akhmatova, Anna!\emph{Poem without a hero}} has the following epigraph:

\begin{verse}

\ldots a jasmine bush that bloomed  \\
Where once Dante walked and the air is mute.

\end{verse}

These two lines by Akhmatova are a variation on (and a truncation of) the following lines lines in a poem about Akhmatova's silence during several years, by Nikolai Klyuev:\footnote{The symbolist poet Nikolai Alekseevich Klyuev\index{Klyuev, Nikolai Alekseevich} (1884-1937) was arrested and sent to the camps in 1933 (the same year as Florensky).  He was shot to death in 1937 in Tomsk (two months before Florensky) and rehabilitated posthumously in 1957 (one year before Florensky). The translation here is from the article \cite{Mager}. Although the details would take us far beyond the scope of this chapter, let me mention that I learned from Plotnitsky that Akhmatova's truncation  and the change of meaning of Klyuev's passage was likely to be mediated by Mandelstam (with Klyuev's poem still as the original source), and that there is a controversy about this because Akhmatova changed her poem several times, see the article \cite{Leving} by Y. Leving.}

\begin{verse}

Akhmatova---a jasmine bush, \\
Scorched by asphalt and gray pitch---\\
Is she lost on the path through the dark pits\\
Where Dante passed and the air was parched,\\
Or is she a nymph spinning flax into crystals?\\
Among Russian women, remote and subtle\\
Anna is a cloud at sunset, she reflects\\
Streaks of gray light like rockets!
\end{verse}

Dante was for Akhmatova,\index{Akhmatova, Anna} like for Florensky, a source of inspiration and reflection.  They both regarded the fourteenth century Florence \emph{Podestà}, which condemned Dante to an exile with no return, as a prelude for the twentieth century brutal Soviet regime. The long journey that Dante recounts in the \emph{Divine Comedy} symbolized for Akhmatova\index{Akhmatova, Anna} all the exiles she had known, people she loved who  emigrated to Europe, and  those who were sent to the Gulag. Among the latter, there were poets, writers, and her son.  
 
 Akhmatova\index{Akhmatova, Anna} had Eurasian origins: her father, like Florensky's, was an engineer, of Greek descent, and her mother was of Tatar descent.
Those who knew her personally described her as a strange woman: shy, gentle, sad, indifferent to worldly things and with an acute taste for contemplation.  All this coincides with what we know of Pavel Florensky.\index{Florensky, Pavel Aleksandrovich} Contradictions were noted in her social behavior.  Her first husband, the prominent poet Nikolai Gumilev,\index{Gumilev, Nikolai Stepanovich} called her an enchantress. It was said that she was both a saint and a harlot, and that she had diabolical charisma.\footnote{On August 14, 1946 the decree signed by  Andrei Zhdanov, who was Stalin's cultural minister, which banned Akhmatova\index{Akhmatova, Anna} from publishing, accused her (in the present case, with an outright malice)
of mixing debauchery with prayer,  describing  her as a "bourgeois lady gone mad, tearing between her boudoir and her altar  \ldots half nun, half harlot, or rather both nun and harlot".} Apart from the harlot part of it, Florensky's social image was similar. 
I already quoted Sabaneeff, who wrote that there was in him  ``something of black magic, dark, devoid of divine grace."  In the same article, we read \cite{Sabaneeff}: ``Beyond any doubt, there was a demoniac
element in him, and as indubitably he was an extraordinary
man, an outstanding personality, quite beyond comparison
with any other prominent man of his time. \ldots Lucifer was
closer to him than Christ." Sabaneeff adds: ``This increased rather than diminished
his fascination for me."

Detached from everything, Akhmatova\index{Akhmatova, Anna} was trying to understand life through poetry. Like Florensky,\index{Florensky, Pavel Aleksandrovich} she refused emigration, even though she had friends in Italy and France. For this, let me quote the opening of her \emph{Requiem}:\footnote{Transl. Robin Kemball.}

\begin{verse}
No, not far beneath some foreign sky then, \\

                        Nor with foreign wings to shelter me, --- \\
                        I was with my people then, close by them,\\
                        Where my luckless people chanced to be.\\
                        \hfill 1961
\end{verse}

Akhmatova was neither exiled nor shot like Florensky, but she experienced exile and execution within herself:  Her first husband, Nikolai Gumilev,\index{Gumilev, Nikolai Stepanovich} accused of conspiration, was shot in 1921. Her son, Lev, endured the Siberian camps. Her lover, the art historian and critic Nikolai Punin,\index{Punin, Nikolay Nikolaevich} died in the Stalin camps, and one could continue this gloomy list.\footnote{An exception to this ``gloomy list" is her last lover, Isaiah Berlin,\index{Berlin, Isaiah}  an English diplomat and academician, born in Riga, which belonged to the Russian Empire, who was twenty years her junior and who came to see her in Leningrad while he was completing an internship at the British Embassy in Moscow. Reportedly, when Stalin learned of this affair, he said: ``Now the nun is receiving foreign spies at night!"
Akhmatova was deeply moved by Berlin, whom she dubbed the ``Guest from the future", and to whom she dedicated four cycles of poems.  
In one of her conversations with him, she declared, ``I have only darkness and Rus' by my side"  \cite{Victoroff}. Berlin  was also deeply moved by Anna. On his return to England, he arranged for her to be made Doctor Honoris Causa of Oxford University. } We already talked about Florensky's very insecure position regarding publication: he was deprived of the opportunity to publish in scientific journals, and his notes and manuscripts were confiscated.
Akhmatova's\index{Akhmatova, Anna} writings were judged by the authorities, but also by the Writers' Union, as ``socially too irrelevant", and her poetry was banned from publication, starting in 1922, for over thirty years. 
She was led to recite her poems secretly to close friends, asking them to memorize them. Thus, it was her works that were deported instead of her. 


 Akhmatova\index{Akhmatova, Anna} had a religious soul. Korney Ivanovich Chukovsky\index{Chukovsky, Korney Ivanovich} dubbed her ``the unique and last Orthodox poet" \cite{Chukovsky}.\footnote{Korney Ivanovich Chukovsky (1882-1969) was a journalist, poet, translator,  literary critic and professor of Russian literature. With Akhmatova, he belonged to the same artistic and intellectual circles, before the revolution.  Chukovsky was the illegitimate son of Emmanuel Solomonovich Levinson, whose (legitimate) daughter, Henrietta Emmanuilovna Levenson, 
was the mother of the outstanding Soviet  mathematician Vladimir Abramovich Rochlin (1919-1984).  Rochlin was the doctoral advisor  of M. L. Gromov, A. M. Vershik, V. G. Turaev, O. Y. Viro,   Y. M. Eliashberg and several other bright mathematicians. 
Chukovsky's daughter, Lydia Korneyevna Chukovskaya  (1907-1996)  was a Russian and Soviet author and literary critic. Her second husband,
 Matvei Petrovich Bronshtein (1906-1938), was arrested as part of Stalin's Great Purge and  was shot on February 18, 1938.
Bronshtein was an outstanding Soviet theoretical physicist, working  on quantum theory, astrophysics, semiconductor theory,
cosmology, quantum gravity theory. His doctoral dissertation, \emph{Quantization of gravitational waves}, defended in 1935, was the first serious work on quantum gravity. He indicated there the main difficulties in this theory. In 1965, Lydia  Chukovskaya   published a novel, \emph{Sophia Petrovna}, in which
 the son of the main character, Sophia Petrovna, the main character,  was imprisoned in Leningrad's jail and later shot. In fact, Sophia Petrovna repeated in many respect the fate of  Lydia Chukovskaya. (I learned most of the content of this note  from V. N. Berestovski\u\i  .)} Like Florensky, she was steeped in Orthodox culture, even though this culture was forbidden. 
 Her poem \emph{Requiem},\index{Akhmatova, Anna!\emph{Requiem}} dedicated to ``seventeen months waiting in prison queues in
Leningrad,"\footnote{The dedication is to the people waiting for hours, often with no result, in this mile long line to get any information about their relatives in prison or in a Gulag. Akhmatova was among these people waiting.} includes whole lines that are variations on Orthodox liturgical texts. 
 
 \begin{verse}
 
 Do not weep for me, Mother, when I am in my grave.  
\\
A choir of angels glorified the hour,
\\
the vault of heaven was dissolved in fire.
\\
Father, why hast Thou forsaken me?
\\
Mother, I beg you, do not weep for me.

\end{verse}
Each line here is taken from a liturgical text, Florensky has a poem on the same theme, titled \emph{Lament of the Mother of God}, which he wrote in 1907, see \cite[p. 89]{Oppo}.
Akhmatova's poem continues:

 \begin{verse}

Mary Magdalene beat her breasts and sobbed,
\\
His dear disciple, stone-faced, stared.
\\
His mother stood apart. No other looked
\\
into her secret eyes. No one dared.
\end{verse}
This is a faithful description of the Orthodox icon of the crucifixion.

%
%
%
%
%
%
%
%
%
%
%
%
%
%

%
%
%

%
%
%
%
%
%
%
%
%
%
%
%
%
%
%
%
%
%
%

  Akhmatova\index{Akhmatova, Anna} was a twin sister of Florensky\index{Florensky, Pavel Aleksandrovich}. The comparison and the close relationship between them could be the subject of another book. 
  
There would still be a third book to write, about the relationship between Florensky and another twin sister : Simone Weil.\index{Weil, Simone} In this respect, I would like to quote a passage from her book\index{Weil, Simone!\emph{La pesanteur et la Grâce}} \emph{La pesanteur et la Grâce} (1947), the chapter titled  \emph{Mystique du travail} (Mysticism of work) \cite{Weil}:  ``The spirituality of work. Work makes us experience in a harassing way the phenomenon of finality bouncing back like a ball; working to eat, eating to work \ldots If we look at one of the two as an end, or at both taken separately, we are lost. The cycle contains the truth."\footnote{Spiritualité du travail. Le travail fait éprouver d'une manière harassante le phénomène de la finalité renvoyée comme une balle ; travailler pour manger, manger pour travailler… Si l'on regarde l'un des deux comme une fin, ou l'un et l'autre pris séparément, on est perdu. Le cycle contient la vérité.} This passage also contains another incarnation of the idea of eternal return, dear to Florensky, which I discussed in \S \ref{s:time} above. Among the strong ties between Simone Weil and Pavel Florensky are a love of mathematics, a strong sense of self-observation, a passion for history and especially for Greek antiquity, extreme mysticism and an intense religious feeling.\footnote{I would like to thank S. Tanabé who suggested to me several relations between Pavel Florensky and Simone Weil.} 

\section{Contradiction and antinomy in mathematics, philosophy and theology} \label{s:Contradiction}

I have tried to allude to Florensky's philosophy, and I shall try to say more about it  in this section.

The Russian Silver Age,\index{Russian Silver Age!philosophy and theology} besides referring to mathematics, poetry and art, was a period of fusion between  philosophy  and theology.
 The representatives of this alliance include
Sergei Nikolaevich Bulgakov (1871-1944),\index{Bulgakov, Sergei Nikolaevich} Lev Isaakovich Shestov\index{Shestov, Lev Isaakovich} (1866-1938), Aleksei Fedorovich Losev\index{ Losev, Aleksei Fedorovich}
(1893–1988), Nikolai\index{ Berdyaev, Nikolai Aleksandrovich} Aleksandrovich Berdyaev\footnote{It has been suggested that the term ``Silver Age", in relation with Russian philosophy, derives from Berdyaev's work \emph{The Brightest Lights of the Silver\index{Russian Silver Age!philosophy and theology} Age} \cite{Berdyaev}.}
 (1874-1948) and
Dmitry Sergeyevich Merezhkovsky\index{Merezhkovsky, Dmitry Sergeyevich} (1866 -1941). The latter was one of the founders of the symbolist movement, to whom Florensky belonged. The others are also connected to him in various ways.

 Florensky's\index{Florensky, Pavel Aleksandrovich} philosophical worldview is based on antinomy\index{antinomy (philosophical notion)}  and I will mostly talk about this now. 

Florensky considered that the world of ideas is based on antinomy\index{antinomy (philosophical notion)}  rather than on logic as one intends it usually; in fact, truth, for him, is grounded on contradiction,\index{contradiction (philosophical notion)} and this holds in particular for the truths of theology. He considered that contradiction and paradoxes are unavoidable in theology. He believed that religious truth is beyond any formal logic, and more generally, beyond any formal categorisation. 
He considered that the Holy Scriptures, including the Old and New Testament, contain contradictory statements that are nevertheless all true, and that these contradictions do not need to be resolved but on the contrary they are an integral part of Christian teaching.  In Letter VI of the \emph{Pillar}, he reviews some sources of antinomy, in Greek Antiquity, passing by Heraclitus, the Eleatics and Plato.
 To understand Florensky's background, it is good to recall first a few facts about contradiction in mathematics, the first field that he studied in depth.

At the beginning of the twentieth century, at the time when Florensky\index{Florensky, Pavel Aleksandrovich} was suffering his existential crisis, mathematics was traversing a \emph{foundational crisis}. 
A strong indication of this crisis is Hilbert's ``second problem", formulated in 1902 at the second congress of mathematicians in Paris, asking whether one can prove the coherence of arithmetic, i.e., whether one can demonstrate that the axioms of arithmetic are not contradictory. In reality, it was not only a question about arithmetic, but of geometry as well, and, by extension, of all mathematics,  because it was shown before (notably, by Hilbert) that the coherence of the axioms of geometry depended on the coherence of the axioms of arithmetic. For the answer to Hilbert's Problem II,\index{Hilbert Problem II} one had to wait for the work of Kurt Gödel,\index{Gödel, Kurt} who demonstrated that the answer is negative. Gödel published his famous \emph{Incompleteness Theorem} in 1931 \cite{Godel1931},\index{Gödel, Kurt!Incompleteness Theorem} while Florensky was already in the glaciers, far away from it all. In any case, mathematicians who had thought about these questions had always felt that there might be problems of coherence in formal logic, and this had been so since Greek antiquity. (Even though formal logic did not exist yet, it was formally the same logic.) The paradoxes of Zeno of Elea are generally quoted in this respect. Let me  recall here that Zeno was considered by Aristotle as the father of dialectic, that is, the art of asking questions and defending successively opposite theses with the aim of finding the truth. Florensky\index{Florensky, Pavel Aleksandrovich} mentions Zeno\index{Zeno of Elea} in the\index{Florensky, Pavel Aleksandrovich!\emph{The pillar and ground of the truth}} \emph{Pillar} (see p. 345 and p. 574), in relation with the continuous and the discontinuous, with references to his paradoxes of motion. He quotes classical paper on Zeno's\index{Zeno} aporias, by the French mathematician and historian of science Paul\index{Tannery, Paul} Tannery \cite{Tannery}.  As a matter of fact, Zeno's arguments are not contradictory, they are only paradoxical, and based on the paradoxes that are possible in mathematics. But these paradoxes are presented in Florensky's writings in a philosophical language that makes them appear practically contradictory. Regarding antinomy\index{antinomy (philosophical notion)}  in theology, Florensky refers to the philosophical system of the monk Serapion Mashkin\index{Serapion Mashkin (Father Serapion)}, from the Optina monastery.\footnote{Father Serapion was a former naval officer. The name \emph{Optina} refers at the same time to a small village and a monastery, in the Kaluga region, 250 km to the South-West of Moscow. The monastery is linked to episodes in the lives of several Russian writers and philosophers. In the summer of 1878, Vladimir Solovyov\index{Solovyov, Vladimir Sergeyevich} brought Fedor Dostoyevsky\index{Dostoyevsky, Fedor} to the monastery following the tragedy of the latter's son's death. Dostoyevsky stayed at the monastery for three days. Several details of his novel \emph{The Brothers Karamazov} appeared following this visit to Optina. The list of the other visitors of the Optina monastery includes Vasily Zhukovsky, Nikolai Gogol, Ivan Turgenev,  Vasily Rozanov and Leo Tolstoy. After the October Revolution, Soviets closed
the Optina monastery and the monks of the were asked to choose between deportation in the Solovki Archipelago\index{Solovki Islands} or working in a secular agricultural cooperative. The monastery itself was declared a Gulag, and various kinds of prisoners were sent there. After these events, Florensky launched a campaign to ``Save Optina", see \cite[p. 607]{De-F}.}  In Note 837 (p. 574) of the \index{Florensky, Pavel Aleksandrovich!\emph{The pillar and ground of the truth}} \emph{Pillar}, talking about Zeno's antinomies, Florensky\index{Florensky, Pavel Aleksandrovich} writes: ``It is precisely on these antinomies,
and their overcoming that the entire philosophical system of Father Serapion is built." Antinomies, but also discontinuity: in the same page (note 836), Florensky writes that ``one cannot
fail to mention Father Serapion Mashkin's attempt to conceive space and time as
composed of finite, further-indivisible elements."\footnote{For a modern treatment of Zeno's aporias and of contradiction in mathematics from the historical point of view, let me mention the paper \emph{Zeno's arguments and paradoxes are 
not against motion and multiplicity but for the separation of true Beings from sensibles} \cite{Negrepontis} by S. Negrepontis.}


Mathematics since at least Euclid is based on definitions and postulates/axioms. The axioms are assumed to be internally consistent, not containing a contradiction. In case an axiomatic system contains a contradiction, the simple rules of logic (\emph{Modus Ponens}) imply that every statement formed in the language of the system is a formal consequence, thus rendering the system most powerful but useless and mathematically unacceptable. Consequently, if a contradiction is detected within a system, then the system is rejected. Frege's\index{ Frege, Gottlob}
 foundation of set theory\index{set theory (historical development)}  was famously found to be contradictory because it contained the rule of comprehension, from which Bertrand Russell\index{Russell, Bertrand} constructed his paradox. However, it is a consequence of G\"odel's\index{Gödel, Kurt!Incompleteness Theorem} Incompleteness Theorem \cite{Godel1931} that internal consistency cannot be proved within the system if the system is reasonably strong. Let me refer here to a paper by Farmaki and Negrepontis, titled \emph{The paradoxical nature of Mathematics} \cite{FN}, in which these authors stand that in a sense all of Mathematics can be considered as a (hopefully) consistent approximation of the truly contradictory. It is not clear whether Florensky\index{Florensky, Pavel Aleksandrovich} has in mind, even in some very intuitive sense, such a view of Mathematics; it would indeed be fascinating if he did.

In any case, the fact that \emph{truth} is antinomic is one of the main themes of \emph{The pillar and ground of the truth}.\index{Florensky, Pavel Aleksandrovich!\emph{The pillar and ground of the truth}} Florensky explained quite simply in his work that reason alone is sufficient to accept or not accept what is not antinomic, and that if faith is needed, it is precisely because truth is antinomic.
Letter VI of the \emph{Pillar} is titled\index{Florensky, Pavel Aleksandrovich!\emph{The pillar and ground of the truth}} \emph{Contradiction}.\index{contradiction (philosophical notion)} In this letter, Florensky analyses different occurrences of antinomy, in philosophy, theology, logic, mathematics, ethics, linguistics and other fields. I short, he develops the thesis that contradiction and antinomy are the basis of knowledge.  As a matter of fact, each of\index{Florensky, Pavel Aleksandrovich} Florensky, Bulgakov\index{Bulgakov, Sergei Nikolaevich} and Losev developed elaborate philosophical and theological systems based on antinomy.\index{antinomy (philosophical notion)} 
 
 At the time he was writing his book, Florensky was witnessing a renaissance in Russian theology, which found its sources in the Jewish Kabbalah \index{kabbalah} and in the writings of the Desert Fathers, a theology which came to be known as \emph{apophatic} theology,\index{apophatic theology}\index{apophatic} or \emph{negative} theology.\index{negative theology} According to this theology, every time one tries to say something positive about God and the transcendent is inevitably led to contradictions or absurdities.  Here, let me quote again Sabaneeff on Florensky: 
 \begin{quote}\small
 His\index{Florensky, Pavel Aleksandrovich} mind was complex, many-storied, and to some extent even
hostile to simplicity. One might even call it a cabalistic mind \ldots His extravagant and excessively luxuriant
thoughts often contradicted one another, but this did not embarrass
him in the least. \ldots He even asserted that every perceived law
 ``generates" its own negation, inasmuch as every law discerned
by the logical apparatus of man is but a part of a real, comprehensive
synthetic law embracing ``all that exists" and ``all that
is possible and thinkable".
 \end{quote}
 
 let me talk now about the notion of infinity.
  
  Florensky was fascinated by the idea of actual infinity which was highlighted by Cantor,\index{Cantor, Georg} and he saw in it a symbolic vision of God. He proposed a philosophical and symbolic interpretation of this concept based on Bugaev's\index{Nikolai Vasilievich Bugaev} mathematical theory of discontinuity.\index{discontinuity (philosophical notion)} 
In 1904, he published the pamphlet titled \emph{Symbols of infinity: Sketch of the ideas of G. Cantor},\index{infinity (philosophical notion)} which we already mentioned,   dedicated\index{set theory (historical development)}  to the analysis of Cantor's\index{Cantor, Georg} ideas on set theory and infinity. 
    According to I. H. Anellis, Florensky's\index{Florensky, Pavel Aleksandrovich} essay was the first article in Russia written on Cantor's set theory \cite{Anellis}.    
In some sense, Florensky was a disciple of 
Cantor\index{Cantor, Georg}, who had already declared in his 1883 work, \emph{Grudlagen einer allgemeinen Mannigfaltigkeitslehre} that everything that has to do with the infinite has theological, mathematical and philosophical connotations that cannot be separated.
 N. Misler writes in \cite[p. 35]{Misler}  
  that Florensky's proposed  philosophy of symbol was in accordance with the wish of the editors of \emph{Novy put' } to promote the relation between the intelligentsia and the Church. The logic of symbols is used by Florensky in his theological arguments. In his writings, the issues about contradictions in formal logic become issues about  contradictions in theology. He declares in Letter VI of the\index{Florensky, Pavel Aleksandrovich!\emph{The pillar and ground of the truth}} \emph{Pillar} \cite[p. 106]{Flo1}: ``Truth is an antinomy.\index{antinomy (philosophical notion)}  \ldots What is needed is a formal logical
theory of antinomy," and he goes on explaining a way to obtain such a theory.
 In Letter X, titled \emph{Sophia}, he declares: ``faith
clothes its knowledge of dogmatic truth in a symbolic garment, in figurative
language, which covers the higher truth and depth of contemplation
in consistent contradictions" (\cite[p. 244]{Flo1}).    In Letter II, titled \emph{doubt}, he writes that 
   ``[the
rational mind] desires self-contradiction" (\cite[p. 17]{Flo1}) and that 
``contradiction is everywhere and always"  \cite[p. 23]{Flo1}. We can read in the same letter: ``At the same time we get:
A is A;
A is not A" \cite[p. 29]{Flo1}. 
  In Letter VI, titled \emph{Contradiction}, he writes: 
   ``For rationality,
truth is contradiction, and this contradiction becomes explicit as soon as
truth acquires a verbal formulation. \ldots.  truth is an antinomy,\index{antinomy (philosophical notion)} and it cannot fail to be
so." \cite[p. 109]{Flo1}.

      Florensky's\index{Florensky, Pavel Aleksandrovich} references include Heraclitus,\index{Heraclitus}  Zeno,\index{Zeno} Plato\index{Plato} and Origen\index{Origen}. 
      Regarding Plato, Florensky writes:       ``The majority of his dialogues
are nothing but a gigantic antinomy,\index{antinomy (philosophical notion)}  developed with all care and
artistically dramatized. Plato's\index{Plato} very predilection for the dialogic form of
exposition, i.e., the form of the contraposition of convictions, already
hints at the antinomic nature of his thinking"  \cite[p. 116]{Flo1}.

 Let me make a remark here. It is difficult to be certain of what Florensky\index{Florensky, Pavel Aleksandrovich} had in mind when he referred to Plato. According to Stelios Negrepontis,\footnote{Personal communication.} one might reasonably assume that Florensky had in mind, at least in part,  Plato's most contradictory-looking statements of the dialogue \emph{Parmenides}, especially in the second hypothesis (\emph{Parmenides} 142-155). These statements are indeed scandalously antinomic/contradictory in their appearance. For instance, Plato's intelligible Being  (i) is infinite and finite, (ii) is One and Many, (iii) is in motion and at rest, etc., and even worse: (iv) is older and younger than itself and the others, and is neither older nor younger than itself or the others! However, it would be a rush conclusion that these statements are in fact contradictory. They have a fully rational and mathematically consistent interpretation. According to an interpretation proposed by Negrepontis, 
Plato exploits in these passages the paradoxical, but fully consistent, nature of the infinite but periodic anthyphairesis, and describes it employing what might be considered an abuse of language, to cause surprise and amazement. But Plato is always safe within consistent non-contradictory mathematics. The same remarks apply to Zeno's scandalous and fascinating paradoxes since they are of the same nature as Plato's paradoxes, and in fact form the primary source for Plato's paradoxical arguments. I refer the reader to the interpretation of Zeno's paradoxes by  Negrepontis in \cite{Negrepontis}. 
In the same page (\cite[p. 116]{Flo1}), Florensky declares:  ``Knowledge of contradiction and love of contradiction, along with ancient
skepticism, appear to be the highest achievement of antiquity."  In the same chapter, he gives a detailed list of dogmatic antinomies.      ``The Holy Scripture, he says, is full of antinomies."   The Apostle Paul's epistle to the Romans, he writes,  is an ``antinomy-charged bursting bomb against the
rational mind."
  The paradoxical method of reasoning is justified by symbolic logic    \cite[p. 112]{Flo1}.
 In Chapter XIV of the \emph{Pillar},\index{Florensky, Pavel Aleksandrovich!\emph{The pillar and ground of the truth}} he writes:  ``There are an uncountable number of antinomies. There are as many of
them as there can be acts of the rational mind. But, as we have already
noted, antinomies are essentially reducible to the dilemma: ``finiteness or
infiniteness" \cite[p. 344]{Flo1}. We are back to mathematics.
We already quoted Lebesgue\index{Lebesgue, Henri} saying that there might exist a set which is nameable and which is neither finite nor infinite (Footnote \ref{Leb}.)

 In Letter VI, (\emph{Contradiction}), Florensky makes a list of eleven theological antinomies\footnote{``I deliberately exclude certain antinomies
that I intend to examine in a future book", he writes.} In this list, we find for example:
  \begin{itemize}

 \item   \emph{Retribution}.--- Thesis: Retribution applied to all according to their works (Rom.
2:6–10; 2 Cor 5:10). Antithesis: Free forgiveness of the redeemed
(Rom. 4:4, 9:11, 11:6)."

\item \emph{Grace}.---  Thesis: ``Where sin abounded, grace did much more abound"
(Rom. 5:20). Antithesis: ``Shall we continue in sin, that
grace may abound? God forbid" (Rom. 6:1–2). ``If we say that
we have no sin, we deceive ourselves, and the truth is not in
us" (1 John 1:8).
 
 \item  \emph{Faith}.---  Thesis: Faith is free and depends on the free will of man (John
3:16–28). Antithesis: Faith is God's gift and is not found in human will but in the will of the Father Who draws us to Christ (John 6:44).
 \end{itemize}
 
There are eight other such theological themes in his list, with antinomic approaches.  I refer the interested reader  to the detailed list together with Florensky's\index{Florensky, Pavel Aleksandrovich} commentary in p. 121-123 of the \emph{Pillar}.\index{Florensky, Pavel Aleksandrovich!\emph{The pillar and ground of the truth}}
             
             Antinomy is found in Yahweh, the supreme order of the world. The God of the Old Testament, who punishes severely and looks very cruel, has two entirely antinomic faces.  Florensky writes in Chapter VI of the \emph{Pillar}: 
             ``The Book of Job wholly consists of such a concentrated experience of
contradiction. This book is wholly constructed on the idea of antinomialness."\footnote{I am using the wording of the English translation \cite{Flo1}.} C. F. Jung, in his book \emph{Answer to Job}, reached conclusions similar to those of Florensky.\footnote{I owe this remark to S. Tanabé.}

             What can we conclude? 
             
             Mathematics is the most exact science, and with the actual axioms, there is no proof
               that there is no contradiction in mathematics. If ever we find
             contradictions in mathematics, then there would be contradictions in all the other sciences. In theology, we already live with contradictions. 
        This is expressed in a most dramatic style in \emph{The Pillar and ground of the truth}, which is Florensky's major theological written work. The \emph{Pillar} is a treatise on contradiction and antinomies.
                           
             From the more recent period, let me mention that the French philosopher and epistemologist  Jules Vuillemin (192-2001), who was a professor at  Collège de France, dedicated his course there,  
 during 7 years  (in the 1960s-1970s), to the notion of mathematical antinomy,\footnote{The titles of  Vuillemin's  courses at the Collège de France  include: 
 1962-63: La notion d'antinomie : analyse des doctrines de Kant et Hegel; 1963-64: La notion d'antinomie : l'origine des antinomies dans les mathématiques modernes; 
 1964-65: La notion d'antinomie : les solutions de Russell en 1905; 1965-66: Sur la notion d'antinomie : Découverte des antinomies ``épistémologiques"; 
 1966-1967: Sur la notion d'antinomie (1890-1910): Conséquences des antinomies épistémologiques pour la conception du principe d'induction complète, logicisme, formalisme, intuitionnisme;  
 1967-1968: Sur la notion d'antinomie : Analyse de l'article de Russell ``Mathematical Logic as based on the Theory of Types"; 
 1968-1969: Sur la notion d'antinomie : Analyse de Russell : Préface à la première édition des ``Principia mathematica" (1910), ``On the Relation of Universals and Particulars" (1911-1912).}
  after which he wrote a book on antinomy in the concept of God \cite{Vuillemin}.

             \section{The philosophers}
\begin{figure}
\begin{center} 
\includegraphics[width=10 cm]{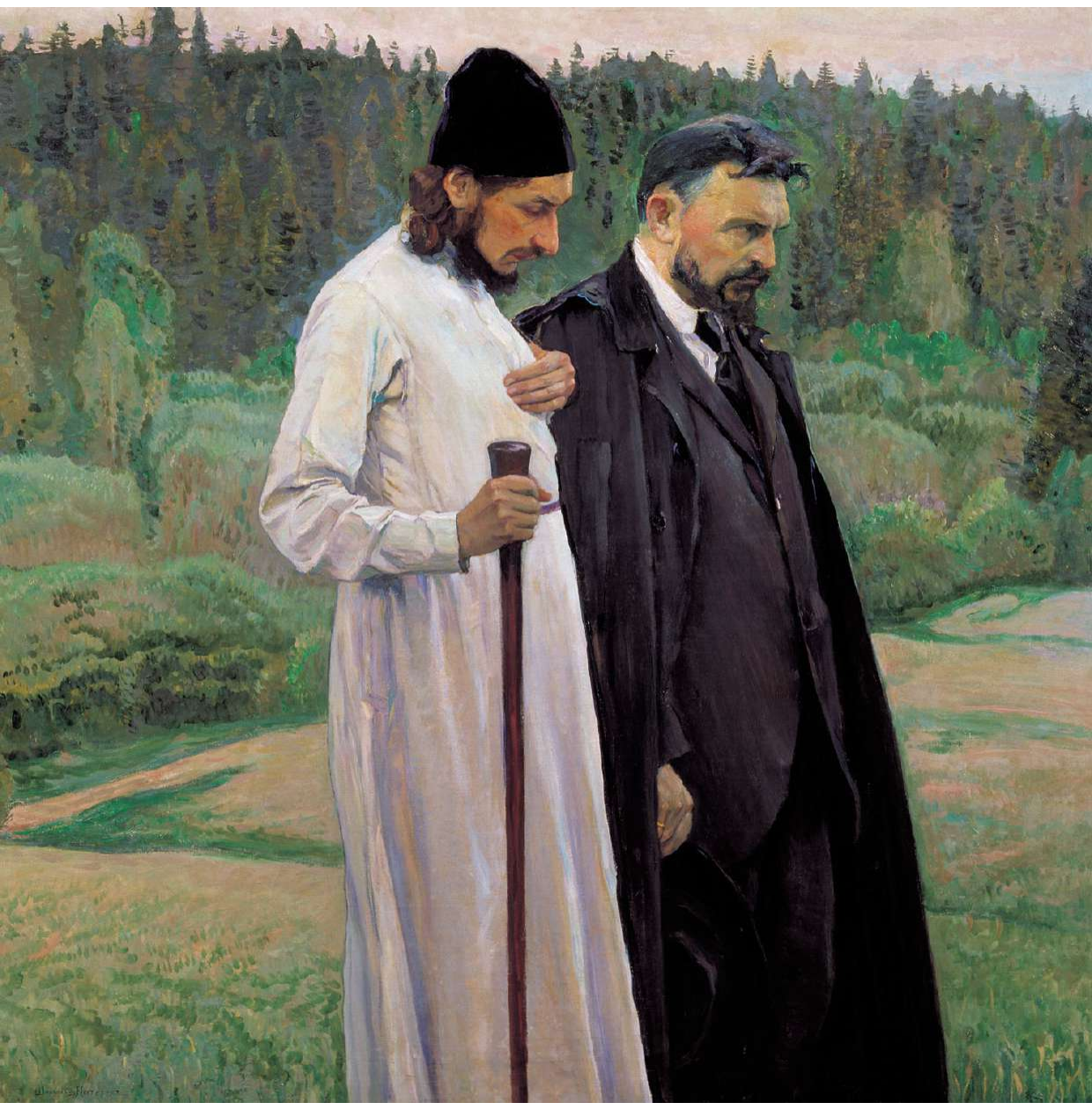} 
\caption{Pavel Florensky and Sergei Bulgakov,  by Mikhail Nesterov, painted in 1917, Tretyakov Gallery, Moscow.}
\label{f:Nesterov}
\end{center}
\end{figure}

             There is a famous canvas, owned by the Moscow Tretyakov Gallery, representing  Pavel Florensky with the philosopher Sergei Bulgakov,\index{Bulgakov, Sergei Nikolaevich} painted in 1917 by Mikhail Nesterov\index{Nesterov, Mikhail} (1862-1942 ), a brilliant representative of the Russian symbolist movement.\index{Russian symbolist movement} The painting is known under the name \emph{The philosophers}.
         It is represented here in Figure \ref{f:Nesterov}. The quiet landscape in the background is the countryside near Sergiev Posad.  The contrast between the two philosophers is  manifest in their face expressions and in the colors of their clothes. Bulgakov became priest in 1918 and was expelled from Russia, together with 160 other intellectuals, in 1922, on one of the \emph{Philosophers' boats}. Florensky decided to stay in Russia.
             Nicoletta Misler wrote a remarkable article on this painting in the Tretyakov gallery Magazine to which I refer the reader \cite{Milser-T}. In this article she quotes a passage by Bulgakov, from a Russian version of \cite{Bulgakov}. I reproduce here an extended passage which I translated from the French version:
             \begin{quote}\small
             A recollection, and with it the foreshadowing of events and achievements to come, never leaves me. It is the portrait of the two of us, painted by our common friend Nesterov on an evening of May 1917 in Father Paul's garden. For the artist, it was not just a portrait of two friends by a third, but a spiritual vision of the times. For the painter, the two faces expressed the same understanding, but for one as a vision of horror, for the other as peace, the joy of overcoming hardship. The painter himself had doubts as to whether he should portray the \emph{first} character, and tried to remake the portrait by replacing the horror with idyll and the tragedy with placidity. But he immediately felt the unbearable falseness of this substitution, and returned to his first vision. On the other hand, he immediately found Father Paul's face, its artistic and spiritual evidence, which he had no need to modify. This was the clear artistic vision of two figures of the Russian apocalypse, of one and the other face of earthly existence, the first in struggle and disarray (and in my heart it was about the fate of my friend), the second in the victorious fulfillment we are now contemplating\ldots
             \end{quote}

       \section{Reverse perspective again} \label{s:reverse}

\begin{figure}
\begin{center} 
\includegraphics[width=6.91 cm]{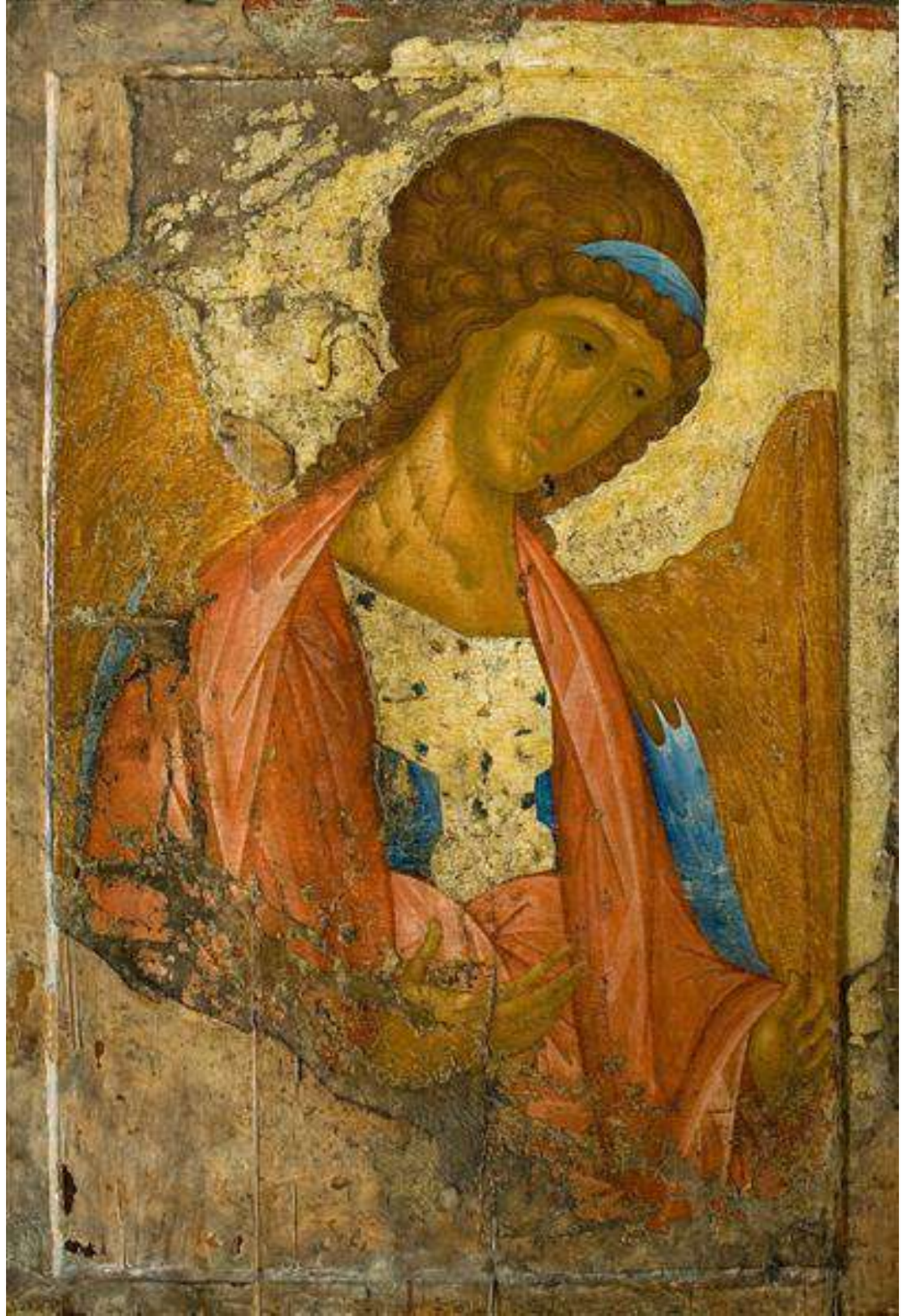} 
\includegraphics[width=7 cm]{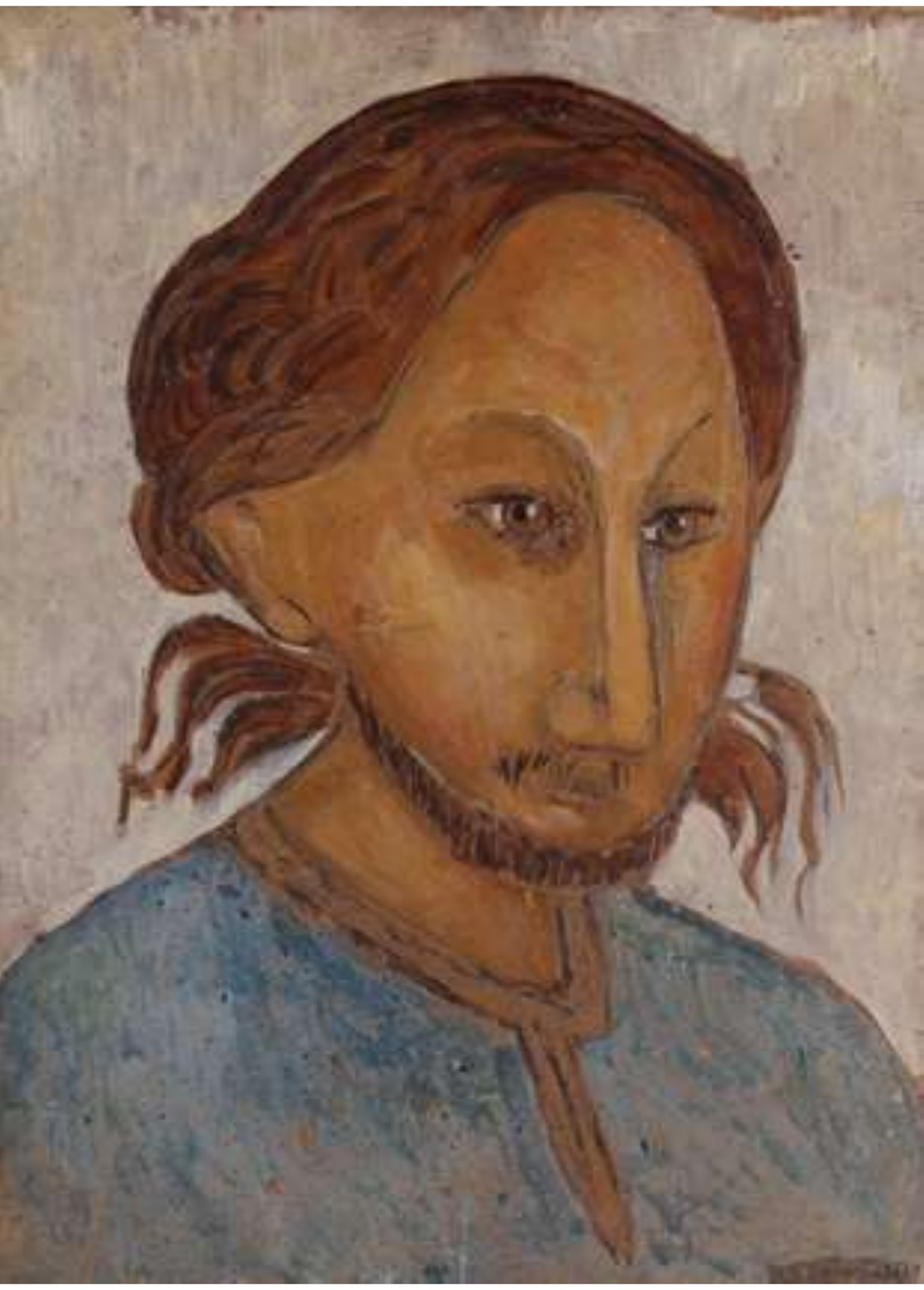} 
\caption{Two paintings drawn in reverse perspective. Left: The Archangel Michael, by Andrei Rublev, ca. 1414, egg tempera. Tretyakov Gallery, Moscow. Right: Portrait  of Pavel Florensky, by Vladimir Alexeyevich Komarovsky. Oil on canvas, 1924. The Pavel Florensky Museum, Moscow.}
\label{f:Michael}
\end{center}
\end{figure}

       Andrei Tarkovsky\index{Tarkovsky, Andrei} was close to Florensky, philosophically, spiritually, and also in his way of making films.    Beatrice Barbalato writes about Tarkovsky: 
        ``Inspired by Florensky's thinking\index{Florensky, Pavel Aleksandrovich!\emph{The reverse perspective}}  (\emph{Reverse Perspective}, 1916), his filmic visual
focuses radiate from the minute detail, both as seen in a dream and as occurs in experiential reality" \cite{Barbalato}.
     Jaromír Blažejovský writes, on the same subject \cite{Bla}: 
     \begin{quote}\small
     Tarkovsky used in \index{Andrei Rublev}\index{Tarkovsky, Andrei!\emph{Andrei Rublev}} what we have
called---following Pavel Florensky---``reverse perspective," typical for Russian Orthodox
icons: narration from various points of view and completed from relatively independent parts
of time and space. The effect of this pattern is a feeling of infinity of the narrative world.
\end{quote} 
The Trinity icon---Abraham's hospitality---which is shown in the film's finale, is an epitome of inverted perspective. I have reproduced in Figure \ref{f:Michael} two paintings, in which the drawing is done in reverse perspective. In such a portrait drawing, the back of the face appears larger than usual, as does the hair.

%
%
%
%
%

              Tarkovsky quotes Florensky at several occasions. In chapter IV of his book \cite{Sculpting}, titled \emph{Cinema's destined role}, he makes some observations on\index{reverse perspective} reverse perspective in Florensky's writings, pointing out relations with his own works:  
\begin{quote}\small
I am reminded of a curious observation of
Father Pavel Florensky's\index{Florensky, Pavel Aleksandrovich} in his book, \emph{The iconostasis}.\index{Florensky, Pavel Aleksandrovich!\emph{The Iconostatis}} He says
that the inverted perspective in the works of that period was not the
result of Russian icon-painters being unaware of the optical laws
which had been assimilated by the Italian Renaissance, after being
developed in Italy by Leon Battista Alberti. Florensky argues,
convincingly, that it was not possible to observe nature without
discovering perspective, it was bound to be noticed. For the time
being, however, it might not be needed---it could be ignored. So the
inverted perspective in ancient Russian painting, the denial of
Renaissance perspective, expresses the need to throw light on certain
spiritual problems which Russian painters, unlike their Italian
counterparts of the Quattrocento, had taken upon themselves. (One
account has it, incidentally, that Andrei Rublev\index{Tarkovsky, Andrei!\emph{Andrei Rublev}}  had actually
visited Venice, in which case he must have been aware of what
Italian painters had been doing with perspective.)
\end{quote}

                Florensky writes in his book on reverse perspective, regarding the drawings of children  (\cite[p. 219 of the English translation]{Flo2}): 
\begin{quote}\small The drawings of children,
in their lack of perspective and especially their use of reverse perspective,
vividly recall mediaeval drawings, despite the efforts of educators to instil
in children the laws of linear perspective. It is only when they lose their spontaneous
relationship to the world that children lose reverse perspective and
submit to the schema with which they have been indoctrinated
\end{quote}

 In my article on Manin\index{Manin, Yuri Ivanovich}  which appears in the present Handbook \cite{Manin}, I mention Tarkovsky\index{Tarkovsky, Andrei} several times. Reverse perspective appears on several occasions in Manin's  writings.\index{Manin, Yuri Ivanovich} It describes a world where there is no pre-determined logic. Manin, like Florensky,  notes that children draw, before they are formatted by the adults' reason, in reverse perspective.  I refer to chapter III.4 of the French version of his book \emph{Mathematics as a metaphor} \cite{Metaphor},\index{Manin, Yuri Ivanovich!\emph{Mathematics as a metaphor}} titled \emph{Spontaneous artistic activity, the origin of logograms and mathematical intuition}.\footnote{The French version of this book is much more complete than the English.}  Let me add that the Swiss biologist, psychologist, logician and epistemologist Jean Piaget\index{Piaget, Jean} (1896-1980) has an article titled \emph{Le développement des perceptions chez l'enfant} (The development of perception in children) \cite{Piaget} whose first section is titled \emph{Structures projectives (perspective)}. Piaget discusses there the evolution of the notion of perspective in child development.  He states in particular that, from the study of children's drawings, it appears that they are unable to reproduce perspective until around the age of 9 or 10, which is the stage when children draw ``what they know about an object, not what they see of it".\footnote{I learned about Piaget's article by S. Tanabé.}

 Finally, let me mention a theory developed by R. K. Luneburg \cite{Luneburg}  according to which  objects  in the
region of binocular vision (which differs substantially from
monocular vision only for nearby objects) are perceived
in terms of Lobachevsky geometry. This is developed by B. V. Rauschenbach in his book \cite{Rauschenbach} and in his paper \cite{Rauschenbach1} titled \emph{On my concept of perceptual perspective that accounts for parallel and inverted
perspective in pictorial art}.\footnote{I learned this from V. Berestovski\u\i . Rauschenbach's theory is commented on and used in Berestovski\u\i 's paper on Lobachevsky geometry \cite{Berestovsky}.}

   \section{In guise of a conclusion} \label{s:conclusion}

       I wrote at the beginning of this essay that Florensky\index{Florensky, Pavel Aleksandrovich} has been compared to Leonardo and Pascal.\index{Pascal, Blaise} The comparison is amply justified. I would like to add to this list two mathematicians of recent times, René Thom\index{Thom, René} and Yuri Manin,\index{Manin, Yuri Ivanovich}  for whom, like Florensky, mathematics, besides being the langage of nature, is that of philosophy, linguistics, poetry and  magic.  
 To conclude,\index{Florensky, Pavel Aleksandrovich} let me quote some verses of a poem by Manin titled \emph{Péniches sur le Rhin (élégie)} in which he mentions Florensky's name \cite[p. 523]{Metaphor}:

 \begin{verse}
Dieu est inconnaissable mais nommable, \\

Le nom de Dieu, c'est Dieu,\\

croyaient les adorateurs du Nom.\\

Le Starets Hilarion, Pavel Florensky ...\\

L'infini non plus n'est pas connaissable, mais il est nommable.\\

Qui connaît les péniches ? Qui leur donne un nom ?\\

On ne sait pas si c'est Lui --- si c'est Elle ...\\

\end{verse}

\bigskip

\noindent {\bf Acknowledgements} I would like to thank V. N.~Berestovski\u\i , Stelios Negrepontis, Ken'ichi Ohshika, Victor Pambuccian, Arkady Plotnitsky,  Valentin Poénaru,  Galina Sinkevich, Alexei Sossinsky and Susumu Tanabé for useful comments and suggestions on previous versions of this chapter. 
This work is supported by the lnterdisciplinary Thematic lnstitute CREAA, as part of the ITl 2021-2028 program of the University of Strasbourg, CNRS and Inserm (ldEx Unistra ANR-10-IDEX-0002), and the French lnvestments for the Future Program.

\end{document}